\documentclass[10pt]{article}
\usepackage{latexsym,amstex}
\usepackage[all]{xy}
\xyoption{v2}
\usepackage{epsfig}
\newtheorem{theorem}{Theorem}[section]
\newtheorem{proposition}[theorem]{Proposition}
\newtheorem{corollary}[theorem]{Corollary}
\newtheorem{lemma}[theorem]{Lemma}

\newtheorem{hypothesis}{IH}

\newenvironment{definition}
{\bigskip\refstepcounter{theorem}\noindent{\bf Definition \thetheorem.}}
{\bigskip}

\newenvironment{remark}
{\bigskip\refstepcounter{theorem}\noindent{\bf Remark \thetheorem.}}
{\bigskip}

\newenvironment{exercise}
{\bigskip\refstepcounter{theorem}\noindent{\bf Exercise \thetheorem.}}
{\bigskip}

\newenvironment{example}
{\bigskip\refstepcounter{theorem}\noindent{\bf Example \thetheorem.}}
{\bigskip}

\newcommand{\MySection}[1]
{\section{ #1}}

\begin{document}
\begin{center}
{{\bf\Large Iterated Monoidal Categories}

\vspace{ 1cm}

{\large C. BALTEANU and Z. FIEDOROWICZ}

\vspace{ .5cm}
{\it Ohio State University\\
Columbus, Ohio 43210-1174, USA}

\vspace{ 1cm}
{\large AND}
\vspace{ 1cm}

{\large R. SCHW\"ANZL and R. VOGT}

\vspace{ .5cm}
{\it Universit\"at Osnabr\"uck\\
49069 Osnabr\"uck, Germany}}
\end{center}

%\vspace{2.5cm}
\newpage

\begin{center}
{ITERATED MONOIDAL CATEGORIES}
\end{center}
\vspace{3cm}

\begin{center}
{Z. Fiedorowicz

Department of Mathematics\\
The Ohio State University\\
100 Mathematics Building\\
231 West 18th Avenue

Columbus, Ohio 43210-1174\\
USA}
\end{center}
%\end{document}
\newpage

\noindent
{\bf\large Introduction}

\vspace{0.5cm}

For many years it has been known that there is a strong connection between coherence theory
of categories and coherence problems in homotopy theory.  Early work of Stasheff \cite{Sta} and
MacLane \cite{Mac} showed that monoidal categories are analogous in a precise way to 1-fold
loop spaces. Later a similar connection was noted between symmetric monoidal categories and
infinite loop spaces. This connection was exploited with great success in
algebraic K-theory.  For instance the group completion of the nerve of a symmetric monoidal
category is an infinite loop space, and the homotopy groups of this infinite loop space are
the Quillen K-groups of that category, which provide algebraic information about the original
category. Conversely this fact has also been used to construct new examples of infinite loop
spaces and infinite loop maps of great interest to topologists.

In recent years many examples of a new kind of algebraic structure on a category
have been discovered: braided monoidal categories, such as categories of representations
of quantum groups (cf. \cite{JS} \& \cite{Kas}). It is striking to note that there appears
to be a very similar connection between braided monoidal categories and 2-fold loop spaces. 
It is shown in \cite{ZF} that the group completion of the nerve of a braided monoidal category
is a 2-fold loop space. This result raises an obvious question: what algebraic structure on
a category corresponds to an $n$-fold loop structure for $3\le n<\infty$? Unfortunately
the proof sheds no light on this matter.

In this paper we provide a comprehensive solution to this problem.  Our solution is
based on pursuing an analogy to the tautology that an $n$-fold loop space is a loop
space in the category of $(n-1)$-fold loop spaces.  Noting the correspondence between
loop spaces and monoidal categories, we iteratively define the notion of $n$-fold
monoidal category as a monoid in the category of $(n-1)$-fold monoidal categories.
There are some subtleties involved in making this definition work: one has to define
``monoidal'' up to a requisite degree of what category theorists call ``laxness''.
If one were to require strict monoidal structures everywhere, then a 2-fold monoidal
category would be strictly commutative and the group completion of its nerve
would be a product of abelian Eilenberg-MacLane spaces.  Another version of this
concept investigated by Joyal and Street \cite{JS} gives a correct analog for 2-fold
loop spaces but for $n\ge 3$ gives a notion equivalent to symmetric monoidal category,
which as noted above is analogous to an infinite loop space.

Our main result is that there is a notion of iterated monoidal category which precisely
corresponds to the notion of an $n$-fold loop space for all $n$.  Firstly the group
completion of the nerve of such a category is an $n$-fold loop space.  Secondly one
can form an operad in the category of small categories which parametrizes the
algebraic structure of an $n$-fold monoidal category.  We show that the nerve of
this categorical operad is a topological operad which is equivalent, as an operad,
to the little $n$-cubes operad, which as shown in \cite{BV1} and \cite{May}
characterizes the notion of $n$-fold loop space.  Thus our result can be regarded
as an algebraic characterization of the notion of $n$-fold loop space.  We also note
that this algebraically defined operad is a finite simplicial operad and is closely
related to the Milgram construction \cite{Mi} for $\Omega^n S^n X$.

We would like to take this opportunity to thank Clemens Berger for several illuminating
email conversations, which provided a key ingredient for the proof of our main result.
We also thank the Deutsche Forschungsgemeinschaft and the Ohio State Mathematical
Research Institute for support during the preparation of this paper.

\pagebreak

\MySection{$n$-fold Monoidal Categories}

In this section we gradually develop the notion of iterated monoidal category.  We start by
recalling the standard notion of monoidal category and defining a slightly nonstandard variant
of the notion of monoidal functor.

\begin{definition}  A {\it (strict) monoidal category} is a category ${\cal C}$
together with a functor 
$\Box: {\cal C}\times{\cal C}\to{\cal C}$  and an object $0$ such that
\begin{enumerate}
\item $\Box$ is strictly associative.
\item $0$ is a strict $2$-sided unit for $\Box$.
\end{enumerate}

A monoidal functor $(F,\eta) :{\cal C}\to{\cal D}$ between monoidal categories consists
of a functor $F$ such that $F(0)=0$ together with a natural transformation
$$
\eta_{A,B}:F(A)\Box F(B)\to F(A\Box B),
$$
which satisfies the following conditions
\begin{enumerate}
\item Internal Associativity: The following diagram commutes
$$
\diagram
F(A)\Box F(B)\Box F(C)
\rrto^{\eta_{A,B}\Box id_{F(C)}}
\dto^{id_{F(A)}\Box \eta_{B,C}}
&&F(A\Box B)\Box F(C)
\dto^{\eta_{A\Box B,C}}\\
F(A)\Box F(B\Box C)
\rrto^{\eta_{A,B\Box C}}
&&F(A\Box B\Box C)
\enddiagram
$$
\item Internal Unit Conditions: $\eta_{A,0}=\eta_{0,A}=id_{F(A)}.$
\end{enumerate}
Given two monoidal functors $(F,\eta) :{\cal C}\to{\cal D}$ and $(G,\zeta) 
:{\cal D}\to{\cal E}$,
we define their composite to be the monoidal functor $(GF,\xi) :
{\cal C}\to{\cal E}$, where
$\xi$ denotes the composite
$$
\diagram
GF(A)\Box GF(B)\rrto^{\zeta_{F(A),F(B)}} 
&& G\bigl(F(A)\Box F(B)\bigr)\rrto^{G(\eta_{A,B})}
&&GF(A\Box B).
\enddiagram
$$
(It is an exercise to check that $\xi$ satisfies the associativity condition above.)
We denote by $\bold{MonCat}$ the category of monoidal categories and monoidal
functors.  Note that the usual product in ${\bf Cat}$ defines a product in
${\bf MonCat}$.
\end{definition}

\begin{remark} It is usually required in standard definitions of the notion of monoidal functor
that $\eta$ be an isomorphism. As we will discuss below, it is crucial for us not to make this
requirement.
\end{remark}

\begin{definition} A $2${\it -fold monoidal category} is a monoid in $\bold{MonCat}$.
This means that we are given a monoidal category $({\cal C},\Box_1,0)$ and a 
monoidal functor
$(\Box_2,\eta):{\cal C}\times{\cal C}\to{\cal C}$ which satisfies
\begin{enumerate}
\item External Associativity: the following diagram commutes in $\bold{MonCat}$

$$
\diagram
{\cal C}\times{\cal C}\times{\cal C}
\rrto^{(\Box_2,\eta)\times id_{\cal C}}
\dto^(0.4){id_{\cal C}\times(\Box_2,\eta)}
&& {\cal C}\times{\cal C}
\dto^{(\Box_2,\eta)}\\
{\cal C}\times{\cal C} 
\rrto^{(\Box_2,\eta)}
&&
{\cal C}
\enddiagram
$$

\item External Unit Conditions: the following diagram commutes in 
$\bold{MonCat}$

$$
\diagram
{\cal C}\times 0 
\rto^{\subseteq}
\dto^{\cong}
& {\cal C}\times{\cal C}
\dto^{(\Box_2,\eta)}
& 0\times{\cal C}
\lto_{\supseteq}
\dto^{\cong}\\
{\cal C} \rto^{=}
&{\cal C} 
&{\cal C} \lto_{=}
\enddiagram
$$
\end{enumerate}

Explicitly this means that we are given a second associative binary operation
$\Box_2:{\cal C}\times{\cal C}\to{\cal C}$, for which $0$ is also a two-sided unit.
Moreover we are given a natural transformation
$$
\eta_{A,B,C,D}: (A\Box_2 B)\Box_1 (C\Box_2 D)\to
(A\Box_1 C)\Box_2(B\Box_1 D).
$$
The internal unit conditions give $\eta_{A,B,0,0}=\eta_{0,0,A,B}=id_{A\Box_2 B}$,
while the external unit conditions give $\eta_{A,0,B,0}=\eta_{0,A,0,B}=id_{A\Box_1 B}$.
The internal associativity condition gives the commutative diagram
$$
\diagram
(U\Box_2 V)\Box_1 (W\Box_2 X)\Box_1 (Y\Box_2 Z)
\xto[rrr]^{\eta_{U,V,W,X}\Box_1 id_{Y\Box_2 Z}}
\dto^{id_{U\Box_2 V}\Box_1 \eta_{W,X,Y,Z}}
&&&\bigl((U\Box_1 W)\Box_2(V\Box_1 X)\bigr)\Box_1 (Y\Box_2 Z)
\dto^{\eta_{U\Box_1 W,V\Box_1 X,Y,Z}}\\
(U\Box_2 V)\Box_1 \bigl((W\Box_1 Y)\Box_2(X\Box_1 Z)\bigr)
\xto[rrr]^{\eta_{U,V,W\Box_1 Y,X\Box_1 Z}}
&&& (U\Box_1 W\Box_1 Y)\Box_2(V\Box_1 X\Box_1 Z)
\enddiagram
$$
The external associativity condition gives the commutative diagram
$$
\diagram
(U\Box_2 V\Box_2 W)\Box_1 (X\Box_2 Y\Box_2 Z)
\xto[rrr]^{\eta_{U\Box_2 V,W,X\Box_2 Y,Z}}
\dto^{\eta_{U,V\Box_2 W,X,Y\Box_2 Z}}
&&& \bigl((U\Box_2 V)\Box_1 (X\Box_2 Y)\bigr)\Box_2(W\Box_1 Z)
\dto^{\eta_{U,V,X,Y}\Box_2 id_{W\Box_1 Z}}\\
(U\Box_1 X)\Box_2\bigl((V\Box_2 W)\Box_1 (Y\Box_2 Z)\bigr)
\xto[rrr]^{id_{U\Box_1 X}\Box_2\eta_{V,W,Y,Z}}
&&& (U\Box_1 X)\Box_2(V\Box_1 Y)\Box_2(W\Box_1 Z)
\enddiagram
$$
\end{definition}

\begin{remark} Notice that we have natural transformations
$$
\eta_{A,0,0,B}:A\Box_1 B\to A\Box_2 B\qquad\mbox{ and }\qquad
\eta_{0,A,B,0}:A\Box_1 B\to B\Box_2 A.
$$
If we had insisted a 2-fold monoidal category be a monoid in the category of monoidal categories
and {\it strictly monoidal\/} functors, this would amount to requiring that $\eta=id$.  In view
of the above, this would imply $A\Box_1 B = A\Box_2 B = B\Box_1 A$ and similarly for morphisms.
Thus the nerve of such a category would be a commutative topological monoid and its group
completion would be equivalent to a product of abelian Eilenberg-MacLane spaces.
\end{remark}

\begin{remark} Recall that a braided monoidal category (also known as braided 
tensor category) is a category ${\cal C}$ together with a functor ${\Box_1 :{\cal C}
\times{\cal C}\to{\cal C}}$ which is strictly associative, has a strict $2$-sided 
unit object $0$ and with a natural commutativity isomorphism 
$c_{A,B}: A\Box_1 B \longrightarrow B\Box_1 A $
 satisfying the following properties:
\begin{enumerate}
\item Unit Condition: $c_{A,0}=c_{0,A}=id_A$.
\item Associativity Conditions: For any three objects $A$, $B$, $C$ the following diagrams
commute:
$$
\diagram
A\Box_1 B\Box_1 C 
\rrto^{c_{A\Box_1 B,C}}
\drto_{id_A\Box_1 c_{B,C}\qquad}
&&C\Box_1 A\Box_1 B\\
& A\Box_1 C\Box_1 B
\urto_{\qquad c_{A,C}\Box_1 id_B}
\enddiagram
$$
$$
\diagram
A\Box_1 B\Box_1 C 
\rrto^{c_{A,B\Box_1 C}}
\drto_{c_{A,B}\Box_1 id_C\qquad}
&&B\Box_1 C\Box_1 A\\
&B\Box_1 A\Box_1 C
\urto_{\qquad id_B\Box_1 c_{A,C}}
\enddiagram
$$
\end{enumerate}
We claim that a braided monoidal category is exactly the same thing as a $2$-fold monoidal category
with $\Box_1=\Box_2$, $\eta$ an isomorphism, and with
$$
\eta_{A,B,0,C}=\eta_{A,0,B,C}=id_{A\Box_1 B\Box_1 C}.
$$
Assuming that $\Box_1 =\Box_2$ and that the natural isomorphism 
$\eta_{A,B,C,D}$ satisfies $\eta_{A,B,0,C}=\eta_{A,0,B,C}=id_{A\Box_1 B\Box_1 C}$, 
one proceeds as follows to show that we have a braided monoidal category.  
In the internal associativity diagram take $V=W=0$ and obtain that 
$\eta_{U,X,Y,Z}=id_U\Box_1 \eta_{0,X,Y,Z}$. Then take $X=Y=0$ and obtain that
$\eta_{U,V,W,Z}=\eta_{U,V,W,0}\Box_1 id_Z$. Combining these two facts, one 
obtains that
$$
\eta_{A,B,C,D}=id_A\Box_1 c_{B,C}\Box_1 id_D,
\eqno(*)
$$
where $c_{B,C}=\eta_{0,B,C,0}$. Then take $U=Z=W=0$ in the internal 
associativity law to get the first associativity law for $c$, and take 
$U=Z=X=0$ to get the other one.  With the additional conditions we have here 
the external associativity law is superfluous.

Conversely given a braided monoidal category, we can define a 2-fold monoidal 
structure by $(*)$.
\end{remark}

\begin{remark}\label{Joyal:Street1} Joyal and Street \cite{JS} considered a very
similar concept to our notion of 2-fold monoidal category.  They loosened our requirement
that the two operations $\Box_1$ and $\Box_2$ be strictly associative with a strict unit by 
only requiring these conditions to hold up to coherent natural isomorphisms.  
More significantly they required the natural transformation $\eta_{A,B,C,D}$ 
to be an isomorphism.  They then showed that such a category is naturally 
equivalent to a braided monoidal category. Briefly given such a category one 
obtains an equivalent braided monoidal category by discarding one of the two 
operations, say $\Box_2$, and defining the commutativity isomorphism for the 
remaining operation $\Box_1$ to be the composite
$$
\diagram
A\Box_1 B\rrto^{\eta_{0,A,B,0}} 
&& B\Box_2 A\rrto^{\eta_{B,0,0,A}^{-1}}
&& B\Box_1 A.
\enddiagram
$$

Our requirement that the operations be strictly associative and unital are not
significant restrictions and were adopted for convenience and simplicity.  One
can always replace categories with operations which are associative and unital up to 
coherent natural isomorphisms by equivalent categories with strictly associative and
unital operations. 
\end{remark}

There is now a pretty obvious way to define the notion of a 2-fold monoidal functor
between 2-fold monoidal categories $F:{\cal C}\to{\cal D}$. It is a functor together
with two natural transformations:
$$\lambda^1_{A,B}:F(A)\Box_1 F(B)\to F(A\Box_1 B)$$
$$\lambda^2_{A,B}:F(A)\Box_2 F(B)\to F(A\Box_2 B)$$
satisfying the same associativity and unit conditions as in the case of monoidal functors.
In addition we require that the following hexagonal interchange diagram commutes:
$$
\diagram
(F(A)\Box_2 F(B))\Box_1(F(C)\Box_2 F(D))
\xto[rrr]^{\scriptstyle {\eta_{F(A),F(B),F(C),F(D)}}}
\dto^{\scriptstyle {\lambda^2_{A,B}\Box_1\lambda^2_{C,D}}}
&&&(F(A)\Box_1 F(C))\Box_2(F(B)\Box_1 F(D))
\dto^{\scriptstyle \lambda^1_{A,C}\Box_2\lambda^1_{B,D}}\\
F(A\Box_2 B)\Box_1 F(C\Box_2 D)
\dto^{\scriptstyle {\lambda^1_{A\Box_2 B,C\Box_2 D}}}
&&&F(A\Box_1 C)\Box_2 F(B\Box_1 D)
\dto^{\scriptstyle {\lambda^2_{A\Box_1 C,B\Box_1 D}}}\\
F((A\Box_2 B)\Box_1(C\Box_2 D))
\xto[rrr]^{\scriptstyle F(\eta_{A,B,C,D})}
&&& F((A\Box_1 C)\Box_2(B\Box_1 D))
\enddiagram
$$

We can now define the category $\bold{2\mbox{\bf -}MonCat}$ of 2-fold monoidal categories and
2-fold monoidal functors, and then define a 3-fold monoidal category as a monoid in
$\bold{2\mbox{\bf -}MonCat}$.  From this point on, the iteration of this notion is quite straightforward
and we arrive at the following definitions.

\begin{definition} An $n${\it -fold monoidal category} is a category ${\cal C}$
with the following structure.
\begin{enumerate}
\item There are $n$ distinct multiplications
$$\Box_1,\Box_2,\dots, \Box_n:{\cal C}\times{\cal C}\to{\cal C}$$
which are strictly associative and ${\cal C}$ has an object $0$ which is a strict unit
for all the multiplications.
\item For each pair $(i,j)$ such that $1\le i<j\le n$ there is a natural
transformation
$$\eta^{ij}_{A,B,C,D}: (A\Box_j B)\Box_i(C\Box_j D)\to
(A\Box_i C)\Box_j(B\Box_i D).$$
\end{enumerate}
These natural transformations $\eta^{ij}$ are subject to the following conditions:
\begin{enumerate}
\item[(a)] Internal unit condition: 
$\eta^{ij}_{A,B,0,0}=\eta^{ij}_{0,0,A,B}=id_{A\Box_j B}$
\item[(b)] External unit condition:
$\eta^{ij}_{A,0,B,0}=\eta^{ij}_{0,A,0,B}=id_{A\Box_i B}$
\item[(c)] Internal associativity condition: The following diagram commutes

$$
\diagram
(U\Box_j V)\Box_i(W\Box_j X)\Box_i(Y\Box_j Z)
\xto[rrr]^{\scriptstyle \eta^{ij}_{U,V,W,X}\Box_i id_{Y\Box_j Z}}
\dto^{\scriptstyle id_{U\Box_j V}\Box_i\eta^{ij}_{W,X,Y,Z}}
&&&\bigl((U\Box_i W)\Box_j(V\Box_i X)\bigr)\Box_i(Y\Box_j Z)
\dto^{\scriptstyle \eta^{ij}_{U\Box_iW,V\Box_i X,YZ}}\\
(U\Box_j V)\Box_i\bigl((W\Box_i Y)\Box_j(X\Box_i Z)\bigr)
\xto[rrr]^{\scriptstyle \eta^{ij}_{U,V,W\Box_i Y,X\Box_i Z}}
&&&(U\Box_i W\Box_i Y)\Box_j(V\Box_i X\Box_i Z)
\enddiagram
$$

\item[(d)] External associativity condition: The following diagram commutes
$$
\diagram
(U\Box_j V\Box_j W)\Box_i(X\Box_j Y\Box_j Z)
\xto[rrr]^{\scriptstyle \eta^{ij}_{U\Box_j V,W,X\Box_j Y,Z}}
\dto^{\scriptstyle \eta^{ij}_{U,V\Box_j W,X,Y\Box_j Z}}
&&& \bigl((U\Box_j V)\Box_i(X\Box_j Y)\bigr)\Box_j(W\Box_i Z)
\dto^{\scriptstyle \eta^{ij}_{U,V,X,Y}\Box_j id_{W\Box_i Z}}\\
(U\Box_i X)\Box_j\bigl((V\Box_j W)\Box_i (Y\Box_j Z)\bigr)
\xto[rrr]^{\scriptstyle id_{U\Box_i X}\Box_j\eta^{ij}_{V,W,Y,Z}}
&&&(U\Box_i X)\Box_j(V\Box_i Y)\Box_j(W\Box_i Z)
\enddiagram
$$
\end{enumerate}
Finally it is required that for each triple $(i,j,k)$ satisfying
$1\le i<j<k\le n$ the (big!) hexagonal interchange diagram commutes.

\setlength{\unitlength}{0.012500in}%
\begingroup\makeatletter\ifx\SetFigFont\undefined
% extract first six characters in \fmtname
\def\x#1#2#3#4#5#6#7\relax{\def\x{#1#2#3#4#5#6}}%
\expandafter\x\fmtname xxxxxx\relax \def\y{splain}%
\ifx\x\y   % LaTeX or SliTeX?
\gdef\SetFigFont#1#2#3{%
  \ifnum #1<17\tiny\else \ifnum #1<20\small\else
  \ifnum #1<24\normalsize\else \ifnum #1<29\large\else
  \ifnum #1<34\Large\else \ifnum #1<41\LARGE\else
     \huge\fi\fi\fi\fi\fi\fi
  \csname #3\endcsname}%
\else
\gdef\SetFigFont#1#2#3{\begingroup
  \count@#1\relax \ifnum 25<\count@\count@25\fi
  \def\x{\endgroup\@setsize\SetFigFont{#2pt}}%
  \expandafter\x
    \csname \romannumeral\the\count@ pt\expandafter\endcsname
    \csname @\romannumeral\the\count@ pt\endcsname
  \csname #3\endcsname}%
\fi
\fi\endgroup
\begin{picture}(520,360)(140,365)
\thinlines
\put(535,700){\vector( 3,-4){ 60}}
\put(200,590){\vector( 0,-1){ 90}}
\put(600,590){\vector( 0,-1){ 90}}
\put(265,700){\vector(-3,-4){ 60}}
\put(200,470){\vector( 3,-4){ 61.800}}
\put(600,470){\vector(-1,-2){ 41}}
\put(140,600){\makebox(0,0)[lb]{\smash{\SetFigFont{10}{12}{rm}
$((A_1 \Box_j B_1) \Box_k (A_2 \Box_j B_2))\Box_i ((C_1 \Box_j D_1) \Box_k
(C_2 \Box_j D_2))$}}}
\put(660,600){\makebox(0,0)[rb]{\smash{\SetFigFont{10}{12}{rm}
$((A_1 \Box_k A_2) \Box_i (C_1 \Box_k C_2)) \Box_j
((B_1 \Box_k B_2) \Box_i (D_1 \Box_k D_2) )$}}}
\put(240,650){\makebox(0,0)[lb]{\smash{\SetFigFont{12}{14.4}{rm}
$\eta^{jk}_{A_1,A_2,B_1,B_2}\Box_i
\eta^{jk}_{C_1,C_2,D_1,D_2}$}}}
\put(560,650){\makebox(0,0)[rb]{\smash{\SetFigFont{12}{14.4}{rm}
$\eta^{ij}_{A_1 \Box_k A_2,B_1 \Box_k B_2,
C_1 \Box_k C_2,D_1 \Box_k D_2}$}}}
\put(400,710){\makebox(0,0)[b]{\smash{\SetFigFont{10}{12}{rm}
$((A_1 \Box_k A_2) \Box_j (B_1 \Box_k B_2)) \Box_i ((C_1 \Box_k C_2) \Box_j
(D_1 \Box_k D_2))$}}}
\put(140,480){\makebox(0,0)[lb]{\smash{\SetFigFont{10}{12}{rm}
$((A_1 \Box_j B_1) \Box_i (C_1 \Box_j D_1)) \Box_k
((A_2 \Box_j B_2) \Box_i (C_2 \Box_j D_2))$}}}
\put(660,480){\makebox(0,0)[rb]{\smash{\SetFigFont{10}{12}{rm}
$((A_1 \Box_i C_1) \Box_k (A_2 \Box_i C_2)) \Box_j ((B_1 \Box_i D_1) \Box_k
(B_2 \Box_i D_2))$}}}
\put(210,540){\makebox(0,0)[lb]{\smash{\SetFigFont{12}{14.4}{rm}
$\eta^{ik}_{A_1 \Box_j B_1,A_2 \Box_j B_2,C_1 \Box_j D_1,
C_2 \Box_j D_2}$}}}
\put(590,540){\makebox(0,0)[rb]{\smash{\SetFigFont{12}{14.4}{rm}
$\eta^{ik}_{A_1,A_2,C_1,C_2}
\Box_j \eta^{ik}_{B_1,B_2,D_1,D_2}$}}}
\put(400,365){\makebox(0,0)[b]{\smash{\SetFigFont{10}{12}{rm}
$((A_1\Box_i C_1) \Box_j (B_1 \Box_i D_1)) \Box_k ((A_2 \Box_i C_2) \Box_j
(B_2 \Box_i D_2))$}}}
\put(240,430){\makebox(0,0)[lb]{\smash{\SetFigFont{12}{14.4}{rm}
$\eta^{ij}_{A_1,B_1,C_1,D_1}
\Box_k \eta^{ij}_{A_2,B_2,C_2,D_2}$}}}
\put(560,430){\makebox(0,0)[rb]{\smash{\SetFigFont{12}{14.4}{rm}
$\eta^{jk}_{A_1 \Box_i C_1,A_2 \Box_j C_2,
B_1 \Box_i D_1,B_2 \Box_j D_2}$}}}
\end{picture}
\end{definition}

\vspace{.5cm}
\begin{definition} An {\it $n$-fold monoidal functor}
$(F,\lambda^1,\dots,\lambda^n):{\cal C}\to{\cal D}$ between $n$-fold monoidal categories
consists of a functor $F$ such that $F(0)=0$ together with natural
transformations
$$\lambda^i_{A,B}:F(A)\Box_i F(B)\to F(A\Box_i B)\quad i=1,2,\dots, n$$
satisfying the same associativity and unit conditions as monoidal functors.
In addition the following hexagonal interchange diagram commutes:
$$
\diagram
(F(A)\Box_j F(B))\Box_i(F(C)\Box_j F(D))
\xto[rrr]^{\scriptstyle {\eta^{ij}_{F(A),F(B),F(C),F(D)}}}
\dto^{\scriptstyle {\lambda^j_{A,B}\Box_i\lambda^j_{C,D}}}
&&&(F(A)\Box_i F(C))\Box_j(F(B)\Box_i F(D))
\dto^{\scriptstyle \lambda^i_{A,C}\Box_j\lambda^i_{B,D}}\\
F(A\Box_j B)\Box_i F(C\Box_j D)
\dto^{\scriptstyle {\lambda^i_{A\Box_j B,C\Box_j D}}}
&&&F(A\Box_i C)\Box_j F(B\Box_i D)
\dto^{\scriptstyle {\lambda^j_{A\Box_i C,B\Box_i D}}}\\
F((A\Box_j B)\Box_i(C\Box_j D))
\xto[rrr]^{\scriptstyle F(\eta^{ij}_{A,B,C,D})}
&&& F((A\Box_i C)\Box_j(B\Box_i D))
\enddiagram
$$
Composition of $n$-fold monoidal functors is defined in exactly the same way
as for monoidal functors.  However there is an additional exercise to check
that the resulting composite satisfies the hexagonal interchange diagram.
\end{definition}

It is pretty straightforward to check that an $(n+1)$-fold monoidal category is
exactly the same thing as a monoid in $\bold{n\mbox{\bf -}MonCat}$, the category of
$n$-fold monoidal categories and functors.  Note that the hexagonal interchange diagrams
for the $(n+1)$-st monoidal operation regarded as an $n$-fold monoidal functor is what
gives rise to the giant hexagonal diagrams involving $\Box_i$, $\Box_j$ and $\Box_{n+1}$.

\begin{remark}\label{symm:remark}
Recall that a symmetric monoidal category is defined in the same way as a
braided monoidal category, subject to the additional requirement that the commutativity
isomorphism
$$c_{A,B}:A\Box B\stackrel{\cong}{\longrightarrow}B\Box A$$
satisfy the symmetry condition
$$c_{B,A}=c_{A,B}^{-1}$$
It is easy to see a symmetric monoidal category is $n$-fold monoidal for all $n$.  One merely
has to take
$$\Box_1=\Box_2=\dots=\Box_n=\Box$$
and define
$$\eta^{ij}_{A,B,C,D}=id_A\Box c_{B,C}\Box id_D$$
for all $i<j$.
\end{remark}

\begin{remark} Joyal and Street \cite{JS} arrived at pretty much the same definitions as we
do in their context.  Because of their insistence that the interchange natural transformations
$\eta^{ij}_{A,B,C,D}$ be isomorphisms, however as they observed, for $n\ge3$ such a notion
is equivalent to the notion of symmetric monoidal category, by an argument similar to
that of Remark \ref{Joyal:Street1}.  Thus the nerves of such categories
have group completions which are infinite loop spaces rather than $n$-fold loop spaces.  In
Remark \ref{inversion:kills} we will give a homotopy theoretic interpretation of this phenomenon.
\end{remark}

\newpage
\MySection{Connection with $n$-fold Loop Spaces}

\vspace{0.5cm}

In this section we sketch a proof of our assertion that the group completion of the nerve of an
$n$-fold monoidal category is a $n$-fold loop space.  The proof closely mimics Thomason's \cite{Th}
proof for the analogous connection between symmetric monoidal categories and infinite loop spaces.
That proof in turn is based on Segal's ideas \cite{Se}.  Our proof sketch omits some important details which
depend on the coherence theorem  for $n$-fold monoidal categories which we will  discuss in
Section 4.  Later on in section Section 6 we will give an alternative
proof of our assertion based on the operad approach to $n$-fold loop spaces due to May \cite{May}.

Segal showed that a space $Y$ is homotopy equivalent to a 1-fold loop space if and only if
one can construct a ``bar construction on $Y$ up to homotopy.''  This means a simplicial space
$X_*:\Delta^{op}\to\bold{Top}$ (where $\bold{Top}$ is the category of compactly
generated spaces) with $X_1=Y$ and satisfying
\begin{enumerate}
\item There is a homotopy equivalence 
$X_n\stackrel{\simeq}{\longrightarrow}(X_1)^n$
induced by certain iterated face maps and $X_0$ is contractible.
\item The multiplication induced by
$(X_1)^2\stackrel{\simeq}{\longleftarrow}X_2\stackrel{d_1}{\longrightarrow}X_1$
admits a homotopy inverse. (This holds if $\pi_0(X_1)$ is a group and if $X_1$ is
numerably contractible, eg. a CW-complex.)
\end{enumerate}

Moreover he showed that the geometric realization $|X_*|$ is an up-to-homotopy delooping of
$X_1=Y$, ie. $\Omega|X_*|\simeq X_1=Y$.  It was subsequently shown \cite{McS} that if condition (2)
is omitted, then under some mild additional homotopy commutativity assumption $H_*(\Omega|X_*|)$
is obtained from $H_*(X_1)$ by inverting the elements of $\pi_0(X_1)\subset H_0(X_1)$.  This
relation is usually referred to as saying that $\Omega|X_*|$ is the {\it group completion\/}
of $X_1$.  Simplicial spaces satisfying condition (1) are referred to as
{\it special $\Delta$-spaces\/}.

Segal also noted that one could formulate categorical versions of these concepts.  For
instance a {\it special $\Delta$-category\/} is a simplicial category
${\cal C}_*:\Delta^{op}\to\bold{Cat}$ satisfying
\begin{enumerate}
\item There is an equivalence of categories 
${\cal C}_n\stackrel{\simeq}{\longrightarrow}({\cal C}_1)^n$
induced by certain iterated face maps and $C_0$ has a initial/terminal object.
\end{enumerate}
Since the nerve construction preserves products and sends categorical equivalences to
homotopy equivalences, the nerve of a special $\Delta$-category is a special $\Delta$-space.

Segal noted that a strictly monoidal category ${\cal C}$ naturally gives rise to a
special $\Delta$-category ${\cal C}_*$ with ${\cal C}_n=({\cal C})^n$ via the bar construction.  If the
monoidal structure is not strictly associative, then one can still construct a
special $\Delta$-category ${\cal C}_*$ but with ${\cal C}_n\simeq({\cal C})^n$. Here 
one has to use the extra flexibility of allowing categorical equivalences 
rather than isomorphisms. 
(This is not critical in this case, since 
monoidal categories are equivalent to strictly associative ones. When one 
attempts to put symmetric monoidal categories in this framework one encounters 
the problem that commutativity can not be made strict.)

Segal's construction of special $\Delta$-categories in the absence of strict
algebraic relations (like associativity) was incomplete and ad hoc.  This was
remedied by Thomason \cite{Th}, who noted that this construction could be done in two steps.
First one can construct a {\it lax functor\/} ${\cal C}_*:\Delta^{op}\to\bold{Cat}$ such that
${\cal C}_n=(C_1)^n$.  Next one could use the result of Street \cite{Str}, which states that
for any category ${\cal I}$ and any lax functor $F:{\cal I}\to\bold{Cat}$, one can construct
an equivalent strict functor $\widehat F:{\cal I}\to\bold{Cat}$.  This functor $\widehat F$
is called the {\it Street rectification\/} of the original lax functor $F$. Applying this
to the lax functor ${\cal C}_*:\Delta^{op}\to\bold{Cat}$, one obtains a strict functor
$\widehat{{\cal C}}_*:\Delta^{op}\to\bold{Cat}$, which is the desired special $\Delta$-category.

While Segal never explicitly considered $n$-fold loop spaces except in the special cases $n=1$ and
$n=\infty$, as noted by Dunn \cite{D}, his ideas can easily be adapted to this case.  One needs to
consider special $(\Delta)^n$-spaces. These are the same thing
as $n$-simplicial spaces
$X_{**\dots*}:\Delta^{op}\times\Delta^{op}\times\dots\times\Delta^{op}\to\bold{Top}$ satisfying the condition
\begin{enumerate}
\item There is a homotopy equivalence 
$X_{p_1,p_2,\dots,p_n}\stackrel{\simeq}{\longrightarrow}(X_{11\dots1})^{p_1p_2\dots p_n}$
induced by certain iterated face maps.
\end{enumerate}
We call such functors {\it special $(\Delta)^n$-spaces\/}. From Segal's results
in the 1-fold loop case, we easily see that for a special  $(\Delta)^n$-space
$X_{**\dots*}$ that $\Omega^n|X_{**\dots*}|$ is a group completion of $X_{**\dots*}$. The notion of special
 $(\Delta)^n$-category can be formulated similarly.

\begin{theorem} An $n$-fold monoidal category ${\cal C}$ determines a lax functor
${\cal C}_{**\dots*}:\Delta^{op}\times\Delta^{op}\times\dots\times\Delta^{op}\to\bold{Cat}$
such that ${\cal C}_{p_1,p_2,\dots,p_n}={\cal C}^{p_1p_2\dots p_n}$.
\end{theorem}

\vspace{0.5cm}
\noindent
{\bf Proof:} The lax functor ${\cal C}_{**\dots*}$ is already specified on objects of
 $(\Delta^{op})^n$.  We begin to define the lax functor on morphisms of $(\Delta^{op})^n$ by
 first considering morphisms of the special form
 $$(id,\dots,id,\alpha,id,\dots,id): (p_1,\dots,p_{i-1},q_i,p_{i+1},\dots,p_n)\longrightarrow
 (p_1,\dots,p_{i-1},p_i,p_{i+1},\dots,p_n)$$
 which have only one nontrivial component $\alpha: q_i\to p_i$ in $\Delta$.

 Recall that given a morphism $\alpha: p_i\to q_i$ in $\Delta^{op}$ and a strict monoidal category
 ${\cal A}$, the bar construction defines a corresponding functor ${\cal A}^{p_i}\to{\cal A}^{q_i}$.
 Now consider the category ${\cal A} = {\cal C}^{p_{i+1}p_{i+2}\dots p_n}$ as a monoidal category
with respect to the $i$-th operation $\Box_i$ applied componentwise.  This defines a functor
$$\diagram
{\cal C}^{p_ip_{i+1}\dots p_n} = {\cal A}^{p_i}\rrto^{\alpha^*}
&&{\cal A}^{q_i}={\cal C}^{q_ip_{i+1}\dots p_n}
\enddiagram$$
Now taking the $p_1p_2\dots p_{i-1}$-fold product of this functor with itself gives a functor
$$\diagram
{\cal C}_{p_1,p_2,\dots,p_n}={\cal C}^{p_1p_2\dots p_n}
\xto[rrr]^{(\alpha^*)^{p_1p_2\dots p_{i-1}}}
&&&{\cal C}^{p_1\dots q_i\dots p_n} = {\cal C}_{p_1,\dots, q_i,\dots, p_n}
\enddiagram$$
which we define to be the value of the lax functor ${\cal C}_{**\dots*}$ on the morphism
$(id,\dots,id,\alpha,id,\dots,id)$.

Now given an arbitrary morphism $(\alpha_1,\alpha_2,\dots,\alpha_n):(p_1,p_2,\dots,p_n)
\longrightarrow(q_1,q_2,\dots,q_n)$ in $(\Delta^{op})^n$ we define its value under
${\cal C}_{**\dots*}$ to be the resulting composite functor.
$${\cal C}_{p_1,p_2,\dots,p_n}\longrightarrow{\cal C}_{q_1,p_2,\dots,p_n}\longrightarrow
\dots\longrightarrow{\cal C}_{q_1,q_2,\dots,q_{n-1},p_n}\longrightarrow
{\cal C}_{q_1,q_2,\dots,q_{n-1},q_n}$$

We claim that this defines a lax functor
${\cal C}_{**\dots*}:\Delta^{op}\times\Delta^{op}\times\dots\times\Delta^{op}\to\bold{Cat}$.
To see this suppose we are given two morphisms
$$\underline{\alpha} = (\alpha_1,\alpha_2,\dots,\alpha_n)$$
$$\underline{\beta} = (\beta_1,\beta_2,\dots,\beta_n)$$
in $(\Delta^{op})^n$ and consider the composite
$$\underline{\beta}\underline{\alpha} = (\beta_1\alpha_1,\beta_2\alpha_2,\dots,\beta_n\alpha_n)$$

By definition the value of ${\cal C}_{**\dots*}$ on $\underline{\beta}\underline{\alpha}$ is the functor
given by the composite of the induced functor
$$(\beta_1\alpha_1,id,\dots,id)^*$$
followed by
$$(id,\beta_2\alpha_2,id,\dots,id)^*$$
etc. Since ${\cal C}_{**\dots*}$ is a functor when restricted to morphisms having the form that all but
one component is trivial, $(\underline{\beta}\underline{\alpha})^*$ can be further decomposed as the
composite
$$(\beta_1,,id,\dots,id)^* (\alpha_1,,id,\dots,id)^*$$
followed by
$$(id,\beta_2,id,\dots,id)^*(id,\alpha_2,id,\dots,id)^*$$
etc.  Similarly $\underline{\beta}^*\underline{\alpha}^*$
breaks up as a composite of exactly the same functors, but composed in a different order.

Thus to construct a natural transformation
$\underline{\beta}^*\underline{\alpha}^*\to(\underline{\beta}\underline{\alpha})^*$
it suffices to construct natural transformations
$$(id,\dots,id,\kappa_i,id,\dots,id)^*(id,\dots,id,\lambda_j,id,\dots,id)^*\longrightarrow
(id,\dots,id,\lambda_j,id,\dots,id)^*(id,\dots,id,\kappa_i,id,\dots,id)^*$$
where the indices $i,j$ indicating the location of the nontrivial components satisfy $i<j$ and
$\kappa_i:p_i\to q_i$, $\lambda_j:p_j\to q_j$ are arbitrary morphisms in $\Delta^{op}$.

Thus we have to construct a natural transformation from the top-right to the left-bottom of the
following diagram:
$$\diagram
\Bigl(({\cal C}^{p_{j+1}\dots p_n})^{p_j}\Bigr)^{p_1p_2\dots p_{j-1}}
\xto[rrr]^{(\lambda_j^*)^{p_1p_2\dots p_{j-1}}}
\ddouble
&&&\Bigl(({\cal C}^{p_{j+1}\dots p_n})^{q_j}\Bigr)^{p_1p_2\dots p_{j-1}}
\ddouble\\
\Bigl(({\cal C}^{p_{i+1}\dots p_n})^{p_i}\Bigr)^{p_1p_2\dots p_{i-1}}
\dto^{(\kappa_i^*)^{p_1p_2\dots p_{i-1}}}
&&&\Bigl(({\cal C}^{p_{i+1}\dots q_j\dots p_n})^{p_i}\Bigr)^{p_1p_2\dots p_{i-1}}
\dto^{(\kappa_i^*)^{p_1p_2\dots p_{i-1}}}\\
\Bigl(({\cal C}^{p_{i+1}\dots p_n})^{q_i}\Bigr)^{p_1p_2\dots p_{i-1}}
\ddouble
&&&\Bigl(({\cal C}^{p_{i+1}\dots q_j\dots p_n})^{q_i}\Bigr)^{p_1p_2\dots p_{i-1}}
\ddouble\\
\Bigl(({\cal C}^{p_{j+1}\dots p_n})^{p_j}\Bigr)^{p_1\dots q_j\dots p_{j-1}}
\xto[rrr]^{(\lambda_j^*)^{p_1\dots q_j\dots p_{j-1}}}
&&&\Bigl(({\cal C}^{p_{j+1}\dots p_n})^{q_j}\Bigr)^{p_1\dots q_j\dots p_{j-1}}
\enddiagram$$
(The vertical equality signs in the diagram are actually canonical permutations.)

Now to construct a natural transformation between two functors taking values in a product
of categories, it suffices to construct a natural transformation separately in each component
of the product.  Thus we may as well assume that
$$p_1=\dots=p_{i-1}=q_i=p_{i+1}=\dots=p_{j-1}=q_j=p_{j+1}=\dots=p_n=1$$
To simplify the notation a little bit, we denote $p_i=r$ and $p_j=s$.  Then the previous
diagram simplifies to
$$\diagram
({\cal C}^r)^s\rto^{\cong}\dto^{(\kappa_i^*)^s}&({\cal C}^s)^r\rrto^{(\lambda_j^*)^r}&&{\cal C}^r\dto^{\kappa_i^*}\\
{\cal C}^s\xto[rrr]^{\lambda_j^*}&&&{\cal C}
\enddiagram$$

Now $\kappa_i^*$ and $\lambda_j^*$, being induced maps in the bar construction, have to
have the general form:
$$\kappa_i^*(A_1,\dots, A_r) = A_k\Box_i A_{k+1}\Box_i\dots\Box_i A_{k+u}
:= {\textstyle\prod_{k\le x\le k+u}^{\Box_i}A_x}$$
$$\lambda_j^*(B_1,\dots, B_s) = B_l\Box_j B_{l+1}\Box_j\dots\Box_j B_{l+v}
:= {\textstyle\prod_{l\le y\le l+v}^{\Box_j}B_y}$$

Now if we track an arbitrary object across the top and right of this diagram we obtain
$$\xymatrix{
\left(C_{xy}\right)_{\begin{array}{c}{\scriptstyle 1\le x\le r}
\\{\scriptstyle 1\le y\le s}\end{array}}^{\begin{array}{c}{\scriptstyle\quad}\\{\scriptstyle\quad}\end{array}}
\ar@{|->}[rr]
&&\left({\textstyle\prod_{l\le y\le l+v}^{\Box_j}C_{xy}}\right)_{1\le x\le r}\ar@{|->}[d]\\
&&{\textstyle\prod_{k\le x\le k+u}^{\Box_i}\left(\prod_{l\le y\le l+v}^{\Box_j}C_{xy}\right)}
}$$
On the other hand if we track it across the left and bottom we obtain:
$$\xymatrix{
\left(C_{xy}\right)_{\begin{array}{c}{\scriptstyle 1\le x\le r}\\{\scriptstyle 1\le y\le s}\end{array}}
\ar@{|->}[d]\\
\left({\textstyle\prod_{k\le x\le k+u}^{\Box_i}C_{xy}}\right)_{1\le y\le s}^{\quad}\ar@{|->}[rr]
&&{\textstyle\prod_{l\le y\le l+v}^{\Box_j}\Bigl(\prod_{k\le x\le k+u}^{\Box_i}C_{xy}\Bigr)}
}$$

It is clear that there is a natural transformation, built out of repeated applications of $\eta^{ij}$,
$${\textstyle\prod_{k\le x\le k+u}^{\Box_i}\Bigl(\prod_{l\le y\le l+v}^{\Box_j}C_{xy}\Bigr)}
\longrightarrow
{\textstyle\prod_{l\le y\le l+v}^{\Box_j}\Bigl(\prod_{k\le x\le k+u}^{\Box_i}C_{xy}\Bigr)}$$

Lastly we must verify that the natural transformations
$\underline{\beta}^*\underline{\alpha}^*\to(\underline{\beta}\underline{\alpha})^*$
we have just constructed satisfy a certain associativity condition.  To see this we rely on the
coherence theorem for $n$-fold monoidal categories, which we will state and prove in the following
two sections.  That theorem states that any diagram built out of the natural transformations
$\eta^{ij}$ must commute.

\begin{theorem} The group completion of the nerve of an $n$-fold monoidal
category is an $n$-fold loop space.
\end{theorem}

\noindent
{\bf Proof Sketch:} By the preceding theorem, we have a lax functor
$${\cal C}_{**\dots*}:\Delta^{op}\times\Delta^{op}\times\dots\times\Delta^{op}\to\bold{Cat}$$
such that ${\cal C}_{p_1,p_2,\dots,p_n}={\cal C}^{p_1p_2\dots p_n}$.
Now apply Street rectification to obtain a genuine
functor
$$\widehat{{\cal C}}_{**\dots*}:\Delta^{op}\times\Delta^{op}\times\dots\times\Delta^{op}
\to\bold{Cat}, $$
with
$\widehat{{\cal C}}_{p_1,p_2,\dots,p_n}\simeq{\cal C}^{p_1p_2\dots p_n}$.  Taking nerves we obtain
a functor $B\widehat{{\cal C}}_{**\dots*}:\Delta^{op}\times\Delta^{op}\times\dots\times\Delta^{op}
\to\bold{Top}$, with 
$B\widehat{{\cal C}}_{p_1,p_2,\dots,p_n}\simeq(B {\cal C})^{p_1,p_2,\dots,p_n}$.
Thus $B\widehat{{\cal C}}_{**\dots*}$ is a
special $(\Delta)^n$-space, and the result follows.

\newpage

\MySection{Free $n$-fold Monoidal Categories and Their Associated Operad: Statement of Main Results}

\vspace{0.5cm}

In this section we consider an alternative and more precise way of relating $n$-fold monoidal categories to $n$-fold
loop spaces: via operads.  First of all we consider free $n$-fold monoidal categories and
construct an associated operad which acts on nerves of $n$-fold monoidal
categories.  We then discuss the relation of this operad to Milgram's permutohedral construction
used to approximate free loop spaces, and to the little $n$-cubes operad of Boardman and Vogt.

\begin{definition} Let ${\cal C}$ be a small category. By ${\cal F}_n{\cal C}$
we will denote the free $n$-fold monoidal category generated by ${\cal C}$.
${\cal F}_n{\cal C}$ may be constructed as follows.  As objects one takes all
finite expressions generated by the objects of ${\cal C}$ using associative
operations $\Box_1$, $\Box_2$,\dots $\Box_n$. For example
$$(((C_1\Box_1 C_2\Box_1 C_3)\Box_2 C_4\Box_2(C_5\Box_3 C_6))\Box_2 C_7)\Box_3
(C_8\Box_2 C_9)$$
Included among such possible expressions is the vacuous expression, denoted
0, which serves as the unit object.  The morphisms of ${\cal F}_n{\cal C}$ are
finite composites of all possible finite formal expressions  generated by
the morphisms of ${\cal C}$ and symbols $\eta^{ij}_{A,B,C,D}$ with $1\le i<j\le n$ and $A$, $B$, $C$,
$D$ objects of ${\cal F}_n{\cal C}$, using the associative operations $\Box_1$, $\Box_2$, \dots,
and $\Box_n$.  Two such composites of formal expressions are identified
if and only if one can be converted to the other by repeated use of various
functoriality, naturality and associativity diagrams. (This is a special case
of forming a colimit in theories, cf. \cite[p. 33 Prop.2.5]{BV1}.)

As a special case we may take ${\cal C}$ to be a finite set whose elements
are taken to be the objects, with the morphisms understood to be just the
identities of these objects.  We will denote by ${\cal M}_n(k)$ the full
subcategory of ${\cal F}_n\{1, 2, \dots,k\}$ whose objects are expressions
in which each element 1, 2, \dots, $k$ occurs exactly once.  For example
$(2\Box_1 1)\Box_2 3$ is an object of ${\cal M}_n(3)$ but not of ${\cal M}_n(4)$,
whereas $(1\Box_2 2)\Box_1 1$ is not in any ${\cal M}_n(k)$. The symmetric
group $\Sigma_k$ acts freely on ${\cal M}_n(k)$ via functors, by permuting
labels on both objects and morphisms. It is easy to see that for any category
${\cal C}$
$${\cal F}_n{\cal C}\cong \coprod_{k\ge0}{\cal M}_n(k)\times_{\Sigma_k}{\cal C}^k$$
In particular
$${\cal F}_n\{1\}\cong \coprod_{k\ge0}{\cal M}_n(k)/\Sigma_k$$

If ${\cal C}$ is already $n$-fold monoidal, then we have a natural evaluation
functor ${\cal F}_n{\cal C}\to{\cal C}$ which gives rise to functors
$${\cal M}_n(k)\times_{\Sigma_k}{\cal C}^n\longrightarrow{\cal C}$$
As a special case we get maps
$${\cal M}_n(k)\times{\cal M}_n(i_1)\times{\cal M}_n(i_2)\times\dots
\times{\cal M}_n(i_k)\longrightarrow{\cal M}_n(i_1+i_2+\dots+i_k)$$
by replacing the labels $\{1,2,\dots,i_j\}$ in ${\cal M}_n(i_j)$
with the labels $\{i_1+i_2+\dots+i_{j-1}+1,\dots, i_1+i_2+\dots+i_{j-1}+i_j\}$.
This gives $\{{\cal M}_n(k)\}_{k\ge0}$ the structure of an operad in
the category of small categories, with a natural action on $n$-fold
monoidal categories. Since the nerve construction preserves products,
the nerve of this categorical operad is a topological operad, which
we also abusively denote ${\cal M}_n$, and this operad acts on nerves of
$n$-fold monoidal categories.
\end{definition}

\begin{definition}\label{set:asLabel} It will be convenient to be a bit more general and consider categories
${\cal M}_n(S)$, where $S$ is an arbitrary finite set. Again we define ${\cal M}_n(S)$ to
be the full subcategory of the free $n$-fold monoidal category ${\cal F}_n(S)$ whose objects
are expressions in which each element of $S$ occurs precisely once.  Obviously any bijection
$S\cong S'$ extends to an isomorphism of categories ${\cal M}_n(S)\cong{\cal M}_n(S')$.  If
$S\subset T$, there is a restriction functor ${\cal M}_n(T)\to{\cal M}_n(S)$, induced by the
functor ${\cal F}_n(T)\to{\cal F}_n(S)$ which sends the elements of $T-S$ to 0.
\end{definition}

The following is an amusing exercise for the reader:

\begin{exercise}Let $a^n_k$ denote the number of objects in ${\cal M}_n(k)/\Sigma_k$.
Then $a^n_0 = a^n_1 = 1$, $a^n_2=n$, $a^n_3=2n^2-n$, $a^n_4=5n^3-5n^2+n$ and we
have the recurrence relation
$$a^n_k = na^n_1a^n_{n-1}+\sum_{i=2}^{n-1}(n-1)a^n_ia^n_{k-i}$$
The ratios $\frac{a^n_{k+1}}{a^n_k}$ slowly increase to a limit of $2n-1+2\sqrt{n^2-n}$.
Thus the number of objects in ${\cal M}_n(k)$ is $k!a^n_k$.
\end{exercise}

While it may seem from the definition that the operads ${\cal M}_n(k)$ are some
infinite-dimensional abstract algebraic monstrosities, this is not the case.
They are actually nice compact polyhedra.

\begin{example}\label{octahedron:example}
It is not difficult to see that ${\cal M}_n(2)$ is the $(n-1)$-dimensional
octahedron and thus homeomorphic to $S^{n-1}$.  (By $(n-1)$-dimensional octahedron we
mean the boundary of the convex hull of $\{\pm\mbox{\boldmath $e_1$}, \pm\mbox{\boldmath $e_2$},\dots,
\pm\mbox{\boldmath $e_n$}\}$, where $\mbox{\boldmath $e_1$}, \mbox{\boldmath $e_2$},\dots,
\mbox{\boldmath $e_n$}$
is the standard basis for $\Bbb R^n$.)  Clearly ${\cal M}_n(2)$ is generated
by the morphisms $\eta^{ij}_{a,0,0,b}$ and $\eta^{ij}_{0,a,b,0}$,
where $a\ne b \in\{1,2\}$ and $1\le i<j\le n$.  The ``Giant Hexagon''
shows that this set of morphisms is closed under composition.  We thus obtain the following
picture for the nerve ${\cal M}_3(2)$
$$
\diagram
&&&{1\Box_3 2}\\
\\
\\
&&{2\Box_2 1}\xto[uuur]^(.3){\eta^{23}_{0210}}
\xto[dddddl]|>>>>>>>>>>>>>>>>>>>>>>>>>>>>>>>>>>>\hole_{\eta^{23}_{2001}}\\
{2\Box_1 1}\urrto_{\eta^{12}_{2001}}\drrto^{\eta^{12}_{0210}}\xto[uuuurrr]^{\eta^{13}_{0210}}
\xto[ddddr]_{\eta^{13}_{2001}}
&&&& {1\Box_1 2}\ullto|>>>>\hole^{\eta^{12}_{0120}}\dllto^{\eta^{12}_{1002}}
\xto[uuuul]_{\eta^{13}_{1002}}\xto[ddddlll]^{\eta^{13}_{0120}}\\
&&{1\Box_2 2}\xto[uuuuur]_{\eta^{23}_{1002}}\xto[dddl]^(.3){\eta^{23}_{0120}}\\
\\
\\
&{2\Box_3 1}
\enddiagram
$$
and this picture obviously generalizes to ${\cal M}_n(2)$ for all $n$.
\end{example}

If we hope to get similar nice pictures of ${\cal M}_n(k)$ for $k>2$, we need a better
description of the categories ${\cal M}_n(k)$ than that given in  the definition.
It is a priori very difficult to determine when two different formal expressions
describe the same morphism in the category.  What we are dealing with, in effect,
is the word problem for a category described by generators and relations.  To present
the solution to this word problem we need the following preliminary definition.

\begin{definition}\label{binaryRelation:definition}
If $a$ and $b$ are distinct elements of $\{1,2,\dots, k\}$ and
$A$ is an object of ${\cal M}_n(k)$, we say that $a\Box_i b$ {\it in\/} $A$ if the
restriction functor ${\cal M}_n(k)\to{\cal M}_n(\{a,b\})$ sends $A$ to $a\Box_i b$.
\end{definition}

\begin{theorem}[Coherence Theorem for $n$-fold Monoidal Categories]\label{coherence:theorem}
Let $A$ and $B$ be objects of ${\cal M}_n(k)$.  Then
\begin{enumerate}
\item There is at most one morphism $A\to B$
\item A necessary and sufficient condition for the existence of a morphism $A\to B$
is that for any two elements $a,b$ in $\{1, 2,\dots, k\}$ if $a\Box_i b$ in $A$, then in $B$
either $a\Box_j b$ for some $j\ge i$ or $b\Box_j a$ for some $j>i$.
\end{enumerate}
\end{theorem}

\begin{example}
There is a morphism $A=(2\Box_2 3)\Box_1 1\longrightarrow 2\Box_2 1\Box_2 3 = B$
in ${\cal M}_2(3)$ since
\begin{itemize}
\item $2\Box_1 1$ in $A$ and $2\Box_2 1$ in $B$
\item $3\Box_1 1$ in $A$ and $1\Box_2 3$ in $B$
\item $2\Box_2 3$ in $A$ and $2\Box_2 3$ in $B$
\end{itemize}
but there is no morphism from $A$ to $C=1\Box_2 3\Box_2 2$ since $2\Box_2 3$ in $A$, while
$3\Box_2 2$ in $C$.
\end{example}

\begin{remark}\label{octo:coherence} 
The first part of the Coherence Theorem asserts that any diagram built out of the natural transformations
$\eta^{ij}$ must commute.  The necessity of the conditions  in the second part of the Coherence Theorem is forced by existence of the restriction functors $R_{\{a,b\}}: M_n(k)\to M_n(\{a, b\})$,
ie. if there is a morphism $A\to B$ in $M_n(k)$, then there must be a morphism $R_{\{a,b\}}(A)\to R_{\{a,b\}}(B)$
in $M_n(\{a, b\})$.  It is far from obvious however, that these conditions are sufficient to insure the existence of a morphism $A\to B$.
\end{remark}

\begin{remark} The coherence theorem implies that the topological operad spaces ${\cal M}_n(k)$ are
nerves of finite posets, and hence are compact polyhedra.
\end{remark}

\begin{definition} We define the Milgram subspace $\overline{\cal J}_n(k)$ to be the full subcategory of
${\cal M}_n(k)$ whose objects are contained in the free monoid with respect to $\Box_1$ on the
free monoid with respect to $\Box_2$ \dots on the free monoid with respect to $\Box_n$ on the
set $\{1, 2,\dots, k\}$.  Thus the objects of $\overline{\cal J}_n(k)$ look like
$$((\dots\Box_3\dots)\Box_2\dots)\Box_1\dots
\Box_1((\dots\Box_3\dots)\Box_2\dots)$$
ie. the operation $\Box_1$ can only occur at the outermost level,
the operation $\Box_2$ can only occur at the next level,\dots,
the operation $\Box_n$ can only occur at the innermost level.
Equivalently we can define the Milgram subspace to be the full subcategory of ${\cal M}_n(k)$ consisting
of objects which can be written without parentheses using the operation precedence rules: $\Box_n$ has
the highest precedence, $\Box_{n-1}$ has the next highest precedence, \dots, $\Box_1$ has the lowest
precedence.
\end{definition}

\begin{remark} The collection of Milgram subspaces $\{\overline{\cal J}_n(k)\}_{k\ge 0}$ is {\bf not} a suboperad of the categorical
operad ${\cal M}_n$.  It is only closed under the actions of the symmetric groups and the unit
maps
$$s_j:\overline{\cal J}_n(k)\longrightarrow\overline{\cal J}_n(k-1)\quad j=1, 2,\dots, k.$$
In other words $\overline{\cal J}_n(k)$, or rather its nerve which we also denote by $\overline{\cal J}_n(k)$, is a
preoperad in the sense of Berger \cite{Be}.  This structure is sufficient to define the premonad
construction
$$\overline{J}_n(X)=\left.\left(\coprod_{k\ge 0} \overline{\cal J}_n(k)\times_{\Sigma_k}X^k\right)\right/\approx$$
where $X$ is any based space.
If $\overline{\cal J}_n$ were an operad, this construction would be a monad, but this isn't the case here.
The notion of preoperad and the associated premonad construction were introduced in \cite{CMT}, where
preoperads are called ``coefficient systems'' (also cf. \cite{May}).
\end{remark}

In \cite{Mi} Milgram defined a construction
$$J_n(X)=\left.\left(\coprod_{k\ge 0}\left(P_k\right)^{n-1}\times X^k\right)\right/\approx$$
on based spaces $X$, where $P_k$ denotes the permutohedron: the convex hull in $\Bbb R^k$
of the $\Sigma_k$ orbit of a point such as $(1,2,\dots,k)$, all of
whose coordinates are distinct.  ($P_k$ is a $k-1$-dimensional cell.  Milgram uses the notation $C(k+1)$ to
denote $P_k$.) He showed that if $X$ is connected,
then $J_n(X)$ has the weak homotopy type of $\Omega^n\Sigma^n(X)$.

\begin{theorem}\label{Milgram2:theorem} For all spaces $X$, there is a natural homeomorphism
$$\overline{J}_2(X)\cong J_2(X).$$
\end{theorem}

Unfortunately for $n>2$ this does not hold. It turns out that our construction
$\overline{J}_n(X)$ is a natural quotient of Milgram's construction, and may be
thought of as a sort of {\it thin version\/} of the Milgram construction.

To understand the connection between $\overline{J}_n(X)$ and $J_n(X)$, we have
to consider yet another variant form of the Milgram construction, which we will call the
{\it thick Milgram construction\/} and denote $\widetilde{J}_n(X)$.
This is defined as the premonad construction on a preoperad
$\{\widetilde{\cal J}_n(k)\}_{k\ge 0}$ where
$$\widetilde{\cal J}_n(k)=\left.P_k^{n-1}\times\Sigma_k\right/\approx$$
where the equivalence relation glues together the $k!$ copies of $P_k^{n-1}$
along certain codimension 1 faces in the boundary.

\begin{theorem}\label{MilgramN:theorem}
\begin{enumerate}
\item There are natural quotient maps
$$\widetilde{J}_n(X)\stackrel{q_1}{\longrightarrow} J_n(X)\stackrel{q_2}{\longrightarrow}
\overline{J}_n(X)$$
which are homotopy equivalences.
\item Each of the variant forms of the Milgram construction arises from a preoperad
having the generic form
$$\left.D^{(n-1)(k-1)}\times\Sigma_k\right/\approx$$
where the equivalence relation glues together $k!$ copies of the $(n-1)(k-1)$-dimensional disk $D^{(n-1)(k-1)}$ along certain codimension 1 faces in the
boundary. The quotient maps $q_1$ and $q_2$ induce equivalences of preoperads.
\end{enumerate}
\end{theorem}

We would like to note that an earlier version of this paper suffered from some
confusion about the relation between the Milgram construction and the preoperad
$\{\overline{\cal J}_n(k)\}$. We would like to thank Clemens Berger for clearing up
this point.

Our main result is the following

\begin{theorem}\label{main:theorem} There is a chain of operad equivalences
$$
{\cal M}_n(k)
\stackrel{\simeq}{\longleftarrow}
\mbox{hocolim}_{{\cal M}_n(k)} F 
\stackrel{\simeq}{\longrightarrow}
\mbox{colim}_{{\cal M}_n(k)} F \stackrel{\simeq}{\hookrightarrow}{\cal C}_n(k),
$$
where ${\cal C}_n(k)$ denotes the little $n$-cubes operad of Boardman and Vogt
(and $F:{\cal M}_n(k)\to\bold{Top}$ is a functor we construct in Chapter
\ref{little:cubessection}). Moreover the inclusion of the Milgram preoperad
$\overline{\cal J}_n(k)$ in the operad ${\cal M}_n(k)$ is an equivalence of preoperads.
\end{theorem}

This gives a more definitive way of showing that the group completion of the nerve of an $n$-fold monoidal
category is an $n$-fold loop space.  For the proof we have given in the preceding section leaves open
the possibility that the group completion of the nerve of an $n$-fold monoidal category might have more
structure than that of an $n$-fold loop space (eg. perhaps it might be an infinite loop space).  This is a
serious possibility, since as we have noted in Section 1, slightly variant definitions of the notion of $n$-fold
monoidal category do indeed correspond to infinite loop spaces rather than $n$-fold loop spaces. The proof
based on Theorem \ref{main:theorem} rules out this possibility, since it shows that the free $n$-fold loop spaces 
$\Omega^n\Sigma^nX$, where $X$ is a discrete space do arise as group completions of $n$-fold monoidal
categories. In a subsequent paper we will show that in fact any $n$-fold loop space can be realized in this way.

\begin{remark}\label{inversion:kills}
Joyal and Street \cite{JS} noted that their theory of iterated monoidal categories collapses to the theory
of symmetric monoidal categories when $n>2$.  The reason for this is that their theory requires that
the interchange natural transformations $\eta^{ij}_{A,B,C,D}$ be isomorphisms.  Hence the categorical
operad for their theory is essentially obtained from our operad ${\cal M}_n$ by inverting all the morphisms.
But inverting all the morphisms in a category has the effect of killing off all the higher homotopy groups
of its nerve, leaving only the fundamental group intact (cf. Quillen \cite{Q}).  But according to Theorem
\ref{main:theorem} the homotopy groups of ${\cal M}_n$ are isomorphic to those of ${\cal C}_n$.  The
spaces of ${\cal C}_2$ are $K(\pi,1)$'s whereas the spaces of ${\cal C}_n$ are simply connected for
$n>2$.  Thus inverting all the morphisms in ${\cal M}_2$ does not change its homotopy type, since all its
higher homotopy groups are trivial anyway.  But inverting the morphisms of ${\cal M}_n$ for $n>2$ kills
off all the homotopy, rendering them into trivial categories, which endows iterated monoidal categories on
which they act with a symmetric monoidal structure.
\end{remark}

A related result is

\begin{theorem}\label{graphs:theorem}There is a chain of operad equivalences
$${\cal M}_n\stackrel{\simeq}{\longrightarrow}{\cal K}^{(n)}\stackrel{\simeq}{\longrightarrow}\Gamma^{(n)}$$
where ${\cal K}^{(n)}$ denotes the $n$-th filtration of Berger's complete graph operad
(cf. \cite{Be}) and $\Gamma^{(n)}$ is the $n$-th Smith filtration of the operad which
parametrizes (strict) symmetric monoidal categories (cf. \cite{Sm}).
\end{theorem}

\begin{remark} The homotopy type of the Milgram preoperad in the case $n=2$
was determined by Salvetti \cite{Sal}, in the more general context
of complements of hyperplane arrangements (also cf. \cite{CD}).
\end{remark}

\newpage
\MySection{The Coherence Theorem for $n$-fold Monoidal Categories}

This section  is devoted to the proof of Theorem \ref{coherence:theorem}, the coherence theorem for $n$-fold monoidal categories.
Before we proceed to the proof however, it will be convenient to reformulate the theorem.

\begin{definition} Let $\widehat{{\cal M}}_n(k)$ denote the category with the same objects as ${\cal M}_n(k)$,
but whose morphisms are as given in Theorem \ref{coherence:theorem}.  That is, there is a (unique) morphism between objects
$A\to B$ if and only if for any two elements $a,b$ in $\{1, 2,\dots, k\}$ if $a\Box_i b$ in $A$, then in $B$
either $a\Box_j b$ for some $j\ge i$ or $b\Box_j a$ for some $j>i$.  Note that by definition $\widehat{{\cal M}}_n(k)$
is a poset. More generally, following Definition \ref{set:asLabel}, we can define a similar category $\widehat{{\cal M}}_n(S)$
for any finite set $S$.  Note that if $S$ and $T$ are disjoint, then there are induced functors
$$\Box_i:\widehat{{\cal M}}_n(S)\times\widehat{{\cal M}}_n(T)\longrightarrow\widehat{{\cal M}}_n(S\amalg T)$$
for $i=1,2,\dots,n$.
\end{definition}

It follows immediately from Remark \ref{octo:coherence} that there is a functor
$$\Lambda_S^n: {\cal M}_n(S)\longrightarrow\widehat{{\cal M}}_n(S)$$
given by the identity on objects (which we will denote simply as $\Lambda_k^n$ if $S=\{1,2,\dots,k\}$).
Then the following is an obvious reformulation of Theorem \ref{coherence:theorem}

\begin{theorem}[Reformulation of the Coherence Theorem for $n$-fold Monoidal Categories]
\label{coherenceTheorem:too} The \linebreak functor
$$\Lambda_S^n: {\cal M}_n(S)\longrightarrow\widehat{{\cal M}}_n(S)$$
is an isomorphism of categories.
\end{theorem}

As we noted in Definition \ref{set:asLabel}, since ${\cal M}_n(S)$ only depends on the cardinality of $S$, it suffices to
prove the coherence theorem for $\Lambda_k^n$.  However it is convenient to recast our basic induction
hypothesis in terms of $\Lambda_S^n$:

\begin{hypothesis}
We assume that $\Lambda_S^n$ is an isomorphism for every proper subset $S\subset\{1,2,\dots,k\}$.
\end{hypothesis}

We note that the coherence theorem is trivially true when $k=1$ and the octahedral picture of ${\cal M}_n(2)$ given
in the preceding section shows that it is also true when $k=2$.  This starts our induction going.

\begin{definition}
If $A\Box_i B$ is an object in ${\cal M}_n(S)$, we denote by $|A|$ the subset of $S$ consisting
of all the generators present in $A$.  Thus by definition
$$S=|A|\amalg |B|$$
We will say that $A\in{\cal M}_n(T)$ is a {\it partial\/} object of ${\cal M}_n(S)$ if $T\subset S$.
\end{definition}

We begin with a few basic observations about the categories ${\cal M}_n(k)$.

\begin{lemma}\label{object:decomposition}
Suppose that $X$ is an object of ${\cal M}_n(S)$ and that there is a partition $S=S_1\amalg S_2$
such that for any $x\in S_1$ and any $y\in S_2$, $x\Box_i y$ in $X$.  Then $X$ has a decomposition
$$X=X_1\Box_i X_2$$
with  $|X_1|=S_1$ and $|X_2|=S_2$.
\end{lemma}

The proof is left as an easy exercise for the reader.  (Hint: use induction on the cardinality of $S$.)

\begin{definition}\label{partial:objects}
Let $A$ and $B$ be two partial objects in ${\cal M}_n(k)$.
\begin{enumerate}
\item The {\it difference} $A - |B|$ is the restriction of $A$ to ${\cal M}_n(|A|-|B|)$, ie. it
is the object obtained from $A$ by zeroing out all generating objects of $A$ which are
also present in $B$.
\item The {\it intersection} $A \cap |B|$ is defined to be 
$A - \left|A - |B|\right|$, ie. the object obtained from $A$ by zeroing out all generating objects of $A$
which are not present in $B$.
\end{enumerate}
\end{definition}

\begin{proposition}\label{coherence:basicObservation}
Let $f : A \Box_i B \longrightarrow C \Box_j D$ be a morphism in
${\cal M}_n(k)$.
\begin{enumerate}
\item If $A$, $B$, $C$ and $D$ are all different from $0$, then
$j \geq i$;

\item If $j = i$ and $\mbox{card}(|A|) = \mbox{card}(|C|)$ then there exist two morphisms
$g : A \longrightarrow C$ and $h : B \longrightarrow D$ in ${\cal M}_n(|A|)$ and
${\cal M}_n(|B|)$, respectively, such that $f = g \Box_i h$ (we
shall call such a morphism $f$ a {\it $\Box_i$--split} morphism).

\item If $j=i$ and $\mbox{card}(|A|) >\mbox{card}(|C|)$ then there exist two morphisms
$g:A\to C\Box_i (A-|C|)$ and $h:(A-|C|)\Box_i B\to D$ so that $f$ factors as
the composite
$$\diagram
A\Box_i B \rrto^(.4){g\Box_i id_B} &&C\Box_i (A-|C|)\Box_i B
\rrto^(.6){id_C \Box_i h} &&C\Box_i D
\enddiagram$$

\item If $j=i$ and $\mbox{card}(|A|) < \mbox{card}(|C|)$ then there exist two morphisms
$g:B\to (C-|A|)\Box_i D$ and $h:A\Box_i (C-|A|)\to C$ so that $f$ factors as
the composite
$$\diagram
A\Box_i B \rrto^(.4){id_A\Box_i g} &&A\Box_i (C-|A|)\Box_i D
\rrto^(.6){h\Box_i id_D} &&C\Box_i D
\enddiagram$$

\end{enumerate}
\end{proposition}

\noindent
{\bf Proof:}
By definition any morphism in ${\cal M}_n(k)$ is a composition of nontrivial morphisms of the
form $\eta^{ij}_{X,Y,Z,W}$ and $f_1\Box_i f_2$ where exactly one of $f_1$ or $f_2$ is
an identity map (of a nonzero object).  We shall refer to such morphisms as
{\it indecomposable morphisms\/}.  To prove part (1) it suffices to prove it for indecomposable
morphisms.  Now the assertion is evidently true for indecomposable morphisms of the form
$f_1\Box_i f_2$.  For nontrivial morphisms of the form
$$\eta^{ij}_{X,Y,Z,W}: (X\Box_j Y)\Box_i(Z\Box_j W)\longrightarrow (X\Box_i Z)\Box_j (Y\Box_i W)$$
the outer operation in the source object is $i$ and the outer operation in the target object is $j>i$,
since by the unit conditions $\eta^{ij}_{X,Y,Z,W}$ is the identity if any of the objects $X\Box_j Y$,
$Z\Box_j W$, $X\Box_i Z$, $Y\Box_i W$ are equal to $0$.

To check part (2),
note first that the conditions $j = i$ and $\mbox{card}(|A|) = \mbox{card}(|C|)$ imply that $|A|=|C|$ and $|B|=|D|$.
For otherwise there would have to exist elements $x\in |A|\cap|D|$ and $y\in |B|\cap|C|$ and then
we would have $x\Box_i y$ in the source object $A\Box_i B$ but $y\Box_i x$ in the target object
$C\Box_i D$, which is precluded by the very existence of the functor $\Lambda_k^n$.
If we then factor $f$ into indecomposable morphisms
$$A\Box_i B\longrightarrow X_1\longrightarrow X_2\longrightarrow\dots\longrightarrow X_{m-1}
\longrightarrow C\Box_i D,$$
it follows directly from Lemma \ref{object:decomposition} that each intermediate object $X_r$
has a decomposition $X_r=X_r'\Box_i X_r''$ with $|X_r'|=|A|=|C|$ and $|X_r''=|B|=|C|$.
This reduces proving part (2) to the case when $f$ is indecomposable.  By the argument of
the preceding paragraph, $f$ would then have to have the form $f=f_1\Box_i f_2$, for some
possibly different $\Box_i$ decomposition of the objects $A\Box_i B$ and $C\Box_i D$.  But
in that case an easy argument using induction hypothesis (IH.1) shows that the decomposition
$f=f_1\Box_i f_2$ can be reparanthesized to a decomposition $f=g\Box_i h$ of the requisite form.

To check part (3), we first demonstrate that $f$ factors through an
object $X\Box_i Y \Box_i Z$ such that $|X|=|C|$, $|Y|=|A|-|C|$ and
$|Z|=|B|$.  Begin by factoring $f$ as
$$A\Box_i B \stackrel{f_1}{\longrightarrow} W
\stackrel{f_2}{\longrightarrow} C\Box_i D$$
with $W$ having a maximal number of $\Box_i$ summands. Now factor
$f_2$ into indecomposable morphisms as
$$W=U_0\longrightarrow U_1\longrightarrow U_2\longrightarrow\dots
\longrightarrow U_m=C\Box_i D$$
We claim that for each $U_p$ and any decomposition $U_p=U_p'\Box_i U_p''$
there is a corresponding decomposition $W=W'\Box_iW''$ with $|W'|=|U_p'|$
and $|W''|=|U_p''|$.  If not, let $U_p$ be the first object in the chain
having a decomposition $U_p=U_p'\Box_i U_p''$ incompatible with $W$.  Since
the morphism
$$U_{p-1}\longrightarrow U_p$$
must be $\Box_i$-split, there is another decomposition $U_p=V_p'\Box_i V_p''$
which is compatible with $U_{p-1}$ and hence with $W$.  Let $W=W'\Box_iW''$
be the compatible decomposition with $|W'|=|V_p'|$ and $|W''|=|V_p''|$.
Then according to part (2) $f_2$ factors as
$$\diagram
W=W'\Box_iW''\rrto^{f_2'\Box_i f_2''}
&&U_p=V_p'\Box_i V_p''\rrto &&C\Box_i D
\enddiagram$$
Then the incompatible decomposition $U_p=U_p'\Box_i U_p''$ must give
either a decomposition of $V_p'$ which is incompatible with that of
$W'$ or a decomposition of $V_p''$ which is incompatible with that of
$W''$.  In the first case IH.1 produces an intermediate object $G'$ between
$W'$ and $V_p'$ with more $\Box_i$ summands than $W'$. Hence $f$ factors
through $G'\Box_i W''$ which has more $\Box_i$ summands than $W$, which
contradicts our choice of $W$.  In the second case we obtain a similar
contradiction. This proves the claim.  Applying this to the decomposition
$U_m=C\Box_i D$, we obtain a compatible decomposition $W=X\Box_i T$ with
$|X|=|C|$ and $|T|=|D|$.

A similar argument on $f_1$ yields another
decomposition $W=S\Box_i Z$ with $|S|=|A|$ and $|Z|=|B|$. Combining the
two decompositions yields the desired decomposition $W=X\Box_i Y\Box_i Z$
with $|Y|=|A|-|C|$.  Then (2) yields a factorization of $f$ as
$$\diagram
A\Box_i B\rrto^{g_1\Box_i g_2} && X\Box_i Y\Box_i Z\rrto^{h_1\Box_i h_2}
&&C\Box_i D
\enddiagram$$
for some morphisms $g_1:A\longrightarrow X\Box_i Y$, $g_2:B\longrightarrow Z$,
$h_1:X\longrightarrow C$ and $h_2:Y\Box_i Z\longrightarrow D$. This yields
another factorization of $f$ as
$$\diagram
A\Box_i B\rrto^{g'\Box_i id_B} && C\Box_i Y\Box_i B\rrto^{id_C\Box_i h'}
&&C\Box_i D
\enddiagram$$
where $g'=(h_1\Box_i id_Y)g_1:A\longrightarrow C\Box_i Y$ and
$h'=h_2(id_Y\Box_i g_2):Y\Box_i B\longrightarrow D$.  Then IH.1 yields
a further factorization of $g'$ as
$$\diagram
A\rrto^g &&C\Box_i(A-|C|)\rrto^{id_C\Box_i l} && C\Box_i Y
\enddiagram.$$
Then the desired factorization
$$\diagram
A\Box_i B \rrto^(.4){g\Box_i id_B}
&&C\Box_i (A-|C|)\Box_i B\rrto^(.6){id_C \Box_i h} &&C\Box_i D
\enddiagram$$
is obtained by setting $h=h'(l\Box_i id_B)$.  This concludes the proof of
(3).  The proof of (4) is similar.
\hfill
$\Box$

\begin{remark}\label{basicObservation:remark}
The results listed in Proposition \ref{coherence:basicObservation} are also true in the category
$\widehat{{\cal M}}_n(k)$, but in this case they follow immediately
from the conditions that have to be satisfied by any two objects
which are, respectively, the source and the target of a certain
morphism.
\end{remark}

\begin{remark}\label{compatible:decomposition}
By similar arguments one can show that given any morphism
$f: A\Box_i B\longrightarrow C\Box_i D$ in ${\cal M}_n(k)$ with $i=1$ or
$i=n$, there are compatible $\Box_i$ decompositions of the source and
target, and hence by (2) $f$ has a nontrivial decomposition $f=f_1\Box_i f_2$.
\end{remark}

\begin{definition}
Let $\mu : A \Box_i B \longrightarrow C$ be a morphism in
$\widehat{{\cal M}}_n(k)$ with $|A|$ having cardinality $p$ and $|B|$ having cardinality $q$ (so $p+q=k$).
We say that $\mu$ is a {\it strong $(p, q)$--shuffle} if :
$$
C - |B| = A \qquad \mbox{and} \qquad C - |A| = B
$$
Note that this means that the order in which the generating objects appear in $C$ is a $(p, q)$--shuffle
(in the standard sense) of the order in which they appear in $A$ and $B$.  However it means that in
addition the operations appearing in $C$ are in some sense the operations appearing in $A$ and $B$ shuffled
together.
\end{definition}

\begin{remark}
The notion of strong shuffle defined above assumes implicitly the
existence of {\it at most one} morphism between any two objects of
the category $\widehat{{\cal M}}_n(k)$. This is why we cannot define it
{\it a priori} in the category ${\cal M}_n(k)$.
\end{remark}

\begin{proposition}\label{strong:shuffle}
Let $\mu : A \Box_i B \longrightarrow C$ be a morphism in
$\widehat{{\cal M}}_n(k)$. Then the following conditions are
equivalent :
\begin{enumerate}
\item $\mu$ is a strong shuffle;

\item There is no nontrivial factorization of $\mu$ as
$$\diagram
A \Box_i B \rrto^{\mu_1 \Box_i \mu_2}
&&X \Box_i Y \rrto^{\xi} &&C
\enddiagram$$
with $\mu_1:A\longrightarrow X$ and $\mu_2:B\longrightarrow Y$.
\end{enumerate}
\end{proposition}

\noindent
{\bf Proof:} Let $A'$ and $B'$ denote, respectively, the objects
$C - |B|$ and $C - |A|$. Then $\mu$ obviously factors as
$$
A \Box_i B \longrightarrow A' \Box_i B' \longrightarrow C.
$$
Therefore conditions (1) and (2) are equivalent.
\hfill
$\Box$

\begin{definition}
An object (or a partial object) $A$ in ${\cal M}_n(k)$ is called :

\begin{enumerate}
\item $\Box_i$--{\it reducible} if it can be expressed
nontrivially as $A_1 \Box_i A_2$;

\item $\Box_i$--{\it irreducible} if it is not
$\Box_i$--reducible.
\end{enumerate}
\end{definition}

\begin{definition}
A morphism $f : A \Box_i B \longrightarrow C \Box_r D$ in ${\cal M}_n(k)$
(or a morphism $\mu : A \Box_i B \longrightarrow C \Box_r D$ in
$\widehat{{\cal M}}_n(k)$) is called :

\begin{enumerate}
\item{\it irreducible} if $r > i$ and all the objects through
which $f$ ( or $\mu$) factors nontrivially are
$\Box_j$--irreducible for all $j \in \{i,\ i+1,\ ...,\ r-1\}$;

\item {\it reducible} if it is not irreducible.
\end{enumerate}
\end{definition}

As we shall see below, we can't get very far with our basic induction
hypothesis (IH.1). We have to use {\it  double
induction}, the second inductive hypothesis being related to the
outermost operation in the targets of the morphisms to be considered.
More precisely, we need :

\begin{hypothesis}
Let $r \geq 2$ be a positive integer.  Then
$$\Lambda_S^n:\mbox{Hom}_{{\cal M}_n(k)}(A,B)\longrightarrow
\mbox{Hom}_{\widehat{{\cal M}}_n(k)}(A,B)$$
is a bijection, whenever $B$ is $\Box_j$--reducible
with $j < r$.
\end{hypothesis}

\begin{remark}\label{compatible:splitting}
Note that according to Remark \ref{compatible:decomposition} and Proposition \ref{coherence:basicObservation}
if $j=1$ and $B=B_1\Box_1 B_2$, then
there is a compatible splitting $A=A_1\Box_1 A_2$ and any morphism $f:A\to B$
in ${\cal M}_n(S)$ must also split $f=f_1\Box_1 f_2$.  Thus in this case, our first
induction hypothesis (IH.1) implies that $\Lambda_S^n$ is a bijection.  This starts
our second induction hypothesis.
Note also that this argument proves more: namely that $\Lambda_S^n$ is bijective
on morphisms where the source and target have the same outermost operation $\Box_r$ and
compatible $\Box_r$ splittings.  Moreover this also holds when the source and target have
the same outermost operation even when there are no compatible splittings.  For by
Proposition \ref{coherence:basicObservation} parts (3) and (4), in that case one can
insert a canonical intermediate object having splittings compatible with both source and
target through which all morphisms must factor.
\end{remark}

\begin{lemma}\label{extAss:Lemma}
\begin{enumerate}
\item Let $\mu:A\Box_i B\longrightarrow C$ be a strong shuffle
in $\widehat{{\cal M}}_n(k)$, with $C$ $\Box_r$-reducible. Then for any splitting
$C=C_1\Box_r C_2$ there are (possibly trivial) splittings
$A=A_1\Box_r A_2$ and $B=B_1\Box_r B_2$ and morphisms
$g_1:A_1\Box_i B_1\longrightarrow C_1$ and $g_2:A_2\Box_i B_2\longrightarrow C_2$
in ${\cal M}_n(k)$ such that $\mu$ lifts to the composite
$$\diagram
A\Box_i B=(A_1\Box_r A_2)\Box_i(B_1\Box_r B_2)\rrto^(.55){\eta^{ir}_{A_1,A_2,B_1,B_2}}
&&(A_1\Box_i B_1)\Box_r(A_2\Box_i B_2)\rrto^(.6){g_1\Box_r g_2}&&C_1\Box_r C_2
\enddiagram$$
\item Any diagram of the form
$$\diagram
A\Box_i B\rrto^(.4){\eta^{ir}_{A_1,A_2,B_1,B_2}}\dto^{\eta^{ir}_{A'_1,A'_2,B'_1,B'_2}}
&&(A_1\Box_i B_1)\Box_r(A_2\Box_i B_2)\dto^g\\
(A'_1\Box_i B'_1)\Box_r(A'_2\Box_i B'_2)\rrto^(.6)h &&C
\enddiagram$$
with both $\eta^{ir}_{A_1,A_2,B_1,B_2}$ and $\eta^{ir}_{A'_1,A'_2,B'_1,B'_2}$ nontrivial,
commutes in ${\cal M}_n(k)$.
\end{enumerate}
\end{lemma}

\noindent
{\bf Proof:} Note first that the objects $A$, $B$, $C$ can be decomposed into
$\Box_r$--irreducible objects as follows :
$$
\begin{array}{cc}
A & = \overline A_1 \Box_r \overline A_2 \Box_r \dots
\Box_r \overline A_s \\
B & = \overline B_1 \Box_r \overline B_2 \Box_r \dots
\Box_r \overline B_t \\
C & = \overline C_1 \Box_r \overline C_2 \Box_r \dots
\Box_r \overline C_u
\end{array}
$$
with $s,\ t \ge 1$ and $u \geq 2$. Since $\mu$ is a strong shuffle
there exist nondecreasing functions
$$
\begin{array}{cc}
\sigma & : \{1,\ 2,...,\ s\} \longrightarrow \{1,\ 2,...,\ u\} \\
\tau & : \{1,\ 2,...,\ t\} \longrightarrow \{1,\ 2,...,\ u\}
\end{array}
$$
defined, respectively, by the relations
$$
\begin{array}{cc}
|\overline A_j| & \subset |\overline C_{\sigma (j)}| \text{, for all}
\ j \in \{1,\ 2,...,\ s\} \\
|\overline B_j| & \subset |\overline C_{\tau (j)}| \text {, for all}
\ j \in \{1,\ 2,...,\ t\}
\end{array}
$$
Then 
$$C_1=\overline C_1\Box_r \overline C_2 \Box_r\dots \Box_r\overline C_v$$
$$C_2=\overline C_{v+1}\Box_r\overline C_{v+2}\Box_r\dots\Box_r\overline C_u$$
for some $v$. Now define the objects $A_1$, $A_2$, $B_1$, $B_2$ by
$$
\begin{array}{cc}
A_1 & := \Box_r \{\overline A_j\ |\ \sigma(j)\le v \}\\
A_2 & := \Box_r \{\overline A_j\ |\ \sigma(j)> v \} \\
\end{array}
\qquad
\qquad
\qquad
\begin{array}{cl}
B_1 & := \Box_r \{\overline B_j\ |\ \tau(j)\le v \}\\
B_2 & := \Box_r \{\overline B_j\ |\ \tau(j)> v \} \\
\end{array}
$$
Then clearly there are morphisms $A_1\Box_i B_1\longrightarrow C_1$
and $A_2\Box_i B_2\longrightarrow C_2$ in $\widehat{{\cal M}}_n(k)$.
By IH.1 these morphisms have lifts $g_1:A_1\Box_i B_1\longrightarrow C_1$
and $g_2:A_2\Box_i B_2\longrightarrow C_2$ in ${\cal M}_n(k)$. It is
immediate that
$$\diagram
A\Box_i B=(A_1\Box_r A_2)\Box_i(B_1\Box_r B_2)\rrto^(.55){\eta^{ir}_{A_1,A_2,B_1,B_2}}
&&(A_1\Box_i B_1)\Box_r(A_2\Box_i B_2)\rrto^(.6){g_1\Box_r g_2}&&C_1\Box_r C_2
\enddiagram$$
is a lift of $\mu$.

To prove part (2) we begin with a definition. We say that the chain in ${\cal M}_n(k)$
$$\diagram
A\Box_i B=(A'_1\Box_r A'_2)\Box_i(B'_1\Box_r B'_2)\rrto^(.55){\eta^{ir}_{A'_1,A'_2,B'_1,B'_2}}
&&(A'_1\Box_i B'_1)\Box_r(A'_2\Box_i B'_2)\rrto^(.7)h&&C
\enddiagram$$
is {\it subordinate\/} to the chain
$$\diagram
A\Box_i B=(A_1\Box_r A_2)\Box_i(B_1\Box_r B_2)\rrto^(.55){\eta^{ir}_{A_1,A_2,B_1,B_2}}
&&(A_1\Box_i B_1)\Box_r(A_2\Box_i B_2)\rrto^(.7)g&&C
\enddiagram$$
if there exists a splitting $A'_1=A_1\Box_r A''_1$ and $B'_1=B_1\Box_r B''_1$.
We first show that the diagram in part (2) is commutative in this case.  In other
words, if one chain is subordinate to the other then their composites are equal.

Note that our hypothesis implies that $A_2=A''_1\Box_r A'_2$ and
$B_2=B''_1\Box_r B'_2$. The existence of the morphisms
$$g:(A_1\Box_i B_1)\Box_r\bigl((A''_1\Box_r A'_2)\Box_i(B''_1\Box_r B'_2)\bigr)
\longrightarrow C$$
$$h:\bigl((A_1\Box_r A''_1)\Box_i(B_1\Box_i B''_1)\bigr)\Box_r(A'_2\Box_i B'_2)
\longrightarrow C$$
implies the existence of a morphism
$$(A_1\Box_i B_1)\Box_r(A''_1\Box_i B''_1)\Box_r(A'_2\Box_i B'_2)
\longrightarrow C$$
in $\widehat{{\cal M}}_n(k)$ which according to Remark \ref{compatible:splitting} has a unique lift
$$l:(A_1\Box_i B_1)\Box_r(A''_1\Box_i B''_1)\Box_r(A'_2\Box_i B'_2)
\longrightarrow C$$
in ${\cal M}_n(k)$. This then yields the following diagram in ${\cal M}_n(k)$
$$
\spreaddiagramrows{-.8pc}
\spreaddiagramcolumns{-5pc}
\diagram
&&(A_1\Box_r A''_1 \Box_r A'_2)\Box_i(B_1\Box_r B''_1 \Box_r B'_2)
\dllto\drrto\\
(A_1\Box_i B_1)\Box_r\bigl((A''_1\Box_r A'_2)\Box_i(B''_1\Box_r B'_2)\bigr)
\drrto\xto[dddddddrr]_g
&&{\squarify{1}\circled}
&&\bigl((A_1\Box_r A''_1)\Box_i(B_1\Box_i B''_1)\bigr)\Box_r(A'_2\Box_i B'_2)
\dllto\xto[dddddddll]^h\\
&&(A_1\Box_i B_1)\Box_r(A''_1\Box_i B''_1)\Box_r(A'_2\Box_i B'_2)\xto[dddddd]^l\\
&{\squarify{2}\circled} &&{\squarify{3}\circled}\\
\\
\\
\\
\\
&&C
\enddiagram$$
where the unlabelled arrows are those which occur in the external associativity
diagram. Then the left hand side of the diagram is one of our given chains and
the right hand side is the other (subordinate) one. This diagram commutes because
all the inner diagrams commute: the square (1) by the external associativity
law and the two triangles (2) and (3) by Remark \ref{compatible:splitting}.

Now given two arbitrary chains
$$\diagram
A\Box_i B=(A_1\Box_r A_2)\Box_i(B_1\Box_r B_2)\rrto^(.55){\eta^{ir}_{A_1,A_2,B_1,B_2}}
&&(A_1\Box_i B_1)\Box_r(A_2\Box_i B_2)\rrto^(.7)g&&C
\enddiagram$$
$$\diagram
A\Box_i B=(A'_1\Box_r A'_2)\Box_i(B'_1\Box_r B'_2)\rrto^(.55){\eta^{ir}_{A'_1,A'_2,B'_1,B'_2}}
&&(A'_1\Box_i B'_1)\Box_r(A'_2\Box_i B'_2)\rrto^(.7)h&&C
\enddiagram$$
one can construct a third chain
$$\diagram
A\Box_i B\xto[dd]^{\scriptstyle\eta^{ir}_{A_1\cap |A'_1|,A-|A_1|\cap |A'_1|,B_1\cap |B'_1|,B-|B_1|\cap |B'_1|}}\\
\\
\bigl((A_1\cap |A'_1|)\Box_i(B_1\cap |B'_1|)\bigr)\Box_r
\bigl((A-|A_1|\cap |A'|)\Box_i(B-|B_1|\cap |B'_1|)\bigr)\rrto^(.8)l
&&C
\enddiagram$$
provided that $A_1\cap |A'_1|$ and $B_1\cap |B'_1|$ are not simultaneously 0.
(Again we use Remark \ref{compatible:splitting} to construct $l$.)  In this case both of the given
chains are subordinate to this third one, hence the composites of the two
are both equal to the composite of the third, and thus equal to each other.

The remaining case is that both $A_1\cap |A'_1|=0=B_1\cap |B'_1|$. In this case we
may assume without loss of generality that $A'_1=0=B_1$.  Then again using
Remark \ref{compatible:splitting} we can again construct a third chain
$$\diagram
A\Box_i B\rrto^(.35){\eta^{ir}_{A_1,A_2,B'_1,B'_2}}
&&(A_1\Box_i B'_1)\Box_r(A_2\Box_i B'_2)\rrto^(.7)l &&C
\enddiagram$$
This time the third chain is subordinate to both of the original ones,
and again we get that their composites are equal. This completes the proof
of part (2).

\begin{lemma}\label{intAss:Lemma}
If $f : A \Box_i B \longrightarrow C$ is an irreducible morphism in
${\cal M}_n(k)$, with $C$ $\Box_r$-reducible, then :

\begin{enumerate}
\item $\ \ A$ and $B$ are $\Box_i$--irreducible objects
(in ${\cal M}_n(|A|)$ and ${\cal M}_n(|B|)$, respectively);

\item $f$ factors as a composite :
$$
\diagram
(A_1 \Box_r A_2) \Box_i (B_1 \Box_r B_2)
\xto[rrr]^{\eta^{ir}_{A_1, A_2, B_1, B_2}}
&&&(A_1 \Box_i B_1) \Box_r (A_2 \Box_i B_2)\dto^g \\
&&&C_1 \Box_r C_2
\enddiagram
$$
with $\Lambda_{|C_1|}^n(g_1)$ and $\Lambda_{|C_2|}^n(g_2)$
strong shuffles, where $g_1$, resp. $g_2$ are the restrictions of $g$
to $A_1 \Box_i B_1$, resp. $A_2 \Box_i B_2$;

\item $\ \ \Lambda^n_k(f)$ is a strong shuffle in
$\widehat{{\cal M}}_n(k)$.
\end{enumerate}
\end{lemma}

\noindent
{\bf Proof:}
Note first that the irreducible morphism $f$ can be obviously written
as a composite
$$
\diagram
(A_1 \Box_r A_2) \Box_i (B_1 \Box_r B_2)
\xto[rrr]^{\eta^{ir}_{A_1, A_2, B_1, B_2}}
&&&(A_1 \Box_i B_1) \Box_r (A_2 \Box_i B_2)\dto^g \\
&&&C
\enddiagram
$$
due to the structure of its source.

Next we shall prove that $A$ is $\Box_i$--irreducible. If not,
then one of the objects $A_1$ and $A_2$ must be equal to $0$
and the other one must be $\Box_i$--split. Without loss of
generality we may assume $A_1 = 0$ and
$A_2 = A_{21} \Box_i A_{22}$ with both $A_{21}$ and $A_{22}$
different from $0$ (since the other case follows from a similar
argument). Then the morphism $\eta^{ir}_{A_1, A_2, B_1, B_2}
= \eta^{ir}_{0, A_{21} \Box_i A_{22}, B_1, B_2}$ can be written as
the composite
$$
\diagram
A_{21} \Box_i A_{22} \Box_i(B_1 \Box_r B_2)
\xto[rrrr]^{id_{A_{21}} \Box_i \eta^{ir}_{0, A_{22}, B_1, B_2}}
&&&&A_{21} \Box_i (B_1 \Box_r (A_{22} \Box_i B_2))
\dto^{\eta^{ir}_{0, A_{21}, B_1, A_{22} \Box_i B_2}} \\
&&&&B_1 \Box_r (A_{21} \Box_i A_{22} \Box_i B_2)
\enddiagram
$$
according to the internal associativity law, contradicting the
irreducibility of $f$. Thus $A$ must be $\Box_i$--irreducible. In
a similar way one can obtain the same property for $B$, finishing the
proof of (1).

Next, suppose that at least one of the morphisms
$\Lambda_{|C_1|}^n(g_1)$, $\Lambda_{|C_2|}^n(g_2)$ is not a
strong shuffle. Then, as in the proof of Proposition \ref{strong:shuffle} and 
using Remark \ref{compatible:splitting}, $g$ can be factored as
$$\diagram
(A_1 \Box_i B_1) \Box_r (A_2 \Box_i B_2)
\xto[rrr]^{(a_1 \Box_i b_1) \Box_r (a_2 \Box_i b_2)}
&&&(X_1 \Box_i Y_1) \Box_r (X_2 \Box_i Y_2)\rrto^(.7){g'} &&C
\enddiagram$$
with at least one of the morphisms $a_1$, $b_1$, $a_2$, $b_2$
different from the corresponding identity, and we obtain the
following commutative diagram (by the naturality of $\eta^{ir}$) :
$$
\diagram
(A_1 \Box_r A_2) \Box_i (B_1 \Box_r B_2)
\xto[rrrr]^{\eta^{ir}_{A_1, A_2, B_1, B_2}}
\dto^{(a_1 \Box_r a_2) \Box_i (b_1 \Box_r b_2)}
&&&&(A_1 \Box_i B_1) \Box_r (A_2 \Box_i B_2)
\dto^{(a_1 \Box_i b_1) \Box_r (a_2 \Box_i b_2)} \\
(X_1 \Box_r X_2) \Box_i (Y_1 \Box_r Y_2)
\xto[rrrr]^{\eta^{ir}_{X_1, X_2, Y_1, Y_2}}
&&&&(X_1 \Box_i Y_1) \Box_r (X_2 \Box_i Y_2)
\enddiagram
$$
contradicting the irreducibility of $f$ and completing the proof of
(2).

To prove (3) assume that $\Lambda_k^n(f)$ is not a strong shuffle.  Then
we can factor $\Lambda_k^n(f)$ as
$$A\Box_i B\longrightarrow A'\Box_i B'\stackrel{\mu}{\longrightarrow} C.$$
From the fact that $\Lambda_{|C_1|}^n(g_1)$ and $\Lambda_{|C_2|}^n(g_2)$
are strong shuffles and Proposition \ref{coherence:basicObservation}, it follows that any splitting
$A'=A'_1\Box_r A'_2$ corresponds to a splitting $A=\tilde{A}_1\Box_r\tilde{A}_2$,
and similarly that any splitting $B'=B'_1\Box_r B'_2$ corresponds to a splitting
$B=\tilde{B}_1\Box_r\tilde{B}_2$.  According to Lemma \ref{extAss:Lemma} (1), $\mu$ can be
lifted to a composite
$$\diagram
A'\Box_i B'=(A'_1\Box_r A'_2)\Box_i(B'_1\Box_r B'_2)\rrto^(.55){\eta^{ir}_{A'_1,A'_2,B'_1,B'_2}}
&&(A'_1\Box_i B'_1)\Box_r(A'_2\Box_i B'_2)\rrto^(.7)h&&C
\enddiagram$$
Pick the corresponding splittings of $A=\tilde{A}_1\Box_r\tilde{A}_2$ and
$B=\tilde{B}_1\Box_r\tilde{B}_2$, then use IH.1 to lift to morphisms
$$l_1:\tilde{A}_1\to A'_1,\qquad l_2:\tilde{A}_2\to A'_2,\qquad
l_3:\tilde{B}_1\to B'_1,\qquad l_4:\tilde{B}_2\to B'_2.$$
Then we have the following diagram in ${\cal M}_n(k)$
$$
\spreaddiagramcolumns{-.42pc}
\diagram
{A\Box_i B}\xto[rrrr]^{\eta^{ir}_{A_1,A_2,B_1,B_2}}
\drrto^{\eta^{ir}_{\tilde{A}_1,\tilde{A}_2,\tilde{B}_1,\tilde{B}_2}}
\xto[ddd]^{(l_1\Box_r l_2)\Box_i(l_3\Box_r l_4)}
&&&&{(A_1\Box_i B_1)\Box_r(A_2\Box_i B_2)}\xto[dddd]^g\\
&&{(\tilde{A}_1\Box_i \tilde{B}_1)\Box_r(\tilde{A}_2\Box_i \tilde{B}_2)}
\xto[ddd]^{(l_1\Box_i l_3)\Box_r(l_2\Box_i l_4)}\\
&{\squarify{1}\circled}&&{\squarify{2}\circled}\\
{A'\Box_i B'}\drrto^{\eta^{ir}_{A'_1,A'_2,B'_1,B'_2}}\\
&&{(A'_1\Box_i B'_1)\Box_r(A'_2\Box_i B'_2)}\rrto^h && C
\enddiagram$$
This diagram commutes since the two inner diagrams commute: (1) by
naturality of $\eta^{ir}$ and (2) according to Lemma \ref{extAss:Lemma}(2). Since the
composite across the top and right is $f$, this contradicts the supposed
irreducibility of $f$.
\hfill
$\Box$

\begin{remark}\label{clearly:irreducible}
Clearly if $f: A \Box_i B \longrightarrow C$ is a morphism in ${\cal M}_n(k)$
such that $\Lambda_k^n(f)$ is irreducible, then $f$ is also
irreducible.
\end{remark}

\begin{lemma}\label{irreducible:uniqueLift}
If $\mu : A \Box_i B \longrightarrow C$ is an irreducible morphism in
$\widehat {{\cal M}}_n(k)$, then $\mu$ has a unique preimage
$f : A \Box_i B \longrightarrow C$ in ${\cal M}_n(k)$
\end{lemma}

\noindent
{\bf Proof:} By induction hypothesis (IH.2) we may as well assume that $C$ is $\Box_r$-reducible.
By Proposition \ref{strong:shuffle}, $\mu$ is a strong shuffle.
By Lemma \ref{extAss:Lemma}(1), $\mu$ has at least one preimage $f$.  But
$\mu$ irreducible implies that any preimage $f$ is also irreducible.
This in turn implies that any preimage $f$ must be a composite
of the form:
$$\diagram
A\Box_i B=(A_1\Box_r A_2)\Box_i(B_1\Box_r B_2)\rrto^(.55){\eta^{ir}_{A_1,A_2,B_1,B_2}}
&&(A_1\Box_i B_1)\Box_r(A_2\Box_i B_2)\rrto^(.7)g&&C
\enddiagram$$
By Lemma \ref{extAss:Lemma}(2) it follows that $f$ is unique.
\hfill
$\Box$

\begin{remark}\label{standard:lift}
The reader might wonder why the same argument doesn't show
that any strong shuffle $\mu: A\Box_i B\longrightarrow C$ has a unique preimage
in ${\cal M}_n(k)$. According to Lemma \ref{extAss:Lemma}, $\mu$ has a preimage of the form
$$
\diagram
(A_1 \Box_r A_2) \Box_i (B_1 \Box_r B_2)
\xto[rrr]^{\eta^{ir}_{A_1, A_2, B_1, B_2}}
&&&(A_1 \Box_i B_1) \Box_r (A_2 \Box_i B_2)\dto^g \\
&&&C
\enddiagram
$$
and any two such preimages are equal.  However at this point we can't rule
out the possibility the $\mu$ has other preimages which do not decompose in
this way. We will refer to the unique preimage of the first kind as the
{\it standard lift\/} of the strong shuffle $\mu$. For example, by Lemma \ref{intAss:Lemma}
any irreducible morphism in ${\cal M}_n(k)$ is automatically a standard lift.
\end{remark}

\begin{lemma}\label{stdLift:Lemma}
Suppose the following diagram is given in ${\cal M}_n(k)$ :
$$
\diagram
A \Box_i B \Box_i C \rrto^{f \Box_i id_C}
\dto^{id_A \Box_i g}
&&D \Box_i C\dto^h \\
A \Box_i G \rrto^l &&F
\enddiagram
$$
with $F$ $\Box_r$-reducible, with $\Lambda_k^n(h)$, $\Lambda_k^n(l)$,
$\Lambda_k^n(f)$ and $\Lambda_k^n(g)$
all strong shuffles, and with $h$ and $l$ being standard lifts.
Then the diagram is commutative.
\end{lemma}

\noindent
{\bf Proof:}
Let us first decompose the objects $A$, $B$, $D$ $C$ and $F$ into
$\Box_r$--irreducible objects:
$$
\begin{array}{cc}
A & = \overline A_1 \Box_r \overline A_2 \Box_r \dots \Box_r
\overline A_s \\
B & = \overline B_1 \Box_r \overline B_2 \Box_r \dots \Box_r
\overline B_t \\
D & = \overline D_1 \Box_r \overline D_2 \Box_r \dots \Box_r
\overline D_u \\
C & = \overline C_1 \Box_r \overline C_2 \Box_r \dots \Box_r
\overline C_v \\
F & = \overline F_1 \Box_r \overline F_2 \Box_r \dots \Box_r
\overline F_w
\end{array}
$$
Then a similar argument to the one used in Lemma \ref{extAss:Lemma} gives the
nondecreasing functions
$$
\begin{array}{cc}
\sigma & : \{1,\ 2,...,\ u\} \longrightarrow \{1,\ 2,...,\ w\} \\
\tau & : \{1,\ 2,...,\ v\} \longrightarrow \{1,\ 2,...,\ w\}
\end{array}
$$
defined, respectively, by the relations
$$
\begin{array}{cc}
|\overline D_j| & \subset |\overline F_{\sigma(j)}| \text{, for all}\ j
\in \{1,\ 2,...,\ u\} \\
|\overline C_j| & \subset |\overline F_{\tau(j)}| \text{, for all}\ j
\in \{1,\ 2,...,\ v\}
\end{array}
$$

Since $h$ is the standard lift, according to Lemma \ref{extAss:Lemma} it factors as the composite
$$\diagram
(D_1 \Box_r D_2) \Box_i (C_1 \Box_r C_2)
\rrto^{\eta^{ir}_{D_1, D_2, C_1, C_2}}
&&(D_1 \Box_i C_1) \Box_r (D_2 \Box_i C_2)
\rrto^(.67){h_1 \Box_r h_2} &&F
\enddiagram$$
with
$$
\begin{array}{cc}
D_1 & := \Box_r \{\overline D_j\ |\ j \in \sigma^{-1}(1)\} \\
D_2 & := \Box_r \{\overline D_j\ |\ j \notin \sigma^{-1}(1)\}
\end{array}
\qquad
\qquad
\qquad
\begin{array}{cc}
C_1 & := \Box_r \{\overline C_j\ |\ j \in \tau^{-1}(1)\} \\
C_2 & := \Box_r \{ \overline C_j\ |\ j \notin \tau^{-1}(1)\}
\end{array}
$$

Moreover, the splitting of $D$ as $D_1 \Box_r D_2$
gives the nondecreasing functions
$$
\begin{array}{cc}
\xi & : \{1,\ 2,...,\ s\} \longrightarrow \{1,\ 2\} \\
\zeta & : \{1,\ 2,...,\ t\} \longrightarrow \{1,\ 2\}
\end{array}
$$
defined, respectively, by the relations
$$
\begin{array}{cc}
|\overline A_j| & \subset |D_{\xi(j)}| \text{, for all}\ j \in
\{1,\ 2,...,\ s\} \\
|\overline B_j| & \subset |D_{\zeta(j)}| \text{, for all}\ j \in
\{1,\ 2,...,\ t\}
\end{array}
$$
Therefore by IH.1 $f$ factors as the composite
$$\diagram
(A_1 \Box_r A_2) \Box_i (B_1 \Box_r B_2)
\rrto^{\eta^{ir}_{A_1, A_2, B_1, B_2}}
&&(A_1 \Box_i B_1)\Box_r (A_2 \Box_i B_2)
\rrto^(.63){f_1 \Box_r f_2} &&D_1 \Box_r D_2
\enddiagram$$
with
$$
\begin{array}{cc}
A_1 & := \Box_r \{\overline A_j\ |\ j \in \xi^{-1}(1)\} \\
A_2 & := \Box_r \{\overline A_j\ |\ j \notin \xi^{-1}(1)\}
\end{array}
\qquad
\qquad
\qquad
\begin{array}{cc}
B_1 & := \Box_r \{\overline B_j\ |\ j \in \zeta^{-1}(1)\} \\
B_2 & := \Box_r \{\overline B_j\ |\ j \notin \zeta^{-1}(1)\}
\end{array}
$$
Since $\Lambda_k^n(f)$ and $\Lambda_k^n(h)$ are strong shuffles, we have that
$$A_1=A\cap |\overline F_1|,\quad B_1=B\cap |\overline F_1|,\quad C_1 = C\cap |\overline F_1|$$
and that
$$A_2=A\cap |F_2|,\quad B_2=B\cap |F_2|,\quad C_2=C\cap |F_2|$$
where $F_2=\overline F_2\Box_r \overline F_3\Box_r\dots\Box_r\overline F_w$.

Similar arguments give decompositions
$$\diagram
(A_1\Box_r A_2)\Box_i(G_1\Box_r G_2)\rrto^{\eta^{ir}_{A_1,A_2,G_1,G_2}}
&&(A_1\Box_i G_1)\Box_r(A_2\Box_i G_2)\rrto^(.65){l_1\Box_r l_2}
&&F
\enddiagram$$
and
$$\diagram
 (B_1\Box_i B_2)\Box_i (C_1\Box_i C_2)\rrto^{\eta^{ir}_{B_1,C_1,B_2,C_2}}
&&(B_1\Box_i C_1)\Box_r(B_2\Box_i C_2)\rrto^(.65){g_1\Box_r g_2}&&G_1\Box_r G_2
\enddiagram$$
of $l$, respectively $g$.

Thus we obtain the following diagram in ${\cal M}_n(k)$:
$$\spreaddiagramrows{-.8pc}
\spreaddiagramcolumns{-.8pc}
\def\objectstyle{\scriptstyle}
\def\labelstyle{\scriptstyle}
\diagram
(A_1\Box_r A_2)\Box_i (B_1\Box_i B_2)\Box_i (C_1\Box_i C_2)
\rrto^{\eta^{ir}\Box_i id_C}\ddto^{id_A\Box_i\eta^{ir}}
&&((A_1\Box_i B_1)\Box_r (A_2\Box_i B_2))\Box_i (C_1\Box_i C_2)
\xto[rrrr]^(.55){(f_1\Box_r f_2)\Box_i id_C}\ddto^{\eta^{ir}}
&&&&(D_1\Box_r D_2)\Box_i (C_1\Box_i C_2)\ddto_{\eta^{ir}_{D_1,D_2,C_1,C_2}}\\
&{\squarify{1}\circled}&&&{\squarify{2}\circled}\\
(A_1\Box_r A_2)\Box_i( (B_1\Box_i C_1)\Box_r(B_2\Box_i C_2))
\rrto^{\eta^{ir}}\ddto^{id_A\Box_i(g_1\Box_r g_2)}
&&(A_1\Box_i B_1\Box_i C_1)\Box_r(A_2\Box_i B_2\Box_i C_2)
\xto[rrrr]^(.55){(f_1\Box_i id_{C_1})\Box_r(f_2\Box_i id_{C_2})}
\ddto^{(id_{A_1}\Box_i g_1)\Box_r(id_{A_2}\Box_i g_2)}
&&&&(D_1\Box_i C_1)\Box_r(D_2\Box_i C_2)\ddto_{h_1\Box_r h_2}\\
&{\squarify{3}\circled}&&&{\squarify{4}\circled}\\
(A_1\Box_r A_2)\Box_i(G_1\Box_r G_2)\rrto^{\eta^{ir}_{A_1,A_2,G_1,G_2}}
&&(A_1\Box_i G_1)\Box_r(A_2\Box_i G_2)\xto[rrrr]^{l_1\Box_r l_2}
&&&&F
\enddiagram$$
The outer square of this diagram is the original diagram we want to show
commutes.  This follows from the fact that all the inner subdiagrams commute:
(1) by the internal associativity diagram, (2) and (3) by naturality of $\eta^{ir}$,
and (4) by Remark \ref{compatible:splitting}.
\hfill
$\Box$

\begin{lemma}\label{extending:diagram}
Suppose the following diagram is given in $\widehat{{\cal M}}_n(k)$ :
$$
\diagram
&B \Box_i C\dto^{\varphi}\\
D \Box_j F \rto^{\psi} &G
\enddiagram
$$
with both morphisms strong shuffles and $i \leq j$. Let $X_s$, $s =
1,\ 2,\ 3,\ 4$, be the objects defined, respectively, by
$$
\begin{array}{cc}
X_1 & = D - |C| \\
X_2 & = F - |C|
\end{array}
\qquad
\qquad
\qquad
\begin{array}{cc}
X_3 & = D - |B| \\
X_4 & = F - |B|
\end{array}
$$

Then :
\begin{enumerate}
\item There exist two morphisms
$$
\begin{array}{c}
(X_1 \Box_j X_2) \Box_i (X_3 \Box_j X_4)
\stackrel{\varphi_1\Box_i \varphi_2}{\longrightarrow} B \Box_i C \\
(X_1 \Box_i X_3) \Box_j (X_2 \Box_i X_4)
\stackrel{\psi_1\Box_j \psi_2}{\longrightarrow} D \Box_j F
\end{array}
$$
in $\widehat{{\cal M}}_n(k)$ extending the given diagram to
$$
\diagram
&(X_1 \Box_j X_2) \Box_i (X_3 \Box_j X_4)\dto^{\overline{\varphi}}\\
(X_1 \Box_i X_3) \Box_j (X_2 \Box_i X_4)
\rto^(.66){\overline{\psi}} &G
\enddiagram
$$

\item The extended diagram can be completed into a commutative
triangle whenever either one of the following conditions are
satisfied :
\begin{enumerate}
\item $i \neq j$;

\item $i = j$ and at least one of the objects $X_2$, $X_3$ is
equal to $0$.
\end{enumerate}
\end{enumerate}
\end{lemma}

\noindent
{\bf Proof:}
Note first that (b) follows immediately from (a). Indeed, the
morphism
$$
\diagram
(X_1 \Box_j X_2) \Box_i (X_3 \Box_j X_4)
\xto[rrr]^{\Lambda^n_k(\eta^{ij}_{X_1, X_2, X_3, X_4})}
&&&(X_1 \Box_i X_3) \Box_j (X_2 \Box_i X_4)
\enddiagram
$$
has the required property if $i \neq j$, while the morphism
$id_{X_1 \Box_i X_2 \Box_i X_4}$ is taking care of the case
$i = j$ (assuming, without loss of generality, $X_3 = 0$). So all we
have to prove is (1).

The condition that both $\varphi$ and $\psi$ are strong shuffles
yields the existence of the following morphisms :
$$
\begin{array}{c}
 (D - |C|) \Box_j (F - |C|) = (D \Box_j F) - |C| 
 \stackrel{\varphi_1}{\longrightarrow} G - |C| = B \\
(D - |B|) \Box_j (F - |B|) = (D \Box_j F) - |B| 
\stackrel{\varphi_2}{\longrightarrow} G - |B| = C \\
(B - |F|) \Box_i (C - |F|) = (B \Box_i C) - |F|
\stackrel{\psi_1}{\longrightarrow} G - |F| = D \\
(B - |D|) \Box_i (C - |D|) = (B \Box_i C) - |D|
\stackrel{\psi_2}{\longrightarrow} G - |D| = F
\end{array}
$$
and therefore the only thing still to prove is the following set of
equalities :
$$
\begin{array}{c}
D - |C| = B - |F|; \\
F - |C| = B - |D|;
\end{array}
\qquad
\qquad
\qquad
\qquad
\begin{array}{c}
D - |B| = C - |F|; \\
F - |B| = C - |D|.
\end{array}
$$
But this can be easily done -- by using again the fact that both
$\varphi$ and $\psi$ are strong shuffles -- as follows :
$$
\begin{array}{cl}
D - |C| & = (D - |F|) - |C| = ((D \Box_j F) - |F|) - |C| = (G - |F|) - |C| = \\
& = (G - |C|) - |F| = ((B \Box_i C) - |C|) - |F| = (B - |C|) - |F| = B - |F|; \\
F - |C| & = (F - |D|) - |C| = ((D \Box_j F) - |D|) - |C| = (G - |D|) - |C| = \\
& = (G - |C|) - |D| = ((B \Box_i C) - |C|) - |D| = (B - |C|) - |D| = B - |D|; \\
D - |B| & = (D - |F|) - |B| = ((D \Box_j F) - |F|) - |B| = (G - |F|) - |B| = \\
& = (G - |B|) - |F| = ((B \Box_i C) - |B|) - |F| = (C - |B|) - |F| = C - |F|; \\
F - |B| & = (F - |D|) - |B| = ((D \Box_j F) - |D|) - |B| = (G - |D|) - |B| = \\
& = (G - |B|) - |D| = ((B \Box_i C) - |B|) - |D| = (C - |B|) - |D| = C - |D|.
\end{array}
$$
and the proof is completed.
\hfill
$\Box$

\begin{lemma}\label{giantHex:Lemma}
Suppose the following diagram is given in ${\cal M}_n(k)$ :
$$
\diagram
A \rto\dto &B \Box_i C\dto^h \\
D \Box_j F \rto^l &G
\enddiagram
$$
with $G$ $\Box_r$-reducible and with $i < j < r$. If $\Lambda_k^n(h)$ and $\Lambda_k^n(l)$ are both
strong shuffles and $h$ and $l$ are standard lifts (cf. Remark \ref{standard:lift}),
then the diagram is commutative.
\end{lemma}

\noindent
{\bf Proof:}
Let $G_1$ and $G_2$ be the objects defined (uniquely) by the
equality $G = G_1 \Box_r G_2$ and the condition that $G_1$ is
$\Box_r$--irreducible. By Lemma \ref{extAss:Lemma} we can replace the original given
decompositions of $h$ and $l$ by new ones compatible with this
splitting of $G$:
$$
\diagram
B \Box_i C = (X_1 \Box_r X_2) \Box_i (X_3 \Box_r X_4)
\rrto^(.55){\eta^{ir}_{X_1, X_2, X_3, X_4}}
&&(X_1 \Box_i X_3) \Box_r (X_2 \Box_i X_4)
\rrto^(.65){h_1 \Box_r h_2} &&G \\
D \Box_j F = (Z_1 \Box_r Z_2) \Box_j (Z_3 \Box_r Z_4)
\rrto^(.55){\eta^{ir}_{Z_1, Z_2, Z_3, Z_4}}
&&(Z_1 \Box_j Z_3) \Box_r (Z_2 \Box_j Z_4)
\rrto^(.65){l_1 \Box_r l_2} &&G.
\enddiagram
$$

Then the morphisms $\varphi_1 := \Lambda_{|G_1|}^n(h_1)$,
$\varphi_2 := \Lambda_{|G_2|}^n(h_2)$,
$\psi_1 := \Lambda_{|G_1|}^n(l_1)$ and
$\psi_2 := \Lambda_{|G_2|}^n(l_2)$ are strong shuffles and we
are within the hypotheses of Lemma \ref{stdLift:Lemma} with the following two
diagrams (in $\widehat{{\cal M}_n}(|G_1|)$ and
$\widehat{{\cal M}_n}(|G_2|)$, respectively) :
$$
\diagram
&X_1 \Box_i X_3\dto^{\varphi_1} \\
Z_1 \Box_j Z_3 \rto^{\psi_1} &G_1
\enddiagram
\qquad \qquad \qquad
\diagram
&X_2 \Box_i X_4\dto^{\varphi_2}\\
Z_2 \Box_j V_4 \rto^{\psi_2} &G_2
\enddiagram
$$
Therefore there exist the objects $Y_u,\ u = 1,\ 2,\ ...,\ 8$
together with the following morphisms (in the corresponding
components of $\widehat M_n$) :
$$
\begin{array}{l}
Y_1 \Box_j Y_3 \stackrel{\xi_1}{\longrightarrow} X_1 \\
Y_5 \Box_j Y_7 \stackrel{\xi_3}{\longrightarrow} X_3 \\
Y_1 \Box_i Y_5 \stackrel{\zeta_1}{\longrightarrow} Z_1 \\
Y_3 \Box_i Y_7 \stackrel{\zeta_3}{\longrightarrow} Z_3
\end{array}
\qquad \qquad \qquad
\begin{array}{l}
Y_2 \Box_j Y_4 \stackrel{\xi_2}{\longrightarrow} X_2 \\
Y_6 \Box_j Y_8 \stackrel{\xi_4}{\longrightarrow} X_4 \\
Y_2 \Box_i Y_6 \stackrel{\zeta_2}{\longrightarrow} Z_2 \\
Y_4 \Box_i Y_8 \stackrel{\zeta_4}{\longrightarrow} Z_4
\end{array}
$$

which -- according to (IH.1) -- can be lifted, respectively, to
$$
\begin{array}{l}
Y_1 \Box_j Y_3 \stackrel{f_1}{\longrightarrow} X_1 \\
Y_5 \Box_j Y_7 \stackrel{f_3}{\longrightarrow} X_3 \\
Y_1 \Box_i Y_5 \stackrel{g_1}{\longrightarrow} Z_1 \\
Y_3 \Box_i Y_7 \stackrel{g_3}{\longrightarrow} Z_3
\end{array}
\qquad \qquad \qquad
\begin{array}{l}
Y_2 \Box_j Y_4 \stackrel{f_2}{\longrightarrow} X_2 \\
Y_6 \Box_j Y_8 \stackrel{f_4}{\longrightarrow} X_4 \\
Y_2 \Box_i Y_6 \stackrel{g_2}{\longrightarrow} Z_2 \\
Y_4 \Box_i Y_8 \stackrel{g_4}{\longrightarrow} Z_4
\end{array}
$$
in the corresponding components of ${\cal M}_n$ since the cardinalities
of all the targets are smaller than $k$.

Also note that there exists a unique morphism $u : A \longrightarrow (Y_{12}^r
\Box_j Y_{34}^r) \Box_i (Y_{56}^r \Box_j Y_{78}^r)$ in
${\cal M}_n(k)$, since such a morphism exists in
$\widehat{{\cal M}_n}(k)$ and its target is $\Box_i$--split, with
$Y_{12}^r$ denoting the object $Y_1 \Box_r Y_2$ and so on.

This gives rise to the following diagram in ${\cal M}_n(k)$:
$$
\spreaddiagramrows{-.8pc}
\spreaddiagramcolumns{-.8pc}
\def\objectstyle{\scriptstyle}
\def\labelstyle{\scriptstyle}
\diagram
&&{B\Box_iC}\xto[rrrrrr]^{\eta^{ir}}
&&&&&&{X_{13}^i\Box_rX_{24}^i}\xto[ddddrr]^{h_1\Box_r h_2}\\
&&&&&{\squarify{2}\circled}\\
&&{\squarify{1}\circled}&&{Y_{13245768}^{jri}}\xto[uull]_f\xto[rr]^{\eta^{ir}}
&&{Y_{13572468}^{jir}\xto[uurr]^{\hat{f}}}\xto[ddrr]\\
\\
{A}\rrto^u\xto[uuuurr]\xto[ddddrr] &&Y_{12345678}^{rji}\uurrto\ddrrto
&&&{\squarify{6}\circled}
&&&{Y_{15372648}^{ijr}}\xto[dddd]_{\hat{g}}&{\squarify{3}\circled}&{G}\\
\\
&&{\squarify{5}\circled}&&{Y_{12563478}^{rij}}\rrto
&&{Y_{15263748}^{irj}}\uurrto^{\eta^{jr}}\xto[ddllll]^g\\
&&&&&&{\squarify{4}\circled}\\
&&{D\Box_jF\xto[rrrrrr]^{\eta^{jr}}}
&&&&&&{Z_{13}^j\Box_rZ_{24}^j}\xto[uuuurr]_{l_1\Box_r l_2}
\enddiagram
$$

Here we denote
$$Y_{abcdxyzw}^{pqs}:=(((Y_a\Box_p Y_b)\Box_q(Y_c\Box_p Y_d))\Box_s
((Y_x\Box_p Y_y)\Box_q(Y_z\Box_p Y_w)))$$
$$X_{ab}^i := X_a\Box_i Y_b \qquad\qquad\qquad Z_{ab}^j := Z_a\Box_j Z_b$$
$$f:=(f_1\Box_r f_2)\Box_i(f_3\Box_r f_4) \qquad\qquad
\hat{f}:=(f_1\Box_i f_3)\Box_r(f_2\Box_i f_4)$$
$$g:=(g_1\Box_r g_2)\Box_j(g_3\Box_r g_4)\qquad\qquad
\hat{g}:=(g_1\Box_j g_3)\Box_r(g_2\Box_j g_4)$$

Then the outer hexagon is an expansion of the original diagram which
we want to show commutes.  To show this we observe that all the inner
subdiagrams commute: diagrams (1) and (5) by IH.2, diagrams (2) and (4)
by naturality of $\eta^{ir}$ and $\eta^{jr}$ respectively, diagram (6)
is the ``Giant Hexagon'', and diagram (3) by Remark \ref{compatible:splitting}.

This concludes the proof. \hfill
$\Box$

\begin{lemma}\label{giantHex:too}
Let the following diagram be given in ${\cal M}_n(k)$ :
$$\diagram
A_1 \Box_i A_2
\xto[rrr]^{\eta^{ij}_{A'_{11}, A'_{12}, A'_{21}, A'_{22}}}
\dto^{\eta^{ir}_{A''_{11}, A''_{12}, A''_{21}, A''_{22}}}
&&&B_1 \Box_j B_2 \dto^g \\
C_1 \Box_r C_2 \xto[rrr]^{h} &&&D_1 \Box_r D_2
\enddiagram$$
with $i < j < r$,
$\eta^{ij}_{A'_{11}, A'_{12}, A'_{21}, A'_{22}}$ and
$\eta^{ir}_{A''_{11}, A''_{12}, A''_{21}, A''_{22}}$ nontrivial
and
$$
\begin{array}{cl}
A_1 & = A'_{11} \Box_j A'_{12} = A''_{11} \Box_r A''_{12}\\
A_2 & = A'_{21} \Box_j A'_{22} = A''_{21} \Box_r A''_{22}
\end{array}
$$
If  
$\Lambda_k^n(g\eta^{ij}_{A'_{11}, A'_{12}, A'_{21}, A'_{22}})
=\Lambda_k^n(h\eta^{ir}_{A''_{11}, A''_{12}, A''_{21}, A''_{22}})$
is a strong shuffle then the diagram is commutative.
\end{lemma}

\noindent
{\bf Proof:} Rewrite $g:B_1\Box_j B_2\longrightarrow D_1\Box_r D_2$
in the form
$$B_1\Box_j B_2\stackrel{g'}{\longrightarrow}B'_1\Box_{j'} B'_2
\stackrel{g''}{\longrightarrow}D_1\Box_r D_2$$
where $g''$ is irreducible. Then the result follows by applying
Lemma \ref{giantHex:Lemma} to the diagram:
$$\diagram
A_1\Box_i A_2\xto[rrr]^{g'\eta^{ij}_{A'_{11}, A'_{12}, A'_{21}, A'_{22}}}
\dto^{id_{A_1\Box_i A_2}}
&&&B'_1\Box_{j'} B'_2\dto^{g''}\\
A_1\Box_i A_2\xto[rrr]^{h\eta^{ir}_{A''_{11}, A''_{12}, A''_{21}, A''_{22}}}
&&&D_1\Box_r D_2
\enddiagram$$
\hfill
$\Box$

\begin{lemma} $f:A\Box_i B\longrightarrow C$ is irreducible in ${\cal M}_n(k)$
iff $\Lambda_k^n(f):A\Box_i B\longrightarrow C$ is irreducible in
$\widehat{{\cal M}}_n(k)$
\end{lemma}

\noindent
{\bf Proof:} By (IH.2) we may as well assume that $C$ is $\Box_r$-reducible.
As noted in Remark \ref{clearly:irreducible}, the implication
$$\Lambda_k^n(f)\ \mbox{irreducible}\quad\Longrightarrow\quad
f\ \mbox{irreducible}$$
is trivially true.

Now suppose $f$ is irreducible. Then by Lemma \ref{intAss:Lemma}, $\Lambda_k^n(f)$ is a strong
shuffle and $A$ and $B$ are both $\Box_i$--irreducible. Thus we can't have
a nontrivial factorization of $\Lambda_k^n(f)$ of the form
$$A\Box_i B\longrightarrow D\Box_i G\longrightarrow C$$
For if $\mbox{card}(|D|)=\mbox{card}(|A|)$, then this contradicts $\Lambda_k^n(f)$ being a strong shuffle.
If $\mbox{card}(|D|)<\mbox{card}(|A|)$, then this contradicts $\Lambda_k^n(f)$ being a strong shuffle and
$A$ being $\Box_i$--irreducible (cf. Proposition \ref{coherence:basicObservation}(4) and Remark 
\ref{basicObservation:remark}).  Similarly we can rule out $\mbox{card}(|D|)>\mbox{card}(|A|)$.

Thus if $\Lambda_k^n(f)$ were not irreducible in $\widehat{{\cal M}}_n(k)$, then
there would have to be a factorization of $\Lambda_k^n(f)$ of the form
$$A\Box_i B\stackrel{\mu}{\longrightarrow} D\Box_j G
\stackrel{\varphi}{\longrightarrow} C$$
with $i<j<r$ and $\varphi$ irreducible.  Then by IH.2 we can lift $\mu$ to
a morphism $h$ and by Lemma \ref{irreducible:uniqueLift} we can lift $\varphi$ to an irreducible
morphism $l$. But then by Lemma \ref{giantHex:Lemma} we have the following commutative
diagram in ${\cal M}_n(k)$
$$\diagram
A\Box_i B\rrto^h\dto^{id_{A\Box_i B}}
&&D\Box_j G\dto^l\\
A\Box_i B\rrto^f &&C
\enddiagram$$
contradicting the irreducibility of $f$. Thus $\Lambda_k^n(f)$ must be irreducible.
\hfill
$\Box$

\begin{lemma}\label{coherence:hardLemma}
Suppose the following diagram is given in ${\cal M}_n(k)$ :
$$
\diagram
A \rto^f\dto^g &B \Box_i C\dto^h \\
D \Box_i F \rto^l &G
\enddiagram
$$
with $G$ $\Box_r$-reducible.  If $h$ and $l$ are both irreducible then the diagram is commutative.
\end{lemma}

\noindent
{\bf Proof:}
Note first that the given diagram can be projected in
$\widehat{{\cal M}}_n(k)$ via the functor $\Lambda_k^n$, the
result being the commutative diagram
\begin{equation}
{\diagram
A \rto\dto &B \Box_i C\dto^{\varphi}\\
D \Box_i F \rto^{\psi} &G
\enddiagram}
\end{equation}
with $\varphi = \Lambda_k^n(h)$ and $\psi = \Lambda_k^n(l)$.
According to Lemma \ref{intAss:Lemma}, $\varphi$ and $\psi$ are strong shuffles in
$\widehat{{\cal M}}_n(k)$. Therefore the lower right--hand side corner
of $(1)$ is exactly the diagram in Lemma \ref{extending:diagram} with $i = j$.

\underbar{\it Case 1.}
Suppose the additional hypothesis in Lemma \ref{extending:diagram}(b) is satisfied in
our situation, namely one of the objects $X_2$, $X_3$ is equal to
$0$. Without loss of generality we can assume $X_3 = 0$. Then the
extended diagram in Lemma \ref{extending:diagram} can be written as
\begin{equation}{
\diagram
X_1 \Box_i X_2 \Box_i X_4
\rrto^{\varphi_1 \Box_i \varphi_2}
\dto^{\psi_1 \Box_i \psi_2}
&&B \Box_i C\dto^{\varphi} \\
D \Box_i F \rrto^{\psi} &&G
\enddiagram}
\end{equation}

Next, the fact that $B$ and $C$, on one hand, and $D$ and $F$, on the
other hand, have no common generating objects yields the following
equivalences :
$$
\begin{array}{l}
X_3 = 0 \Longleftrightarrow C - |F| = 0 \Longleftrightarrow F - |B| = C
\Longleftrightarrow X_4 = C \\
X_3 = 0 \Longleftrightarrow D - |B| = 0 \Longleftrightarrow B - |F| = D
\Longleftrightarrow X_1 = D
\end{array}
$$
together with the equalities $\varphi_2 = id_C$ and $\psi_1 = id_D$
(in $\widehat{{\cal M}}_n(k)$). Therefore Lemma \ref{irreducible:uniqueLift} and (IH.1) give the
following (unique) lift of $(2)$ in ${\cal M}_n(k)$ :
$$
\diagram
D \Box_i X_2 \Box_i C \rrto^{f_1 \Box_i id_C}
\dto^{id_D \Box_i g_2} &&B \Box_i C\dto^h \\
D \Box_i F \rrto^l &&G
\enddiagram
$$
satisfying the hypotheses in Lemma \ref{stdLift:Lemma}; hence it is commutative :
\begin{equation}
h \circ (f_1 \Box_i id_C) = l \circ (id_D \Box_i g_2)
\end{equation}

Finally, there exists a morphism
$$
\xi : A \longrightarrow D \Box_i X_2 \Box_i C 
$$
in $\widehat{{\cal M}}_n(k)$. According to (IH.2), $\xi$ has a unique
lift $a$ in ${\cal M}_n(k)$ and the morphisms $f$, $g$ factor,
respectively, as
\begin{equation}
{\diagram
A \rto^(.3)a &D \Box_i X_2 \Box_i C \rrto^{f_1 \Box_i id_C}
&&B \Box_i C \\
A \rto^(.3)a &D \Box_i X_2 \Box_i C \rrto^{id_D \Box_i g_2}
&&D \Box_i F
\enddiagram}
\end{equation}
Now the conclusion follows immediately from $(3)$ and $(4)$.

\underbar{\it Case 2.}
Let us assume now that both $X_2$ and $X_3$ are different from $0$.
In this situation the extended diagram in Lemma \ref{extending:diagram} cannot be closed
to a commutative triangle. Nevertheless,  we can consider the
objects $B\cap D:=B\cap |D|=D\cap |B|$ (since both $B$ and $D$ are restrictions of the
object $G$ as both $\varphi$ and $\psi$ are strong shuffles) and
$C\cap F:= C\cap |F|=F\cap |C|$.

{\it Subcase 2.1.}
Suppose that at least one of the objects $B \cap D$ and $C \cap F$ is
not equal to $0$. (Without loss of generality we may consider
$Y := B \cap D \neq 0$.)

Then the morphism $\varphi \circ \Lambda_k^n(f) = \psi \circ
\Lambda_k^n(g)$ factors (in $\widehat{{\cal M}}_n(k)$) through
the object $Y \Box_i Z$, with $Z := G - |Y|$. Suppose this
factorization is
$$
A \stackrel{\xi}{\longrightarrow} Y \Box_i Z \stackrel{\mu}{\longrightarrow} G
$$
Then $\xi$ has a unique lift $a$ in ${\cal M}_n(k)$, according to
(IH.2),while $\mu$ is a strong shuffle, by the definition of the
objects $Y$ and $Z$, and therefore it has a unique standard lift $b$ in
${\cal M}_n(k)$, according to Lemma \ref{extAss:Lemma} and Remark \ref{standard:lift}.

Then we have the following diagrams in ${\cal M}_n(k)$ which can be shown to
commute by the same argument as in Case 1 :
$$
\diagram
A \rto^f\dto^a &B \Box_i C\dto^h \\
Y \Box_i Z \rto^b &G
\enddiagram
\qquad
\qquad
\qquad
\qquad
\diagram
A \rto^g\dto^a &D \Box_i F\dto^l \\
Y \Box_i Z \rto^b &G
\enddiagram
$$

{\it Subcase 2.2.}
The remaining situation is $B \cap D = C \cap F = 0$.  In this case we must have $|B|=|F|$ and $|C|=|D|$. 
Since all the objects $B$, $C$, $D$ and $F$ are restrictions of the object $G$ (because $\varphi$
and $\psi$ are strong shuffles), we must have $B=F$ and $C=D$. Therefore the given diagram can be
written as
\begin{equation}
{\diagram
A \rto^f\dto^g &B \Box_i C\dto^h \\
C \Box_i B \rto^l &G
\enddiagram}
\end{equation}

A closer look at the morphism $f$ shows that -- according to (IH.2)
-- it factors through a certain object
\begin{equation}
Z := \overline B_1 \Box_{i-1} \overline C_1 \Box_{i-1}
\overline B_2 \Box_{i-1} \overline C_2 \Box_{i-1} \dots
\Box_{i-1} \overline B_m \Box_{i-1} \overline C_m
\end{equation}
with $m$ a positive integer and the objects $\overline B_t$,
$\overline C_t (t \in \{1,\ 2,\ ...,\ m\}$) given by the procedure
described below.

Consider the following subsets of $\{1,\ 2,\ ...,\ k\}$ :
$$
\begin{array}{cl}
{\cal B}_1 & := \{b \in |B|\ |\ \forall c \in |C|,\ \exists s < i
\mbox { such that } b\Box_s c\mbox{ in }A\} \\
{\cal C}_1 & := \{c \in |C|\ |\ \forall b \in |B| \setminus
{\cal B}_1,\ \exists s < i \mbox { such that } c\Box_s b\mbox{ in }A\} \\
{\cal B}_2 & := \{b \in |B| \setminus {\cal B}_1\ |\ \forall c \in
|C| \setminus {\cal C}_1,\ \exists s < i \mbox { such that } b\Box_s c\mbox{ in }A\} \\
{\cal C}_2 & := \{c \in |C| \setminus {\cal C}_1\ |\ \forall b \in
|B| \setminus ({\cal B}_1 \cup {\cal B}_2),\ \exists s < i \mbox
{ such that } c\Box_s b\mbox{ in }A\} \\
\dots & \\
{\cal B}_{m-1} & := \{b \in |B| \setminus (\bigcup_{t=1}^{m-2}
{\cal B}_t)\ |\ \forall c \in |C| \setminus (\bigcup_{t=1}^{m-2}
{\cal C}_t),\ \exists s < i \mbox { such that } b\Box_s c\mbox{ in }A\} \\
{\cal C}_{m-1} & := \{c \in |C| \setminus (\bigcup_{t=1}^{m-2}
{\cal C}_t)\ |\ \forall b \in |B| \setminus (\bigcup_{t=1}^{m-1}
{\cal B}_t),\ \exists s < i \mbox { such that } c\Box_s b\mbox{ in }A\} \\
{\cal B}_m & := |B| \setminus (\bigcup_{t=1}^{m-1} {\cal B}_t) \\
{\cal C}_m & := |C| \setminus (\bigcup_{t=1}^{m-1} {\cal C}_t)
\end{array}
$$
Then the objects $\overline B_t$ ($t \in \{1,\ 2,\ ...,\ m\}$) are
defined as the results obtained by deleting in $B$ all the generating
objects from $|B| \setminus {\cal B}_t$, respectively. A similar
definition gives the objects $\overline C_t$. Note that {\it only}
$\overline B_1$ or $\overline C_m$ or both can be equal to $0$.
Moreover, if $\overline B_1 = 0$ then $m \geq 2$.

{\it Subcase 2.2.1.}
If {\it only two} of the objects $\overline B_t$, $\overline C_t$ in
the right--hand side of $(6)$ are different from $0$ then $Z$ is
given by one of the equalities
$$
\begin{array}{cl}
Z & = \overline B_1 \Box_{i-1} \overline C_1 \\
Z & = \overline C_1 \Box_{i-1} \overline B_2
\end{array}
$$
and, in order to make a choice, we shall assume {\it the first one}
to hold (the other situation being treated in a similar way). In this
case we obviously have $B_1 = B$ and $C_1 = C$.

According to (IH.2), there exists a unique morphism $a : A \longrightarrow
B \Box_{i-1} C$ in ${\cal M}_n(k)$. Then, again by (IH.2), both
$f$ and $g$ factor through $B \Box_{i-1} C$ as
$$
\begin{array}{cl}
f & = \eta^{ji}_{B, 0, 0, C} \circ a \\
g & = \eta^{ji}_{0, B, C, 0} \circ a
\end{array}
$$
with $j := i - 1$ and the commutativity of $(7)$ is obviously
reduced to the commutativity of the following diagram :
\begin{equation}
{\diagram
B \Box_j C \rrto^{\eta^{ji}_{B, 0, 0, C}}
\dto^{\eta^{ji}_{0, B, C, 0}} &&B \Box_i C\dto^h \\
C \Box_i B \rrto^l &&G
\enddiagram}
\end{equation}

Next, let us have a closer look at the irreducible morphism $h$.
According to Lemma \ref{extAss:Lemma}, the morphism $h$ can be factored as
$$
\diagram
(B_1 \Box_r B_2) \Box_i (C_1 \Box_r C_2)
\rrto^{\eta^{ir}_{B_1, B_2, C_1, C_2}}
&&(B_1 \Box_i C_1) \Box_r (B_2 \Box_i C_2)
\dto^{h_1 \Box_r h_2}\\
&&G_1 \Box_r G_2
\enddiagram
$$
for any decomposition $G_1 \Box_r G_2$ of the
$\Box_r$--reducible object $G$. But this fact implies the
existence of a $\Box_r$--split morphism
$$
\mu_1 \Box_r \mu_2 : (B_1 \Box_j C_1) \Box_r
(B_2 \Box_j C_2) \longrightarrow G_1 \Box_r G_2
$$
in $\widehat{{\cal M}}_n(k)$ which, according to (IH.1), has a unique
lift $h'_1 \Box_r h'_2$ in ${\cal M}_n(k)$. Hence we have obtained
the following two diagrams in ${\cal M}_n(k)$ :
$$
\diagram
B \Box_j C \rrto^{\eta^{ji}_{B, 0, 0, C}}
\dto^{\eta^{jr}_{B_1, B_2, C_1, C_2}}
&&B \Box_i C\dto^h \\
D_1 \Box_r D_2 \rrto^{h'_1 \Box_r h'_2} &&G_1 \Box_r G_2
\enddiagram
\qquad
\qquad
\diagram
B \Box_j C \rrto^{\eta^{ji}_{0, B, C, 0}}
\dto^{\eta^{jr}_{B_1, B_2, C_1, C_2}}
&&C \Box_i B\dto^l \\
D_1 \Box_r D_2 \rrto^{h'_1 \Box_r h'_2} &&G_1 \Box_r G_2
\enddiagram
$$
with $D_1 \Box_r D_2$ denoting the object
$(B_1 \Box_j C_1) \Box_r (B_2 \Box_j C_2)$. Now the
conclusion follows easily by applying Lemma \ref{giantHex:too}.

{\it Subcase 2.2.2.}
If {\it at least three} of the objects $\overline B_t$,
$\overline C_t$ in the right--hand side of $(8)$ are different
from $0$ then $Z$ is given by an equality having one of the forms
$$
\begin{array}{cl}
Z & = \overline B_1 \Box_{i-1} \overline C_1 \Box_{i-1} V \\
Z & = \overline C_1 \Box_{i-1} \overline B_2 \Box_{i-1} V
\end{array}
$$
corresponding to $\overline B_1 \ne 0$ and $\overline B_1 = 0$,
respectively. Again we can assume {\it the first} equality to hold.

It follows that we have the following diagram in $\widehat{{\cal M}}_n(k)$
$$A\stackrel{\mu}{\longrightarrow}\overline B_1\Box_{i-1} (G-|\overline B_1|)\stackrel{\nu}{\longrightarrow}
G.$$
Next factor $\nu$ as
$$
\overline B_1\Box_{i-1} (G-|\overline B_1|) \stackrel{\nu_1}{\longrightarrow} X \Box_j Y
\stackrel{\nu_2}{\longrightarrow} G
$$
with $ j < r$, $\nu_2$ irreducible and $X \Box_j Y$ different from
both $B \Box_i C$ and $C \Box_i B$ (since there are no
morphisms from $\overline B_1\Box_{i-1} (G-|\overline B_1|)$ into either $B \Box_i C$ or
$C \Box_i B$ in $\widehat{{\cal M}}_n(k)$).
Now use (IH.2) to lift $\mu$ and $\nu_1$ and Lemma \ref{irreducible:uniqueLift} to 
lift $\nu_2$.  We denote the lifts by $u$, $v_1$ and $v_2$ respectively.
This gives us a morphism $A \longrightarrow G$ in ${\cal M}_n(k)$ given by
$$
A \stackrel{a}{\longrightarrow} X \Box_j Y \stackrel{b}{\longrightarrow} G
$$
where $a=v_1u$ and $b=v_2$.

Finally, let us consider the diagrams
$$
\diagram
A \rto^f\dto^a &B \Box_i C\dto^h \\
X \Box_j Y \rto^b &G
\enddiagram
\qquad
\qquad
\qquad
\qquad
\diagram
A \rto^g\dto^a &C \Box_i B\dto^l\\
X \Box_j Y \rto^b &G
\enddiagram
$$
which are commutative, either by Lemma \ref{giantHex:Lemma} (for $i \ne j$) or by one of
the cases already discussed during this proof (for $i = j$). Now the
conclusion is immediate and the lemma is completely proven.
\hfill
$\Box$

Finally, we have all the necessary preliminaries for the proof of
Theorem \ref{coherenceTheorem:too}.

\bigskip

\noindent
{\bf Proof of the Coherence Theorem for $n$-fold Monoidal
Categories.}
It remains to show that
$$\Lambda_k^n:\mbox{Hom}_{{\cal M}_n(k)}(A,B)\longrightarrow
\mbox{Hom}_{\widehat{{\cal M}}_n(k)}(A,B)$$
is a bijection when $A$ is $\Box_r$--irreducible and $B$ is $\Box_r$--reducible,
since (IH.2) and Remark \ref{compatible:splitting} take care of all the other possibilities.

Note first that any morphism $\mu : A \longrightarrow B$ in
$\widehat{{\cal M}}_n(k)$ with $A$ $\Box_r$--irreducible and $B$
$\Box_r$--reducible can be factored as
$$
A \stackrel{\mu'}{\longrightarrow} A'_1 \Box_i A'_2
\stackrel{\mu_0}{\longrightarrow} B
$$
with $\mu_0$ irreducible in $\widehat{{\cal M}}_n(k)$. Since both
$\mu'$ and $\mu_0$ have lifts in ${\cal M}_n(k)$, the former by (IH.2)
and the latter by Lemma \ref{irreducible:uniqueLift}, it follows that the morphism $\mu$ has
such a lift. Therefore the functor $\Lambda_k^n$ is {\it surjective}
on morphisms.

Next let us consider two morphisms $f, g : A \longrightarrow B$ in
${\cal M}_n(k)$, with $A$ $\Box_r$--irreducible and $B$
$\Box_r$--reducible. Obviously $f$ and $g$ can be factored,
respectively, as
$$
\diagram
A \rto^(.4){f_0} &A'_1 \Box_i A'_2 \rto^(.6){f'} &B \\
A \rto^(.4){g_0} &A''_1 \Box_j A''_2 \rto^(.6){g''} &B
\enddiagram
$$
with $f'$ and $g''$ irreducible. But in this way we have obtained in
fact the following diagram in ${\cal M}_n(k)$ :
$$
\diagram
A \rto^(.4){f_0}\dto^{g_0} &A'_1 \Box_i A'_2\dto^{f'} \\
A''_1 \Box_j A''_2 \rto^(.6){g''} &B
\enddiagram
$$
which, according to Lemmas \ref{giantHex:Lemma} (if $i \neq j$) or \ref{coherence:hardLemma} (if $i = j$),
is commutative and therefore yields the equality $f = g$.

Thus the functor $\Lambda_k^n$ is also {\it injective} on morphisms
and the coherence theorem is completely proved.
\hfill
$\Box$

\newpage
\MySection{The Milgram Construction}

This section is devoted to a detailed discussion of the relation between the Milgram construction
and the premonad construction with respect to the Milgram subpreoperad
$\overline{\cal J}$ of the $n$-fold monoidal operad ${\cal M}_n$,
introduced in Section 3.

\begin{definition} Let $X$ be an object in the category
$\overline{\cal J}_n(k)$.  We denote by ${\cal S}(X)$
the full subcategory of $\overline{\cal J}_n(k)$ consisting of all the objects
$Y$ in ${\cal J}_n(k)$ which map
into $X$ (including $X$ itself).  As usual abusively we also use the same notation ${\cal S}(X)$ to denote the
nerve of this category.
\end{definition}

For $n=2$ the natural homeomorphism  of Theorem \ref{Milgram2:theorem} between the Milgram construction
and the premonad construction $\overline{J}_n(X)$ is a direct consequence of the following result:

\begin{theorem}\label{permutohedron:theorem}
${\cal S}(1\Box_2 2\Box_2 3\Box_2\dots \Box_2 k)$ is homeomorphic to the permutohedron
$P_k$. More precisely:
\begin{enumerate}
\item The simplicial triangulation of ${\cal S}(1\Box_2 2\Box_2 3\Box_2\dots \Box_2 k)$ arising from its
definition as a nerve is isomorphic to the barycentric subdivision of the natural cell structure on $P_k$.
\item There is a functorial action of the symmetric group $\Sigma_k$ on the category
${\cal S}(1\Box_2 2\Box_2 3\Box_2\dots \Box_2 k)$ inducing an action on its nerve which corresponds
under this isomorphism to the natural action of $\Sigma_k$ on $P_k$.
\item For each $i=1,2,\dots,k$ the functor 
${\cal S}(1\Box_2 2\Box_2\dots\Box_2 k)\longrightarrow{\cal S}\left(1\Box_2 2\Box_2\dots\Box_2 (k-1)\right)$
induced by the map of generating elements:
$$1\mapsto 1,\ 2\mapsto 2,\ \dots,\ (i-1)\mapsto (i-1),\ i\mapsto 0,\ (i+1)\mapsto i,\ (i+2)\mapsto (i+1),\ \dots,\ k\mapsto (k-1)$$
corresponds to the $i$-degeneracy map $D_i:P_k\to P_{k-1}$.
\end{enumerate}
\end{theorem}

Before we go on to the proof of this theorem, we illustrate this for the case $k=3$. Recall that $P_3$
is a hexagon.  Here is a picture of the nerve of ${\cal S}(1\Box_2 2\Box_2 3)$:
$$\diagram
&&1\Box_1 2\Box_1 3\dlto\drto\xto[3,0]\\
&(1\Box_2 2)\Box_1 3\ddrto&&1\Box_1(2\Box_2 3)\ddlto\\
2\Box_1 1\Box_1 3\drrto\urto\dto&&&&
1\Box_1 3\Box_1 2\dllto\ulto\dto\\
2\Box_1(1\Box_2 3)\rrto&&1\Box_2 2\Box_2 3&&(1\Box_2 3)\Box_1 2\llto\\
2\Box_1 3\Box_1 1\urrto\drto\uto&&&&
3\Box_1 1\Box_1 2\ullto\dlto\uto\\
&(2\Box_1 3)\Box_1 1\uurto&&3\Box_1(1\Box_2 2)\uulto\\
&&3\Box_1 2\Box_1 1\ulto\urto\xto[-3,0]
\enddiagram$$
(Here, as elsewhere throughout this section, we rely heavily on the coherence theorem for $n$-fold monoidal
categories. Thus we do not have to worry about labelling the arrows in our diagram, since there can be at most
one between any pair of objects, and the existence of the arrows shown can be easily checked.)

\bigskip
\noindent
{\bf Proof of Theorem \ref{permutohedron:theorem}.}  The coherence theorem implies that the objects of
${\cal S}(1\Box_2 2\Box_2 3\Box_2\dots \Box_2 k)$ have the form $A_1\Box_1 A_2\Box_1\dots\Box_1 A_s$,
with
$$A_r = i_{r1}\Box_2 i_{r2}\Box_2\dots\Box_2 i_{rj_r}\qquad 1\le i_{r1}<i_{r2}<\dots< i_{rj_r}\le k,$$
ie. $\left((i_{rt})_{t=1}^{j_r}\right)_{r=1}^s$ forms a $(j_1,j_2,\dots,j_s)$-shuffle in $\Sigma_k$.

We begin by defining a functorial action of the symmetric group $\Sigma_k$
on ${\cal S}(1\Box_2 2\Box_2 3\Box_2\dots \Box_2 k)$.  Given an element of $\sigma\in\Sigma_k$, there
is a functor
$$
{\cal S}(1\Box_2 2\Box_2 3\Box_2\dots \Box_2 k)\longrightarrow
{\cal S}\left(\sigma(1)\Box_2 \sigma(2)\Box_2 \sigma(3)\Box_2\dots \Box_2 \sigma(k)\right)
$$
given by permuting the generating elements $\{1,2,\dots,k\}$ according to $\sigma$. We compose this with the
functor
$$
{\cal S}\left(\sigma(1)\Box_2 \sigma(2)\Box_2 \sigma(3)\Box_2\dots \Box_2 \sigma(k)\right)
\longrightarrow{\cal S}(1\Box_2 2\Box_2 3\Box_2\dots \Box_2 k)
$$
which reorders the generating elements within the inner parentheses in their natural order when read from left to
right.  (To see that this defines a functor one must use the coherence theorem.)

To illustrate this action consider the totally order reversing permutation $[6,5,4,3,2,1]$ acting on the object
$(2\Box_2 4)\Box_1 (3\Box_2 5\Box_2 6)\Box_1 1\in
{\cal S}(1\Box_2 2\Box_2 3\Box_2 3\Box_2 4\Box_2 5 \Box_2 6)$. We have
$$
(2\Box_2 4)\Box_1 (3\Box_2 5\Box_2 6)\Box_1 1\mapsto
(5\Box_2 3)\Box_1 (4\Box_2 2\Box_2 1)\Box_1 6\mapsto
(3\Box_2 5)\Box_1 (1\Box_2 2\Box_2 4)\Box_1 6$$

We now proceed by induction on $k$ to prove part (1) of the theorem. For $k=1$ this is trivially true,
since both ${\cal S}(1)$ and $P_1$ are consist of a single point.  We then note that by the coherence theorem,
${\cal S}(1\Box_2 2\Box_2 3\Box_2\dots \Box_2 k)$ is the cone, with respect to the vertex
$1\Box_2 2\Box_2 3\Box_2\dots \Box_2 k$, of the union
$$\cup_{p=1}^{k-1}\cup_{\alpha\in Sh_{p,k-p}}\alpha
\left({\cal S}\left((1\Box_2 2\Box_2\dots\Box_2 p)\Box_1\left((p+1)\Box_2(p+2)\Box_2\dots\Box_2 k\right)\right)\right),
$$
where $Sh_{p,k-p}$ denotes the set of $(p,k-k)$ shuffles acting via the symmetric group action defined above.

Moreover by the coherence theorem, any object in
$${\cal S}\left((1\Box_2 2\Box_2\dots\Box_2 p)\Box_1\left((p+1)\Box_2(p+2)\Box_2\dots\Box_2 k\right)\right)$$
must have a canonical splitting $X_1\Box_1 X_2$, with $X_1$ in ${\cal S}\left(1\Box_2 2\Box_2\dots\Box_2 p\right)$.
Thus there is a canonical isomorphism
\begin{eqnarray*}
\lefteqn{{\cal S}\left((1\Box_2 2\Box_2\dots\Box_2 p)\Box_1\left((p+1)\Box_2(p+2)\Box_2\dots\Box_2 k\right)\right)}\\
&\cong&{\cal S}\left(1\Box_2 2\Box_2\dots\Box_2 p\right)\times{\cal S}\left(1\Box_2 2\Box_2\dots\Box_2 (k-p)\right)
\end{eqnarray*}
Hence ${\cal S}(1\Box_2 2\Box_2 3\Box_2\dots \Box_2 k)$ can be identified with the cone on
$$\cup_{p=1}^{k-1}\cup_{\alpha\in Sh_{p,k-p}}\alpha
\left({\cal S}\left(1\Box_2 2\Box_2\dots\Box_2 p\right)\times{\cal S}\left(1\Box_2 2\Box_2\dots\Box_2 (k-p)\right)
\right).
$$
Now according to \cite{Mi}, the boundary of the permutohedron $P_k$ has a similar decomposition as a union:
$$\cup_{p=1}^{k-1}\cup_{\alpha\in Sh_{p,k-p}}\alpha
\left(P_k\times P_{k-p}\right).
$$
Thus we construct our simplicial isomorphism by sending the vertex $1\Box_2 2\Box_2\dots\Box_2 k$ to the
barycenter of $P_k$ and then extending to the boundary by sending
$$\alpha
\left({\cal S}\left(1\Box_2 2\Box_2\dots\Box_2 p\right)\times{\cal S}\left(1\Box_2 2\Box_2\dots\Box_2 (k-p)\right)
\right)$$
to $\alpha\left(P_p\times P_{k-p}\right)$ via the inductively defined isomorphisms
${{\cal S}\left(1\Box_2 2\Box_2\dots\Box_2 p\right)\cong P_p}$ and \linebreak
${{\cal S}\left(1\Box_2 2\Box_2\dots\Box_2(k-p)\right)\cong P_{k-p}}$.

To check that this is well-defined, we note that if two codimension 1 faces 
$$\alpha\left({\cal S}\left(1\Box_2 2\Box_2\dots\Box_2 p\right)\times{\cal S}\left(1\Box_2 2\Box_2\dots\Box_2 (k-p)\right)
\right)$$ and 
$$\alpha'
\left({\cal S}\left(1\Box_2 2\Box_2\dots\Box_2 q\right)\times{\cal S}\left(1\Box_2 2\Box_2\dots\Box_2 (k-q)\right)
\right)$$
have a nonempty intersection, then we must have $p\ne q$ and the intersection must have the form
$$\beta
\left({\cal S}\left(1\Box_2 2\Box_2\dots\Box_2 u\right)\times
{\cal S}\left(1\Box_2 2\Box_2\dots\Box_2 v\right)
\times{\cal S}\left(1\Box_2 2\Box_2\dots\Box_2 w\right)\right),$$
where $u=\mbox{min}(p,q)$, $w=\mbox{min}(k-p,k-q)$, $v=k-u-w$, and $\beta$ is a $(u,v,w)$-shuffle.  Moreover
$\beta$ is determined as the only shuffle such that $\beta(1\Box_1 2\Box_1\dots\Box_1 k)$ is contained in both
codimension 1 faces.  We then note that the analogs of these facts are also true in $P_k$.

The rest of the proof is straightforward and is left as an exercise.

\bigskip

\noindent
{\bf Proof of Theorem \ref{Milgram2:theorem} for $\bold{n=2}$.} The Milgram construction for $n=2$ can be rearranged as the
premonad construction on the preoperad
whose $k$-th space is the quotient space $P_k\times\Sigma_k/\approx$. The equivalence relation $\approx$
identifies the codimension 1 face $\alpha\left(P_p\times P_{k-p}\right)$ in $P_k\times\{\sigma\}$ with the
codimension 1 face $P_p\times P_{k-p}$ in $P_k\times\{\alpha^{-1}\sigma\}$, for
any $(p,k-p)$-shuffle $\alpha$.

The preoperad space $\overline{\cal J}_2(k)$ can be similarly expressed as a similar quotient space
${\cal S}(1\Box_2 2\Box_2 3\Box_2\dots \Box_2 k)\times\Sigma_k/\approx$, where we identify
${\cal S}(1\Box_2 2\Box_2 3\Box_2\dots \Box_2 k)\times\{\sigma\}$ with
${\cal S}\left(\sigma(1)\Box_2 \sigma(2)\Box_2 \sigma(3)\Box_2\dots \Box_2 \sigma(k)\right)$.
The result now follows directly from Theorem \ref{permutohedron:theorem}.

\bigskip
The following is left as an exercise for the interested reader. It gives an
intrinsic description in terms of generators and relations of the categories
$\overline{\cal J}_2(k)$ and ${\cal J}_2(k)/\Sigma_k$. (It is not difficult
to do the exercise with the help of the coherence theorem.)

\begin{exercise} A $J_2$ functor is a functor $F: {\cal A}\longrightarrow {\cal B}$
between monoidal categories which is strongly monoidal in the following sense.
Denote by $\Box_2$, $\Box_1$ the monoid operations in ${\cal A}$, ${\cal B}$ respectively
and by 0 the unit object in either category. Then we require that $F(0)=0$ and that
for each $(p,q)$-shuffle $\sigma\in\Sigma_k$ there is given a natural transformation
$$\xymatrix@R=-5pt @C=-130pt{
\zeta_{A_1,A_2,\dots,A_p;A_{p+1},A_{p+2},\dots,A_k}^\sigma:
F(A_1\Box_2\dots\Box_2 A_p)\Box_1 F(A_{p+1}\Box_2\dots\Box_2 A_k)\\
&\longrightarrow
\quad F(A_{\sigma^{-1}1}\Box_2 A_{\sigma^{-1}2}\Box_2\dots\Box_2 A_{\sigma^{-1}k})
}$$
satisfying the following properties:
\begin{enumerate}
\item (Unit condition) $\zeta_{A_1,A_2,\dots,A_p;0,0,\dots 0}^\sigma=
id_{F(A_1\Box_2 A_2\Box_2\dots\Box_2 A_p)}$ and\newline
$\zeta_{0,0\dots,0;A_{p+1},A_{p+2},\dots,A_k}^\sigma=
id_{F(A_{p+1}\Box_2 A_{p+2}\Box_2\dots\Box_2 A_k)}$.
\item (Substitution Property) If $A_i=B_{ij_1}\Box_2 B_{ij_2}\Box_2\dots\Box_2 B_{ij_i}$,
then
$$\zeta_{A_1,A_2,\dots,A_p;A_{p+1},A_{p+2},\dots,A_k}^\sigma =
\zeta_{B_{11},\dots,B_{1j_1},B_{21},\dots,B_{pj_p};B_{j+1,1},\dots,B_{kj_k}}^{\sigma'}$$
where $\sigma'\in\Sigma_{j_1+j_2+\dots+j_k}$ permutes the blocks $\{1,2,\dots,j_1\}$,
$\{j_1+1,j_1+2,\dots,j_1+j_2\}$, \dots, $\{j_1+\dots+j_{k-1}+1,j_1+\dots+j_{k-1}+2,\dots,
j_1+\dots+j_{k-1}+j_k\}$ the same way that $\sigma$ permutes 1, 2, \dots $k$.
A special case of this is if any $A_i=0$ (ie. a 0-fold $\Box_2$-sum), in which case
the resulting $\zeta$ is the same as one where the corresponding 0 entries have been
deleted.
\item (Associativity) Given $p+q+r=k$, a $(p,q)$-shuffle $\sigma$, a $(p+q,r)$-shuffle
$\tau$, a $(q,r)$-shuffle $\kappa$ and a $(p,q+r)$-shuffle $\lambda$, such that
$\lambda(id\oplus\kappa)=\tau(\sigma\oplus id)=\gamma$ in $\Sigma_k$, then the following
diagram commutes
$$\def\objectstyle{\scriptscriptstyle}\def\labelstyle{\scriptscriptstyle}
\xymatrix@C=-90pt{
&{F(A_1\Box_2\dots\Box_2 A_p)\Box_1 F(A_{p+1}\Box_2\dots\Box_2 A_{p+q})\Box_1 F(A_{p+q+1}\Box_2\dots
\Box_2 A_k)}\ar[ddr]^{\quad\zeta_{A_1,\dots,A_p;A_{p+1},\dots A_{p+q}}^\sigma\Box_1 id}
\ar[ddl]^(.6){\qquad id\Box_1\zeta_{A_{p+1},\dots A_{p+q};A_{p+q+1},\dots,A_k}^{\kappa}}\\
\\
\vbox{\vspace{1ex}}\vtop{\vspace{1ex}}
{F(A_1\Box_2\dots\Box_2 A_p)\Box_1 F(A_{\tilde{\kappa}^{-1}(p+1)}\Box_2\dots\Box_2 A_{\tilde{\kappa}^{-1}n})}
\ar[ddr]^(.4){\ \zeta_{A_1,\dots,A_p;A_{\tilde{\kappa}^{-1}(p+1)},\dots,A_{\tilde{\kappa}^{-1}n}}^\lambda}
&&\vbox{\vspace{1ex}}\vtop{\vspace{1ex}}{F(A_{\sigma^{-1}1}\Box_2\dots\Box_2 A_{\sigma^{-1}(p+q)})\Box_1 F(A_{p+q+1}\Box_2\dots
\Box_2 A_k)}\ar[ddl]^{\qquad\zeta_{A_{\sigma^{-1}1},\dots,A_{\sigma^{-1}(p+q)};A_{p+q+1},\dots,
A_k}^\tau}\\
\\
&\vbox{\vspace{2ex}}{F(A_{\gamma^{-1}1}\Box_2 A_{\gamma^{-1}2}\Box_2\dots\Box_2 A_{\gamma^{-1}n})}
}$$
where $\tilde{\kappa}$ is the translation of $\kappa$ to the set $\{p+1,p+2,\dots,k\}$.
\end{enumerate}

Show that $\overline{\cal J}_2(k)/\Sigma_k$ is the target of the universal
$J_2$ functor from the free monoidal category on one object. Similarly
$\overline{\cal J}_2(k)$ can be described as
a subcategory of the universal $J_2$ functor from the free monoidal category on
$\{1, 2,\dots, k\}$.
\end{exercise}

The basic building block of the Milgram construction $J_n(X)$ for $n>2$
is the product $(P_k)^{n-1}$.
In order to relate $\overline{J}_n(X)$ to the Milgram construction, we have to
relate $(P_k)^{n-1}$ to ${\cal S}(1\Box_n 2\Box_n\dots\Box_n k)$. Unfortunately
the analog of Theorem \ref{permutohedron:theorem} breaks down when $n>2$:
${\cal S}(1\Box_n 2\Box_n\dots\Box_n k)$ is not isomorphic as a cell complex
to $(P_k)^{n-1}$ but rather to a quotient of $(P_k)^{n-1}$.  Nevertheless
as we show below, ${\cal S}(1\Box_n 2\Box_n\dots\Box_n k)$ is homeomorphic
to a disk of dimension $(n-1)(k-1)$, and thus also to $(P_k)^{n-1}$.

\bigskip\noindent
\begin{definition}
Let $A$ and $B$ be two objects of $P_k={\cal S}(1\Box_2 2\Box_2 3\Box_2\dots\Box_2 k)$ (using
Theorem \ref{permutohedron:theorem}). Suppose that $A=A_1\Box_1 A_2\Box_1\dots\Box_1 A_p$ where
each $A_i$ is $\Box_1$-irreducible. Define a new object $\pi_A(B)$ of $P_k$ by
$$\pi_A(B)= (B\cap|A_1|)\Box_1(B\cap|A_2|)\Box_1(B\cap|A_3|)\Box_1\dots\Box_1(B\cap|A_p|)$$
It is obvious that this induces a map of posets
$$\pi_A:P_k\longrightarrow P_k$$
retracting $P_k$ onto the face ${\cal S}(A)$.
\end{definition}

We collect here, for future reference, the following basic properties of the
retractions $\pi_A$:

\begin{lemma}\label{retractions:properties}
If $A$ and $B$ are objects of $P_k={\cal S}(1\Box_2 2\Box_2 3\Box_2\dots\Box_2 k)$, then
\begin{description}
\item[{\rm (i)}] $\pi_A\pi_B=\pi_{\pi_A(B)}$
\item[{\rm (ii)}] If ${\cal S}(A)\cap{\cal S}(B)={\cal S}(C)$, then
$$\pi_A(B)=\pi_B(A)=C,$$
and consequently
$$\pi_A\pi_B=\pi_B\pi_A=\pi_C.$$
\item[{\rm (iii)}] If $\sigma\in\Sigma_k$, then
$$\sigma\pi_A(B)=\pi_{\sigma A}(\sigma B).$$
\end{description}
\end{lemma}

\begin{definition}
Given a based space $X$ we define the {\it thick Milgram construction\/ } $\widetilde{J}_n(X)$
to be the quotient of the disjoint union $\coprod_{k\ge0}(P_k)^{n-1}\times X^k$ by the equivalence  
relation generated by the relations
\begin{description}
\item[{\rm (i)}]$(c_1,c_2,\dots,c_{n-1};x_1,\dots,x_{i-1},*,x_i,\dots,x_{k-1})\approx
(s_i(c_1),s_i(c_2),\dots,s_i(c_{n-1});x_1,\dots,x_{i-1},x_i,\dots,x_{k-1})$
\item[{\rm (ii)}]If some $c_i$ is in a boundary face $\alpha(P_p\times P_q)$,
where $\alpha$ is a $(p,q)$ shuffle in $\Sigma_k$, then
$$(c_1,c_2,\dots,c_{n-1};x_1,x_2,\dots,x_{k})\approx
\left(\alpha^{-1}(c_1),\alpha^{-1}(c_2),\dots,\alpha^{-1}(c_{n-1});
\alpha(x_1,x_2,\dots,x_k)\right)$$
\end{description}
We define the {\it Milgram construction\/} $J_n(X)$ to be the quotient of the
thick Milgram construction by the following additional equivalence relations:
\begin{description}
\item[{\rm (iii)}] If $c_i$ is in a boundary face ${\cal S}(A)$, then
$$(c_1,c_2,\dots,c_i,c_{i+1},\dots,c_{n-1};x_1,x_2,\dots,x_{k})\approx
(c_1,c_2,\dots,c_i,\pi_A(c_{i+1}),\dots,\pi_A(c_{n-1});x_1,x_2,\dots,x_{k})$$
\end{description}
Finally we define the {\it thin Milgram construction\/} $\widehat{J}_n(X)$
by conically extending the relations (iii) to the interior of $P_k^{n-1}$:
\begin{description}
\item[{\rm (iii')}] If $(c_1,c_2,\dots,c_{n-1})$ is in the cone (with respect
to $(1\Box_2 2\Box_2\dots\Box_2 k,1\Box_2 2\Box_2\dots\Box_2 k,\dots,
1\Box_2 2\Box_2\dots\Box_2 k)$) of $P_k\times\dots\times P_k\times{\cal S}(A)
\times P_k\times\dots\times P_k$, then
$$(c_1,c_2,\dots,c_{n-1};x_1,x_2,\dots,x_{k})\approx
(c'_1,c'_2,\dots,c'_{n-1};x_1,x_2,\dots,x_{k}),$$
where $(c'_1,c'_2,\dots,c'_{n-1})$ is the image of $(c_1,c_2,\dots,c_{n-1})$
under the conical extension of the map
$$(d_1.d_2,\dots,d_i,d_{i+1},\dots,d_{n-1})\mapsto
\left(d_1.d_2,\dots,d_i,\pi_A(d_{i+1}),\dots,\pi_A(d_{n-1})\right)$$
on the boundary face.
\end{description}
\end{definition}

\begin{remark}In Milgram's own description of his construction $J_n(X)$, the
relations (ii) and (iii) are combined into a single relation, cf. 
\cite[p. 24]{Mi2}.
\end{remark}

It is clear that each of the above variants of Milgram's construction arises
as the premonad construction on a preoperad. The preoperad $\widetilde{\cal J}_n$
associated with the thick Milgram construction has the form
$$\widetilde{\cal J}_n(k)=P_k^{n-1}\times\Sigma_k/\approx,$$
where the equivalence relation identifies a point
$(c_1,c_2,\dots,c_{n-1};\sigma)$, if some $c_i$ is in $\alpha(P_p\times P_q)$,
with the point $(\alpha^{-1}(c_1),\alpha^{-1}(c_2),\dots,\alpha^{-1}(c_{n-1});
\alpha\sigma)$. The unit maps $s_i:\widetilde{\cal J}_n(k)\to
\widetilde{\cal J}_n(k-1)$ are applied coordinatewise:
$$s_i(c_1,c_2,\dots,c_{n-1};\sigma)=\left(s_i(c_1),s_i(c_2)\dots,s_i(c_{n-1});
s_i(\sigma)\right)$$
The preoperad ${\cal J}_n$, corresponding to the Milgram construction,
is obtained from $\widetilde{\cal J}_n$ by
taking the quotient of $P_k^{n-1}$ by the relations
\begin{description}
\item[{\rm (*)}] If $c_i$ is in a boundary face ${\cal S}(A)$, then
$$(c_1,c_2,\dots,c_i,c_{i+1},\dots,c_{n-1})\approx
(c_1,c_2,\dots,c_i,\pi_A(c_{i+1}),\dots,\pi_A(c_{n-1})).$$
\end{description}
The preoperad $\widehat{\cal J}_n$, corresponding to the thin Milgram
construction, is obtained from $\widetilde{\cal J}_n$ by
taking the quotient of $P_k^{n-1}$ by the conical extension of the
relations (*):
\begin{description}
\item[{\rm (**)}] If $(c_1,c_2,\dots,c_{n-1})$ is in the cone (with respect
to $(1\Box_2 2\Box_2\dots\Box_2 k,1\Box_2 2\Box_2\dots\Box_2 k,\dots,
1\Box_2 2\Box_2\dots\Box_2 k)$) of $P_k\times\dots\times P_k\times{\cal S}(A)
\times P_k\times\dots\times P_k$, then
$$(c_1,c_2,\dots,c_{n-1})\approx(c'_1,c'_2,\dots,c'_{n-1}),$$
where $(c'_1,c'_2,\dots,c'_{n-1})$ is the image of $(c_1,c_2,\dots,c_{n-1})$
under the conical extension of the map
$$(d_1.d_2,\dots,d_i,d_{i+1},\dots,d_{n-1})\mapsto
\left(d_1.d_2,\dots,d_i,\pi_A(d_{i+1}),\dots,\pi_A(d_{n-1})\right)$$
on the boundary face.
\end{description}
Thus all of these preoperads take the generic form
\begin{eqnarray*}
\widetilde{\cal J}_n(k) &= &\widetilde{D}_n(k)\times\Sigma_k/\approx\\
{\cal J}_n(k) &= &D_n(k)\times\Sigma_k/\approx\\
\widehat{\cal J}_n(k) &= &\widehat{D}_n(k)\times\Sigma_k/\approx
\end{eqnarray*}
where $\widetilde{D}_n(k)=P_k^{n-1}$ and $D_n(k)$, resp. $\widehat{D}_n(k)$,
are the quotients of $P_k^{n-1}$ by the relations (*), resp. (**).
(One needs Lemma \ref{retractions:properties} to check that the relations
(*) and (**) commute with the equivariancy relations used to glue together
the $k!$ copies of these quotients.) The
following pictures illustrate these constructions for $n=3$ and $k=2$:

\centerline{\psfig{figure=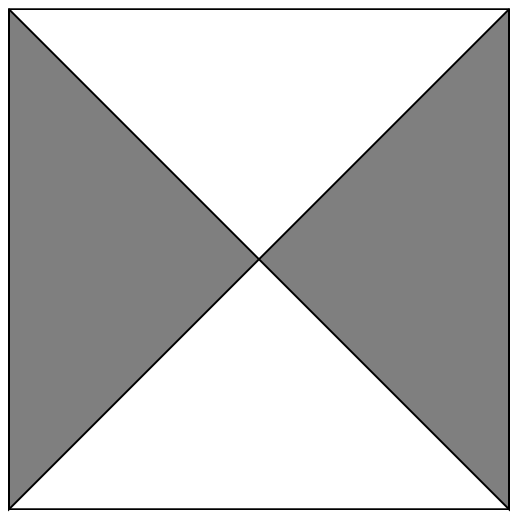}\qquad\psfig{figure=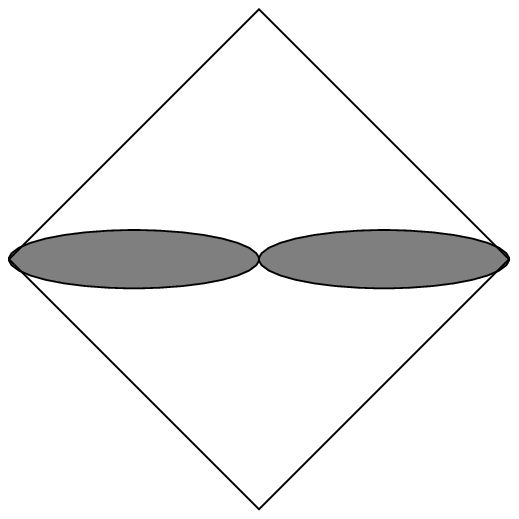}
\qquad\psfig{figure=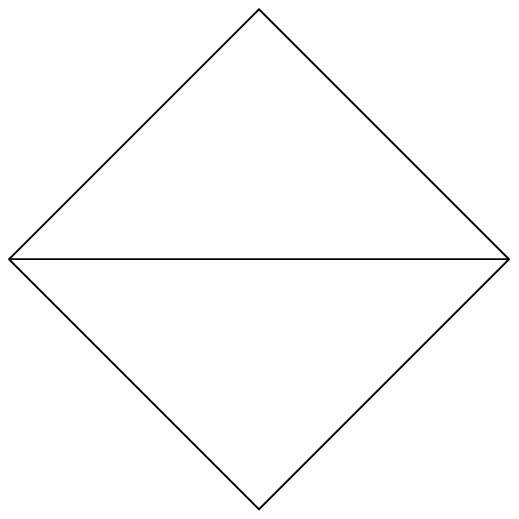}}

\vspace{10pt}
\noindent
The first picture shows $\widetilde{D}_3(2)=P_2^{2}=I\times I$. The second
picture shows $D_3(2)$, which is obtained from the first picture by collapsing
the vertical sides of the shaded triangles to points. The collapsed triangles
become ``polygons'' with two sides. The third picture shows $\widehat{D}_3(2)$,
which is obtained from the first picture by collapsing the shaded triangles
to horizontal lines.

\bigskip
To complete the proof of Theorem \ref{MilgramN:theorem} we will need the
following elementary result from PL topology:

\begin{lemma}\label{elementary:lemma}Let $D^i$ denote the $i$-dimensional disk.
\begin{description}
\item[{\rm (a)}] If $D^m\subset\partial D^n$ is a PL
imbedded disk and $\phi:D^m\to D^k$ is an elementary collapse to
a boundary face, then $D^n\cup_{D^m}D^k\cong D^n$
\item[{\rm (b)}] If $D^n$, $\phi:D^m\to D^k$ are as in (a)
and $C_pX$ denotes the cone with respect to an interior point $p\in D^n$,
then $D^n\cup_{C_pD^m}C_pD^k\cong D^n$.
\end{description}
\end{lemma}

\noindent
{\bf Proof:} We first prove part (a) for the case $m=n-1$. We take as
a model for $D^n$ the prism $\Delta^{n-1}\times I$ and we take the boundary
disk we are collapsing to be the top face $\Delta^{n-1}\times\{1\}$.
(That we can arrange this follows from the Disk Theorem of PL topology,
cf. \cite[p. 56]{RS}.)
Let $K$ denote the convex hull in $\Delta^{n-1}\times I$ of
$\Delta^{n-1}\times\{0\}$ and $\Delta^k\times\{1\}$.  Then $K$ is
obviously an $n$-dimensional topological disk. Now consider the map of pairs
$\lambda:(\Delta^{n-1}\times I,\Delta^{n-1}\times \{1\}\longrightarrow
(K,\Delta^k\times\{1\})$ given by the formula
$$(\mathbf{x},t)\mapsto\left((1-t)\mathbf{x}+t\phi(\mathbf{x}),t\right).$$
This map is a relative homeomorphism, since if $(\mathbf{x_1},t)$ and
$(\mathbf{x_2},t)$ both mapped to the same point for some $t<1$, then
the vectors $\mathbf{x_1}-\mathbf{x_2}$ and
$\phi(\mathbf{x_1})-\phi(\mathbf{x_2})$ would have to point in opposite
directions, which can't happen for a linear retraction $\phi$. Since the
restriction of $\lambda$ to $\Delta^{n-1}\times \{1\}\to\Delta^k\times \{1\}$
is just $\phi$, $\lambda$ induces a homeomorphism
$$D^n\cup_{D^m}D^k=(\Delta^{n-1}\times I)\cup_{\Delta^{n-1}\times\{1\}}
\Delta^k\times \{1\}\cong K\cong D^n$$
We can reduce the general case of part (a) to the special case proved
above as follows:
$$D^n\cup_{D^m}D^k\cong\left(D^n\cup_{D^{n-1}}D^m\right)\cup_{D^m}D^k\cong
D^n\cup_{D^{n-1}}D^k\cong D^n.$$

Finally we can reduce part (b) to part (a) as follows. Cut apart the given
disk $D^n$ along a suitable codimension 1 subdisk passing through the point
$p$. (If $m=n-1$, excise the interior of the cone $C_pD^m$ first.) Then
the we can realize $D^n\cup_{C_pD^m}C_pD^k$ as the result of
a two step process. In the first step we are attaching $D^{k+1}$ to each
of the two $n$-dimensional disks we obtained after the cut, by an attachment
of the form given in part (a). By part (a) we know that the resulting spaces
are homeomorphic to $n$-disks.  In the second step we glue together these
disks along the parts of the boundaries of the two pieces which were originally
$(n-1)$-disks (where we originally made the cut), but where we attached
$D^{k+1}$'s. That the resulting part of the boundaries are still $(n-1)$-disks
follows by noting that the complementary parts of the boundaries are PL
imbedded $(n-1)$-disks.

\bigskip
\noindent
{\bf Proof of Theorem \ref{MilgramN:theorem}:} $\widetilde{\cal J}_n(k)=
P_k^{n-1}$ is evidently a $(k-1)(n-1)$-dimensional disk. By repeated
applications of part (a) of Lemma \ref{elementary:lemma}, so is
${\cal J}_n(k)$.  By repeated applications of part (b) of Lemma
\ref{elementary:lemma}, $\widehat{\cal J}_n(k)$ is also a
$(k-1)(n-1)$-dimensional disk.

We now construct a map of posets $q:\widetilde{\cal J}_n(k)=P_k^{n-1}\to
\overline{\cal J}_n(k)={\cal S}(1\Box_n 2\Box_n\dots\Box_n k)$ as follows.
Given $(A_1,A_2,\dots,A_{n-1})\in P_k^{n-1}=
\left({\cal S}(1\Box_n 2\Box_n\dots\Box_n k)\right)^{n-1}$, first replace
it by
$$(B_1,B_2,\dots,B_{n-1})=\left(A_1,\pi_{A_1}(A_2),\pi_{A_1}\pi_{A_2}(A_3),
\dots,\pi_{A_1}\pi_{A_2}\dots\pi_{A_{n-2}}(A_{n-1})\right)$$
We then have
$$B_{n-1}\le B_{n-2}\le\dots\le B_2\le B_1,$$
from which it follows that the parenthesization of any object $B_i$
induces a (usually redundant) parenthesization of the object
$B_{i+1}$. From this it follows that we can endow $B_{n-1}$ with $n-1$
levels of parentheses: the innermost coming from the original parenthesization
of $B_{n-1}$, the next level coming from the parenthesization of $B_{n-1}$,
\dots, the outermost level coming from the parenthesization of $B_1$.
Now define $q(A_1,A_2,\dots,A_{n-1})$ to be the object constructed from
this heavily parenthesized version of $B_{n-1}$ as follows. Replace each
$\Box_2$ in the innermost level of parentheses by $\Box_n$. Then replace
each $\Box_1$ in the next level of parentheses by $\Box_{n-1}$, etc. At the
penultimate step replace each $\Box_1$ in the next to outermost level of
parentheses by $\Box_2$. At the final step leave the outermost $\Box_1$'s
alone.

The following example (with $n=4$, $k=5$) illustrates this process. Let
$$(A_1,A_2,A_3)=\left((1\Box_2 3)\Box_1(2\Box_2 4\Box_2 5),
(1\Box_2 3\Box_2 4)\Box_1(2\Box_2 5),(1\Box_2 2\Box_2 4\Box_2 5)\Box_1 3
\right).$$
Then
$$(B_1,B_2,B_3)=\left((1\Box_2 3)\Box_1(2\Box_2 4\Box_2 5),
(1\Box_2 3)\Box_1 4\Box_1(2\Box_2 5),3\Box_1 1\Box_1 4\Box_1(2\Box_2 5)
\right),$$
and the resulting redundant parenthesization of $B_3$ is
$$B_3=(((3)\Box_1(1)))\Box_1((4)\Box_1((2\Box_2 5))$$
Thus
$$q(A_1,A_2,A_3)=(3\Box_3 1)\Box_1(4\Box_2(2\Box_4 5)).$$

It is easy to see that this map of posets extends to a map of preoperads
$q:\widetilde{\cal J}_n\to\overline{\cal J}_n$, which factors through a
map of preoperads $q':\widehat{\cal J}_n\to\overline{\cal J}_n$. To check
that $q'$ is a simplicial isomorphism, it is only necessary to note that
$X\le Y$ in ${\cal S}(1\Box_n 2\Box_n\dots\Box_n k)$ if and only if we
can find $A$, $B$ in $(P_k)^{n-1}$ so that $q(A)=X$, $q(B)=Y$ and $A\le B$.

Thus we have quotient maps of preoperads
$$\widetilde{\cal J}_n\stackrel{q_1}{\longrightarrow}
{\cal J}_n\stackrel{q_2}{\longrightarrow}
\widehat{\cal J}_n\cong\overline{\cal J}_n$$
We have that $q_2$ is an equivalence, since $D_n(k)\to\widehat{D}_n(k)$
are both given by elementary collapses of $\widetilde{D}_n(k)$, and since
the gluings of the $k!$ copies of $D_n(k)$, resp. $\widehat{D}_n(k)$ are
along the boundaries where $D_n(k)$ and $\widehat{D}_n(k)$ are isomorphic.

So it remains to show that
$$q=q_2q_1:\widetilde{\cal J}_n\longrightarrow \overline{\cal J}_n$$
is an equivalence. Since this map is given by a map of posets, we use
Quillen's Theorem A: we show that for any object in the poset
$\overline{\cal J}_n$ the overcategory of objects in $\widetilde{\cal J}_n$
is contractible. But this is easy: elementary collapses of this overcategory
given by relations (**) above gives the cone over that object in
$\overline{\cal J}_n$. Since the cone is obviously contractible, and
elementary collapse do not change the homotopy type, the overcategory
must be contractible too.  This completes the proof.

\bigskip
\begin{remark} In \cite[p. 55]{GJ}, Getzler and Jones consider a poset closely
related to $\overline{\cal J}_n(k)$. More precisely their poset is isomorphic
to the ``dual'' $\overline{\cal J}_n(k)^*$. By this we mean the full
subcategory of ${\cal M}_n(k)$ consisting of objects whose nesting of
operations is opposite to those in $\overline{\cal J}_n(k)$: the $\Box_n$
operations are nested on the outermost level, the $\Box_{n-1}$ operations are
nested at the next outermost level,\dots, the $\Box_1$ operations are nested
at the innermost level. Getzler and Jones denote the objects of their poset
as ``multiple bar codes'': permutations $\sigma\in\Sigma_k$ with their elements
separated by multiple bars
$$\sigma(1)|_{i_1}\sigma(2)|_{i_2}\dots|_{i_{k-1}}\sigma(k)$$
where the subscript on each bar is $\le n$ and denotes the number of times
the bar is supposed to be repeated. The poset isomorphism with
$\overline{\cal J}_n(k)^*$ is given by the replacement $|_i\mapsto\Box_i$,
with the resulting object parenthesized according to the operation precedence
rules: $\Box_1$ has the highest precedence, $\Box_2$ has the next highest
precedence,\dots, $\Box_n$ has the lowest precedence.

There is a duality anti-automorphism of ${\cal M}_n(k)$ given by
$\Box_i\mapsto\Box_{n-i+1}$, which is easy to verify using the coherence
theorem. This anti-automorphism takes $\overline{\cal J}_n(k)$ to
$\overline{\cal J}_n(k)^*$.  Thus $\overline{\cal J}_n(k)$ is anti-isomorphic
to $\overline{\cal J}_n(k)^*$ and hence also to the Getzler-Jones poset. It
follows that the nerve of $\overline{\cal J}_n(k)$ is isomorphic to the
nerve of the Getzler-Jones poset.

Getzler and Jones also consider an operad freely generated by these posets.
This operad obviously maps into our operad ${\cal M}_n(k)$. It will be
shown in a forthcoming paper that this map of operads is an equivalence.

There is also an extensive discussion of the Getzler-Jones posets and their
relation to various other constructions in \cite[p. 46]{Be}.
\end{remark}

\newpage
\MySection{Relation to Little $n$-Cubes}\label{little:cubessection}
Boardman and Vogt \cite{BV1} introduced the little $n$-cubes operad to parametrize multiplications
on an $n$-fold loop space.  Later May \cite{May} used these operads to construct small models of
$\Omega^n S^n X$, an alternative to Milgram's models.  This section is devoted to the proof of our
main result Theorem \ref{main:theorem}, relating the $n$-fold monoidal operad ${\cal M}_n$ to the little $n$-cubes
operad ${\cal C}_n$, and then derive some consequences relating $n$-fold monoidal categories to
$n$-fold loop spaces.

We begin by associating to each object of ${\cal M}_n(k)$ a contractible space of $k$-fold configurations
of little $n$-cubes.

\begin{definition}\label{G:definition}
We think of a little $n$-cube $c$ as a product of closed subintervals of the unit
interval.  Thus the elements of the $k$-th space of little $n$-cubes ${\cal C}_n(k)$ have the form
$(c_1,c_2,\dots,c_k)$ where
$$c_j=[u_{j1},v_{j1}]\times[u_{j2},v_{j2}]\times\dots\times[u_{jn},v_{jn}]$$
where the interiors of the little $n$-cubes $c_j$ are required to be pairwise disjoint. If
$$c=[u_1,v_1]\times[u_2,v_2]\times\dots\times[u_n,v_n]\qquad
d=[z_1,w_1]\times[z_2,w_2]\times\dots\times[z_n,w_n]$$
are little $n$-cubes we write $c<_i d$ to mean that $v_i\le z_i$.  Equivalently $c<_i d$ if there
is a hyperplane perpendicular to the $i$-coordinate axis such that the interior of $c$ lies on the
negative side of the hyperplane and the interior of $d$ on the positive side of the hyperplane. (Note
that this condition forces the interiors of $c$ and $d$ to be disjoint.)  If $A$ is an object of
${\cal M}_n(k)$ let $G(A)$ denote the space of all $k$-fold configurations of little $n$-cubes satisfying the
following conditions:
$$G(A)=\{(c_1,c_2,\dots,c_k)\ |\ c_a<_i c_b \mbox{ if } a\Box_i b \mbox{ in } A\}$$
(cf. Definition \ref{binaryRelation:definition}).  Note that we do not have to explicitly require that $G(A)$ be a subspace of
${\cal C}_n(k)$ -- the ordering relations defining $G(A)$ force the little $n$-cubes in a configuration
in $G(A)$ to have pairwise disjoint interiors thus forcing the configuration to be in ${\cal C}_n(k)$.
Because of this, $G(A)$ may be identified with a convex, hence contractible, subspace of
${\Bbb R}^{2k}$ given by a set of inequalities between the coordinates.  For the same reason $G(A)$
is a closed subspace of ${\cal C}_n(k)$ (but not of ${\Bbb R}^{2k}$ and hence not compact, since the
requirement that each little $n$-cube have a nonvacuous interior is an open condition given by strict
inequalities).
\end{definition}

\noindent
\begin{example}
The following two figures represent configurations belonging to
$G((1\Box_2 2)\Box_1(3\Box_2 4))$:

\centerline{\hbox{
\unitlength1cm
\begin{picture}(12, 4.5)
\put(0, 0) {\line(1, 0) {4}}
\put(0, 0) {\line(0, 1) {4}}
\put(0, 4) {\line(1, 0) {4}}
\put(4, 0) {\line(0, 1) {4}}
\put(2, 0) {\line(0, 1) {4}}
\put(0, 2) {\line(1, 0) {2}}
\put(2, 1) {\line(1, 0) {2}}
\put(0.9, 0.7) {1}
\put(2.9, 0.4) {3}
\put(0.9, 2.7) {2}
\put(2.9, 2.2) {4}
\put(8, 0) {\line(1, 0) {4}}
\put(8, 0) {\line(0, 1) {4}}
\put(8, 4) {\line(1, 0) {4}}
\put(12, 0) {\line(0, 1) {4}}
\put(10, 0) {\line(0, 1) {4}}
\put(8, 1) {\line(1, 0) {2}}
\put(10, 2) {\line(1, 0) {2}}
\put(8.9, 0.4) {1}
\put(10.9, 0.7) {3}
\put(8.9, 2.2) {2}
\put(10.9, 2.7) {4}
\end{picture}
}}
More generally such configurations could have the subcube $i$ properly contained in the region marked $i$ $1\le i\le 4$.
\end{example}

\begin{remark}\label{decomposableCubes:remark}
If $A=B\Box_i C$ then for any configuration $(d_1,d_2,\dots,d_k)$ in $G(B\Box_i C)$,
we can find a hyperplane perpendicular to the $i$-th coordinate axis such that all little cubes in the configuration
having labels coming from the generating objects in $B$ have their interiors on the negative side of the
hyperplane and all little cubes in the configuration whose labels come from the generating objects in $C$
have their interiors on the negative side of the hyperplane.  For if
$$d_j=[u_{j1},v_{j1}]\times[u_{j2},v_{j2}]\times\dots\times[u_{jn},v_{jn}]$$
and if we let
$$M=\mbox{max}\{v_{bi}\ |\ b\mbox{ in } B\}\qquad m=\mbox{min}\{u_{ci}\ |\ c\mbox{ in } C\},$$
then the conditions that $(d_1,d_2,\dots,d_k)$ must satisfy in order to be in $G(B\Box_i C)$ imply that
$M\le m$.  Thus we can take $x_i = M$ as a separating hyperplane with the required properties. ($x_i=m$ 
would also work, as would any hyperplane in between those two.) It follows from this observation that
$$\cup_{A\in obj({\cal M}_n(k))}G(A)={\cal D}_n(k)$$
where ${\cal D}_n(k)\subset{\cal C}_n(k)$ is the subspace of {\it decomposable configurations\/} of
little cubes. Decomposability is defined recursively as follows.  First of all a configuration consisting
of a single $n$-cube, ie. an element of ${\cal C}_n(1)$, is declared to be decomposable. For a $k$-fold
configuration to be decomposable, we require that there be a hyperplane perpendicular to one of the
coordinate axes which does not pass through the interior of any little $n$-cube in the configuration and
which divides the configuration into two proper subconfigurations.  We further require that the
subconfigurations on both sides of the separating hyperplane to be themselves decomposable.  It is
trivially true that all ${\cal C}_1(k)$, ${\cal C}_n(1)$, and ${\cal C}_n(2)$ consist entirely of
decomposable configurations.  The same is also true for ${\cal C}_2(3)$, but all other spaces in the
little $n$-cubes operads contain nondecomposable configurations.  For instance $(c_1, c_2, c_3)$ where
$$c_1=[0,\frac{1}{2}]\times[0,1]\times[0,\frac{1}{2}]\quad
c_2 = [0,1]\times[0,\frac{1}{2}]\times[\frac{1}{2},1]\quad
c_3 = [\frac{1}{2},1]\times[\frac{1}{2},1]\times[0,1]$$
is a nondecomposable configuration of little 3-cubes in ${\cal C}_3(3)$ 
and the following figure shows a 4-fold configuration of little 2-cubes which is nondecomposable.
$$\mbox{
\centerline{\hbox{\unitlength1cm
\begin{picture}(4.5, 4.5)
\put(0, 0){\line(0, 1) {4} }
\put(0, 0){\line(1, 0) {4} }
\put(0, 4){\line(1, 0) {4} }
\put(4, 0){\line(0, 1) {4} }
\put(0, 1){\line(1, 0) {3} }
\put(1, 3){\line(1, 0) {3} }
\put(1, 1){\line(0, 1) {3} }
\put(3, 0){\line(0, 1) {3} }
\put(1.4, 0.4){1}
\put(3.4, 1.4){2}
\put(2.4, 3.4){4}
\put(0.4, 2.4){3}
\end{picture}}
}
}$$
The decomposable little $n$-cubes form a suboperad ${\cal D}_n$ of ${\cal C}_n$.  By sufficiently
shrinking every little $n$-cube in a configuration towards its barycenter, we can convert
any configuration into a decomposable one. This shows that the inclusion ${\cal D}_n\subset{\cal C}_n$
is an equivalence of operads.  The operad ${\cal D}_n$ was studied by Dunn \cite{D} who showed it is
homeomorphic to the $n$-fold tensor product ${\cal C}_1\otimes{\cal C}_1\otimes\dots\otimes{\cal C}_1$ of
the little 1-cubes operad.
\end{remark}

The assignment $A\mapsto G(A)$ is only defined on objects, not on morphisms.  In order to construct
a functor on ${\cal M}_n(k)$ we proceed as follows:

\begin{definition}
For any object $A\in{\cal M}_n(k)$ define
$$F(A)=\cup_{X\to A} G(X),$$
where the union is indexed over all objects $X\in{\cal M}_n(k)$ which map into $A$.  Then by definition
given a morphism $B\longrightarrow A$ in ${\cal M}_n(k)$, there is an induced inclusion of subspaces
$F(B)\subset F(A)$.  Thus we have constructed a functor $F:{\cal M}_n(k)\longrightarrow\bold{Top}$
\end{definition}

\begin{remark}
This construction and proof of Theorem \ref{main:theorem}, based on the analysis of the resulting colimits,
was inspired by the work of Clemens Berger on cellular operads. Our original proposed line of proof was to
form similar colimits over the barycentric subdivision of $M_n(k)$, associating to the barycenter the
intersection of the spaces $G(X)$ over all the vertices in the simplex.  This caused a great number of technical
difficulties due to the fact that some of these intersections are empty. 
\end{remark}

\begin{lemma}\label{F:contractible}
For any object $A\in{\cal M}_n(k)$ the inclusion
$$G(A)\subset F(A)$$
is a strong deformation retract.  Thus $F(A)$ is contractible.
\end{lemma}

\noindent
{\bf Proof.} The deformation retraction is constructed in a number of stages.  If $A=B\Box_i C$,
we first show that the subspace
$$\cup_{X_1\Box_i X_2\to A}G(X_1\Box_i X_2) \subset F(A)$$
is a strong deformation retract, where the union is taken over all objects of the form $X_1\Box_i X_2$
which map into $A$, with $X_1$, $X_2$ having the same underlying sets of generating objects as $B$ and
$C$ respectively.

Suppose $X$ is an arbitrary object of ${\cal M}_n(k)$ which maps into $A$.  Then define
$$X_1=X-|C|\qquad X_2=X-|B|,$$
(cf. Definition \ref{partial:objects}).  By the coherence theorem $X_1\Box_i X_2$ maps into $X$.
Now let $(d_1,d_2,\dots,d_k)$ be a configuration of little $n$-cubes contained in $G(X)$,
with
$$d_j=[u_{j1},v_{j1}]\times[u_{j2},v_{j2}]\times\dots\times[u_{jn},v_{jn}],$$
and let
$$M=\mbox{max}\{v_{bi}\ |\ b\mbox{ in } B\}\qquad m=\mbox{min}\{u_{ci}\ |\ c\mbox{ in } C\}.$$
If $M\le m$ then $(d_1,d_2,\dots,d_k)$ is contained in $G(X_1\Box_i X_2)$ and we leave the configuration
alone.  Otherwise if $M>m$, let $D_1$ denote the linear deformation which takes the closed interval
$[0, M]$ onto the closed interval $[0,\frac{M+m}{2}]$ and let $D_2$ denote the linear deformation which
takes the closed interval $[m,1]$ to the closed interval $[\frac{M+m}{2},1]$.  Now apply the deformation
$D_1$ (resp. $D_2$) simultaneously to the $i$-th coordinates of all little cubes $d_b$ (resp. $d_c$) whose labels
correspond to generators $b$ in $B$ (resp. $c$ in $C$).  We claim that this defines a strong deformation
retraction of $G(X)$ onto $G(X)\cap G(X_1\Box_i X_2)$.  The only nonobvious point is that the
retraction stays within $G(X)$.  This follows from the coherence theorem.  For the relative position of
any two little cubes in the configuration can change only if the label of one, say $d_b$ is in $B$ and the label of the
other $d_c$ is in $C$.  Moreover this only happens in the $i$-th coordinate direction and only if $d_b\not<_i d_c$.
So the only trouble which could arise is if $(d_1,d_2,\dots,d_k)\in G(X)$ required that $d_c<_i d_b$.  But
this could only happen if $c\Box_i b$ in $X$.  But if that were the case, by the coherence theorem, there
couldn't be a morphism $X\to A=B\Box_i C$ since $b\Box_i c$ in $A$.

By gluing together the deformations of $G(X)$ onto $G(X)\cap G(X_1\Box_i X_2)$ over all objects $X$
in ${\cal M}_n(k)$ mapping into $A$ one obtains that
$$\cup_{X_1\Box_i X_2\to A}G(X_1\Box_i X_2) \subset F(A)$$
is a strong deformation retract.  In the next stage of the deformation one decomposes $B=B'\Box_r B''$
and $C=C'\Box_s C''$ obtaining a decomposition
$$A=(B'\Box_r B'')\Box_i(C'\Box_s C''),$$
and then one shows by a similar argument that
$$\cup_{(X_1'\Box_r X_2'')\Box_i(X_2'\Box_s X_2'')}
G\left((X_1'\Box_r X_2'')\Box_i(X_2'\Box_s X_2'')\right)
\subset\cup_{X_1\Box_i X_2\to A}G(X_1\Box_i X_2)$$
is a strong deformation retract.  Composing the two retractions, one obtains that
$$\cup_{(X_1'\Box_r X_2'')\Box_i(X_2'\Box_s X_2'')}
G\left((X_1'\Box_r X_2'')\Box_i(X_2'\Box_s X_2'')\right)
\subset F(A)$$
is a strong deformation retract.  One continues this refinement process, restricting to objects $X$ which
map into $A$ and which resemble $A$ to an ever deeper level of parentheses and operations, showing at
each stage that the resulting union of $G(X)$ is a strong deformation retract of the union of $G(X)$ at
the preceding stage and hence is also a strong deformation retract of $F(A)$.  After finitely many stages
the only object $X$ left is $A$ itself.  Thus we obtain that $G(A)\subset F(A)$ is a strong deformation
retract.  Now as we noted in Definition \ref{G:definition}, $G(A)$ can be identified with a convex subspace of Euclidean
space and hence is contractible.  Therefore $F(A)$ is also contractible.

\begin{lemma}\label{intersection:lemma}
For any two objects $A$, $B$ of ${\cal M}_n(k)$,
$$F(A)\cap F(B)=\bigcup_{{\begin{array}{c}\scriptstyle X\to A\\ \scriptstyle X\to B\end{array} } }F(X),$$
where the union is indexed over all objects $X$ in ${\cal M}_n(k)$ which map into both $A$ and $B$.
\end{lemma}

\noindent
{\bf Proof.} First note that the inclusion
$$\bigcup_{{\begin{array}{c}\scriptstyle X\to A\\ \scriptstyle X\to B\end{array} } }F(X)\subseteq
F(A)\cap F(B)$$
is immediate from definition.  To prove equality, we proceed by double induction.  Our primary induction
is on $k$, the number of generating objects, starting with the observation that the lemma holds trivially
if $k=1$.  Building on this induction we first prove the following:

\noindent
{\it Claim.} If there are nontrivial decompositions $A=A_1\Box_i A_2$ and $B=B_1\Box_i B_2$
with $A_1$ and $B_1$ (and hence also $A_2$ and $B_2$) having the same underlying set of generating
objects, then the intersection $F(A)\cap F(B)$ satisfies the lemma.

By our primary induction:
$$
F(A_1)\cap F(B_1)=\bigcup_{{\begin{array}{c}\scriptstyle X_1\to A_1\\ \scriptstyle X_1\to B_1\end{array} } }F(X_1)
\qquad
F(A_2)\cap F(B_2)=\bigcup_{{\begin{array}{c}\scriptstyle X_2\to A_2\\ \scriptstyle X_2\to B_2\end{array} } }F(X_2)
$$
It follows immediately that
$$
F(A)\cap F(B)=\bigcup_{{\begin{array}{c}\scriptstyle X_1\to A_1\\ \scriptstyle X_1\to B_1\\
\scriptstyle X_2\to A_2\\ \scriptstyle X_2\to B_2\\\end{array} } }F(X_1\Box_i X_2)
$$
and thus implies the lemma in this case, proving the claim.

Our secondary induction is on objects $A$ and $B$ with respect to the ordering in the poset ${\cal M}_n(k)$.
If $A$ or $B$, say $A$, is minimal in the poset ${\cal M}_n(k)$, then $A$ has the form
$$A=j_1\Box_1 j_2\Box_1 j_3\Box_1\dots\Box_1 j_k.$$
Now there are two possibilities.  First of all if $m_1\Box_1 m_2$ in $B$ implies that $m_1\Box_1 m_2$
in $A$, then by the coherence theorem $A$ maps into $B$ and we have
$$F(A)\cap F(B)=F(A)$$
and we are done. Conversely if there is a pair of generating objects $m_1$, $m_2$ such that $m_1\Box_1 m_2$
in $A$, whereas $m_2\Box_1 m_1$ in $B$, then by the coherence theorem $m_2\Box_1 m_1$ in $X$ for
any object $X$ mapping into $B$.  Hence $G(A)\cap G(X)=\emptyset$ for all such $X$, since a requirement for a
configuration $(c_1,c_2,\dots,c_k)$ to lie in $G(A)$ is that $c_{m_1}<_1 c_{m_2}$ whereas a requirement for
that configuration to lie in $G(X)$ is that $c_{m_2}<_1 c_{m_1}$, and no configuration can simultaneously
satisfy both requirements. It follows that
$$F(A)\cap F(B)=\emptyset,$$
and the lemma again holds in this case.  This starts the secondary induction.

Now suppose we have shown that the lemma holds for all intersections $F(C)\cap F(D)$ where $C$ maps into
$A$, $D$ maps into $B$, and at least one of $C$, $D$ is not equal to $A$, resp. $B$.  Let us suppose that the
outermost operation in $A$ is $\Box_i$ and the outermost operation in $B$ is $\Box_j$.  Thus
$A=A_1\Box_i A_2$ and $B=B_1 \Box_j B_2$.  Without loss of generality, we may suppose that $i\le j$.

Consider the partial objects
$$A'_1=B\cap |A_1|\qquad A'_2=B\cap |A_2|$$
(cf. Definition \ref{partial:objects}).  We clearly have
$$F(A)\cap F(B) = F(A)\cap F(A'_1\Box_i A'_2)\cap F(B)$$
We can apply the claim above to the intersection $F(A)\cap F(A'_1\Box_i A'_2)$.  We can then distribute
the intersection with $F(B)$ over the resulting union.  If $A= A'_1\Box_i A'_2$, we get no reduction, since
then $F(A)\cap F(A'_1\Box_i A'_2)=F(A)$.  Otherwise we can apply our secondary induction to the resulting
union of intersections.  Similarly we consider
$$B'_1=A\cap |B_1| \qquad B'_2=A\cap |B_2|,$$
note that
$$F(A)\cap F(B) = F(A)\cap F(B'_1\Box_j B'_2)\cap F(B)$$
and apply the claim to the intersection $F(B'_1\Box_j B'_2)\cap F(B)$.  Again using our secondary induction
we obtain that the lemma applies unless $B=B'_1\Box_j B'_2$.  

Thus we are left with the case when both $A= A'_1\Box_i A'_2$ and $B=B'_1\Box_j B'_2$
hold.  But in this case we must have decompositions
$$A=(C\Box_j D)\Box_i (U\Box_j V)\qquad B=(C\Box_i U)\Box_j (D\Box_i V)$$
for some objects $C$, $D$, $U$ and $V$.  Now if $i<j$, then there is a morphism 
$\eta^{ij}_{C,D,U,V}:A\longrightarrow B$.  Hence $F(A)\cap F(B)=F(A)$ and the lemma holds.  If $i=j$
and either $D=0$ or $U=0$, then $A=B$ and again we are done.  Finally if both $D\ne0$ and $U\ne0$,
then $G(A)\cap G(B)=\emptyset$ for a configuration of little $n$-cubes in the intersection would have
to satisfy contradictory specifications on the relative positions of little $n$-cubes with labels in $D$ and
$U$.  This then means that
$$
F(A)\cap F(B) =
\bigcup_{{\begin{array}{c}\scriptstyle X\to A\\ \scriptstyle X\ne A\end{array} } }F(X)\cap F(B)
\ \bigcup \ 
\bigcup_{{\begin{array}{c}\scriptstyle Y\to B\\ \scriptstyle Y\ne B\end{array} } }F(A)\cap F(Y)
$$
and we can apply our secondary induction.  This concludes the induction and proof.

\bigskip
In view of Remark \ref{decomposableCubes:remark} and
the obvious fact that the inclusion of a finite union of closed convex spaces of $\Bbb R^N$ into a bigger
such finite union is a closed cofibration, a direct consequence of the preceding lemma is:

\begin{corollary}\label{colim:corollary}
The natural map induced by inclusions
$$\mbox{colim}_{A\in{\cal M}_n(k)} F(A)\longrightarrow
\cup_{A\in Obj({\cal M}_n(k))} F(A) ={\cal D}_n(k)$$
is a homeomorphism.  Moreover for each object $A$ in ${\cal M}_n(k)$ the induced map
$$\mbox{colim}_{X\stackrel{\ne}{\longrightarrow} A} F(X)\longrightarrow F(A)$$
is a closed cofibration.
\end{corollary}

The main technical ingredient in the proof of Theorem \ref{main:theorem} is the following:

\begin{proposition}\label{hocolim:proposition}
Let ${\cal P}$ be a finite poset and let
$F:{\cal P}\longrightarrow \bold{Top}$ be a functor satisfying the property that
for each object $i$ in ${\cal P}$ the induced map
$$\mbox{colim}_{j<i}F(j)\longrightarrow F(i)$$
is a closed cofibration.
Then the natural map $\mbox{hocolim}_{{\cal P}} F\longrightarrow\mbox{colim}_{{\cal P}} F$ is
an equivalence.
\end{proposition}

\noindent
{\bf Proof:} We first observe that
$$\mbox{hocolim}_{{\cal P}} F=\mbox{colim}_{{\cal P}} G$$
where $G:{\cal P}\longrightarrow \bold{Top}$ is given by
$$G(i)=\mbox{hocolim}_{j\le i} F$$
and that $G$ satisfies the cofibration condition also. We note that $G(i)\longrightarrow F(i)$ is an equivalence for all
objects $i\in{\cal P}$. 

Then we filter the objects of ${\cal P}$ according to the length
of the largest increasing chain of objects which terminates in the given object.  Thus the objects of
filtration 0 are the minimal objects.  We denote by ${\cal P}_k$ the full subcategory of ${\cal P}$
whose objects have filtration $\le k$. We proceed by induction on $k$ to show that
$$\mbox{colim}_{{\cal Q}_k} G\longrightarrow\mbox{colim}_{{\cal Q}_k} F$$
is an equivalence, for any subposet ${\cal Q}_k\subseteq {\cal P}_k$ satisfying the condition that if $j<i$
and $i\in {\cal Q}_k$, then $j\in {\cal Q}_k$. This is true for $k=0$ since in that case the colimits are just
disjoint unions of the values of $G$ and $F$ over minimal objects.

The induction step from $k-1$ to $k$ is based upon the pushout lemma for equivalences:
suppose given a commutative cubical diagram of spaces and maps as shown
$$\diagram &\bullet
\rrto|<\hole|<<\ahook
\dline
\dlto_\alpha
&&\bullet \dlto_\beta \ddto \\
\bullet\rrto|<\hole|<<\ahook
\ddto
&\hole\dto&\bullet \ddto \\
&\bullet\dlto_\gamma \rline|<\hole|<<\ahook
 &\rto &\bullet \dlto_\delta \\
\bullet\rrto|<\hole|<<\ahook&&\bullet
\enddiagram$$
Assume that the front and back faces are pushout squares with the map across the top being a closed cofibration
in each case. (In the sequel we will refer to such pushout squares as cofibration squares. It will also be
useful to note that in such a cofibration square the map across the bottom is also a cofibration.) If the
maps marked $\alpha$, $\beta$, and $\gamma$ are equivalences, then so is the map marked $\delta$.

We note that we have a pushout square
$$\diagram \coprod_{i\in{\cal Q}_k \mbox{\scriptsize\& filt}(i)=k}\mbox{colim}_{j<i}F(j)\rto|<<\ahook\dto
&\coprod_{i\in{\cal Q}_k \mbox{\scriptsize\& filt}(i)=k}F(i) \dto\\
\mbox{colim}_{{\cal Q}_{k-1}}F\rto|<<\ahook&\mbox{colim}_{{\cal Q}_k}F
\enddiagram$$
with ${\cal Q}_{k-1}= {\cal Q}_k\cap {\cal P}_{k-1}$.
The map across the top is a closed cofibration by hypothesis.

The same considerations apply to the functor $G$ and we get an analogous cofibration square.
We thus obtain a commutative cube as in the pushout lemma, with the front face being the
cofibration square for $F$ and the back face being the cofibration square for $G$, and the
maps from the back face to the front face being induced by the natural transformation $G\to F$.
It is immediate that the map corresponding to $\beta$ is an equivalence, while the maps corresponding
to $\gamma$ and $\beta$ are equivalences by the induction hypothesis.  This completes the induction
and proof.

\begin{remark} Proposition \ref{hocolim:proposition} is true for any cofinite strongly directed set ${\cal P}$ (ie. ${\cal P}$
is a directed set such that $a\le b$ and $b\le a$ implies $a=b$, and each $a\in{\cal P}$ has only a
finite number of predecessors).  This statement is a fairly immediate consequence of the closed model
category structure on the category of ${\cal P}$-diagrams in $\bold{Top}$ dual to the one constructed by
Edwards and Hastings in \cite[\S (3.2)]{EH}
\end{remark}

\begin{lemma}\label{hocolim:operad}
Let $\{{\cal M}(n)\}_{n\ge0}$ be an operad in the category of small categories.  Let
$\{F_n:{\cal M}(n)\longrightarrow \bold{Top}\}_{n\ge0}$ be a collection of functors satisfying the following
conditions:
\begin{enumerate}
\item There is an operad ${\cal C}$ such that for each object $A$ of ${\cal M}(n)$ $F_n(A)\subseteq{\cal C}(n)$,
and for each morphism $f:A\to B$ in ${\cal M}(n)$ $F_n(f):F_n(A)\longrightarrow F_n(B)$ is an inclusion.
\item For each permutation $\sigma\in\Sigma_n$, action by $\sigma$ on ${\cal C}(n)$ sends the subspace
$F_n(A)$ to the subspace $F_n(A\sigma)$.
\item Given objects $A\in{\cal M}(n)$, $B_i\in{\cal M}(j_i)$ $1\le i\le n$, the structure map
$${\cal C}(n)\times{\cal C}(j_1)\times{\cal C}(j_2)\times\dots\times{\cal C}(j_n)\longrightarrow
{\cal C}(j_1+j_2+\dots+j_n)$$
sends the subspace $F_n(A)\times F_{j_1}(B_1)\times\dots\times F_{j_n}(B_n)$ to the subspace
$F_{j_1+\dots+j_n}(\gamma(A; B_1,B_2,\dots,B_n))$, where $\gamma$ denotes the structure map
of ${\cal M}$.
\item The unit element in ${\cal C}(1)$ is contained in $F_1(1)$, where 1 denotes the unit element of ${\cal M}(1)$.
\end{enumerate}
Then $\{\mbox{hocolim}_{{\cal M}(n)}F_n\}_{n\ge0}$ is an operad and the natural map
${\{\mbox{hocolim}_{{\cal M}(n)}F_n\}_{n\ge0}\longrightarrow{\cal C}}$ is a map of operads.
\end{lemma}
\vspace{0.3cm}

The proof of this lemma is completely straightforward and will be left as an exercise for the reader. Moreover
we also note that in case the action of $\Sigma_n$ on both ${\cal M}(n)$ and ${\cal C}(n)$ is free, then the same
is the case with the action on $\mbox{hocolim}_{{\cal M}(n)}F_n$.

\bigskip

\noindent
{\bf Proof of Theorem \ref{main:theorem}.} Corollary \ref{colim:corollary}, Proposition \ref{hocolim:proposition},
Lemma \ref{hocolim:operad} and Remark \ref{decomposableCubes:remark} imply that the chain
$$\left\{\mbox{hocolim}_{{\cal M}_n(k)} F\longrightarrow\mbox{colim}_{{\cal M}_n(k)}F\cong {\cal D}_n(k)\subset
{\cal C}_n(k)\right\}_{k\ge 0}$$
is a chain of operad maps which are also equivalences. Similarly by Lemma \ref{F:contractible} the natural map of the diagram 
$F:{\cal M}_n(k)\longrightarrow\bold{Top}$ to the trivial diagram $*:{\cal M}_n(k)\longrightarrow\bold{Top}$
induces a map of operads which is also an equivalence:
$$\left\{\mbox{hocolim}_{{\cal M}_n(k)} F\longrightarrow \mbox{hocolim}_{{\cal M}_n(k)} * ={\cal NM}_n(k)
={\cal M}_n(k)\right\}_{k\ge 0}$$
where the last equality is our usual notational abuse of using the same symbol for a category and its nerve.

It remains to show that  the inclusion of the Milgram preoperad $\overline{\cal J}_n(k)$ in the operad ${\cal M}_n(k)$
is an equivalence.
To do this requires defining a subdiagram of subspaces of the diagram $F$, indexed over the Milgram subcategory
$\overline{\cal J}_n(k)$.  Specifically given an object $A$ in $\overline{\cal J}_n(k)$ we define
$$\overline{F}(A)=\cup_{X\to A} G(X)$$
where the union is indexed over all objects $X$ in $\overline{\cal J}_n(k)$ ({\em not} ${\cal M}_n(k)$) mapping into
$A$.  The inclusion of diagrams then induces a commutative diagram:
$$\diagram
\overline{\cal J}_n(k)\dto &\mbox{hocolim}_{\overline{\cal J}_n(k)} \overline{F}\lto\rto|<<\ahook\dto
&\mbox{colim}_{\overline{\cal J}_n(k)}\overline{F}\dto\\
{\cal M}_n(k)&\mbox{hocolim}_{{\cal M}_n(k)} F\lto_(.55){\simeq}\rto|<<\ahook^(.6){\simeq}
&{\cal C}_n(k)
\enddiagram$$
We have already shown that the maps across the bottom row are equivalences.  Using similar arguments,
first proving the analogs of Lemmas \ref{F:contractible} and \ref{intersection:lemma} and
Corollary \ref{colim:corollary} hold for the diagram $\overline{F}$,
we can show the maps across the top row are also equivalences.  Thus it suffices to show that the right
hand vertical arrow is an equivalence.

By the analog of Corollary \ref{colim:corollary} we can identify
$\mbox{colim}_{\overline{\cal J}_n(k)}\overline{F}$ with the union
$$\bigcup_{A\in\overline{\cal J}_n(k)}G(A)\subset {\cal D}_n(k)\subset {\cal C}_n(k)$$
This in turn is the subspace of Milgram decomposable configurations of little $n$-cubes.
A configuration in ${\cal C}_n(k)$ is said to be {\it Milgram decomposable\/} if one can cut
through the configuration with a finite set of hyperplanes perpendicular to the first coordinate axis which miss
the interiors of all the little $n$-cubes and each of the resulting strips individually can then be cut through
by a finite number of hyperplanes perpendicular to the second coordinate axis (again missing
the interiors of all the little cubes in the strip), and each of those resulting strips can then be cut by
hyperplanes perpendicular to the the third coordinate axis, etc. with the final cuts being done by hyperplanes
perpendicular to the last coordinate axis, so that at the end of this process there is exactly one little cube
in each compartment.

The following two figures in ${\cal C}_2(k)$ illustrate the concept of Milgram decomposability
$$\mbox{
\centerline{\hbox{
\unitlength1cm
\begin{picture}(12, 4.5)
\put(0, 0) {\line(1, 0) {4}}
\put(0, 0) {\line(0, 1) {4}}
\put(0, 4) {\line(1, 0) {4}}
\put(4, 0) {\line(0, 1) {4}}
\put(1, 0) {\line(0, 1) {4}}
\put(3, 0) {\line(0, 1) {4}}
\put(0, 1) {\line(1, 0) {1}}
\put(0, 3) {\line(1, 0) {1}}
\put(1, 2) {\line(1, 0) {2}}
\put(3, 3) {\line(1, 0) {1}}
\put(0.4, 0.2) {3}
\put(0.4, 1.7) {6}
\put(0.4, 3.2) {1}
\put(1.9, 0.7) {2}
\put(3.4, 1.2) {4}
\put(3.4, 3.2) {5}
\put(8, 0) {\line(1, 0) {4}}
\put(8, 0) {\line(0, 1) {4}}
\put(8, 4) {\line(1, 0) {4}}
\put(12, 0) {\line(0, 1) {4}}
\put(8, 2) {\line(1, 0) {4}}
\put(10, 2) {\line(0, 1) {2}}
\put(9.9, 0.7) {2}
\put(8.9, 2.7) {1}
\put(10.9,2.7) {3}
\end{picture}
}}
}$$
The configuration on the left is Milgram decomposable, whereas the one on the right is not (although it is
decomposable).

We now show that that the inclusion of the space of Milgram decomposable configurations of little $n$-cubes
into ${\cal C}_n(k)$ is an equivalence.  Given any configuration of little $n$-cubes in ${\cal C}_n(k)$
let $m$ be the minimum distance between barycenters of different subcubes in the $\ell_\infty$ norm.
Define a map ${\cal C}_n(k)\to {\cal C}_n(k)$ which linearly shrinks (towards their barycenters) those
the little cubes in a configuration whose dimensions are bigger than $\frac{m}{2k}$ by $\frac{m}{2k}$
to subcubes of this size (leaving alone dimensions of cubes which are smaller).  This map is clearly
homotopic to the identity. 

It also takes any configuration to a Milgram decomposable one by the following argument.
We say that two little cubes in a configuration overlap in the first coordinate direction if there is a hyperplane
perpendicular to the first coordinate direction which passes through the interiors of both.  We say that two little cubes
are in the same 1-clump if there is chain of little cubes from one to the other such that any two adjacent
ones in the chain overlap in the first coordinate direction. Clearly the  1-clumps of little cubes can be
separated from each other by hyperplanes perpendicular to the first coordinate direction.  The barycenters of
any two little cubes in the same 1-clump are separated in the first coordinate direction by a distance at
most $\frac{m}{2}$. (There at most $k$ elements in the chain connecting the little cubes, with the
barycenters of adjacent subcubes in the chain having separation in the first coordinate direction 
at most $\frac{m}{2k}$.)  Thus the separation in at least one of the other coordinate directions between
the barycenters of any two little cubes in the 1-clump must be at least $m$.  Now for the little cubes
within a given 1-clump define an analogous notion of 2-clump and repeat.  At the final stage of this process
we will have an $(n-1)$-clump of cubes which overlap in all the coordinate directions except the last.
It will follow that all the little cubes in this $(n-1)$-clump must have barycenters separated in the last
coordinate direction by distances of at least $m$. Since the little cubes have dimensions $\frac{m}{2k}$,
they can then be separated from one another by hyperplanes perpendicular to the last coordinate direction,
proving the configuration is Milgram decomposable.

Moreover the homotopy from the shrinking map to the identity restricts to the subspace of Milgram
decomposable configurations.  It follows that the inclusion of
the space of Milgram decomposable configurations in ${\cal C}_n(k)$ 
is an equivalence, completing the proof of Theorem \ref{main:theorem}.

\bigskip

Before we proceed to the proof of Theorem \ref{graphs:theorem} we recall the basic definitions, due to
Berger \cite{Be2}.

\begin{definition}
An {\it acyclic orientation\/} of the complete graph on the set of vertices $\{1,2,3,\dots,k\}$ is an
assignment of direction to each edge of the graph such that no directed cycles occur.  Equivalently
an acyclic orientation is determined uniquely by a total ordering (ie. a permutation) of the vertices.
A {\it coloring\/} of the complete graph on $k$ vertices is an assignment of colors to each edge of the
graph from the countable set of colors $\{1,2,3,\dots\}$.  The poset ${\cal K}(k)$ has as elements pairs
$(\mu,\sigma)$, where $\mu$ is a coloring and $\sigma$ is an acyclic orientation of the complete graph
on $k$ vertices.  The order relation on ${\cal K}(k)$ is determined as follows: we say that
$(\mu_1,\sigma_1)\le (\mu_2,\sigma_2)$ if for every edge $a\stackrel{i}{\longrightarrow}b$ in
$(\mu_1,\sigma_1)$ the corresponding edge in $(\mu_2,\sigma_2)$ has either orientation and coloring
 $a\stackrel{j}{\longrightarrow}b$ with $j\ge i$ or $b\stackrel{j}{\longrightarrow}a$ with $j>i$.
 Per our usual abuse we also denote by ${\cal K}(k)$ the nerve of this poset.

The action of the symmetric group $\Sigma_k$ on ${\cal K}(k)$ is via permutation of the vertices.
The structure map
$${\cal K}(k)\times{\cal K}(m_1)\times{\cal K}(m_2)\times\dots\times{\cal K}(m_k)\longrightarrow
{\cal K}_(m_1+m_2+\dots+m_k)$$
assigns to a tuple of orientations and colorings in
${\cal K}(k)\times{\cal K}(m_1)\times{\cal K}(m_2)\times\dots\times{\cal K}(m_k)$ the orientation
and coloring obtained by subdividing the set of $m_1+m_2+\dots+m_k$ vertices into $k$ blocks containing
$m_1$, $m_2$, \dots, $m_k$ vertices respectively. The edges connecting vertices within the $i$-th block
are oriented and colored according to the given element in ${\cal K}(m_i)$.  The edges connecting vertices
between blocks $i$ and $j$ are all oriented and colored according to the corresponding edge in the given element
of ${\cal K}(k)$.  It is easy to check that this specification gives ${\cal K}(k)_{k\ge0}$ the structure of an
$E_\infty$ operad. 

The $n$-th filtration ${\cal K}^{(n)}(k)$ is the subposet of ${\cal K}(k)$ where the colorings are restricted
to take values in the subset $\{1,2,3,\dots,n\}$. It is obvious that ${\cal K}^{(n)}(k)_{k\ge0}$ is a
suboperad of ${\cal K}$.  There is an inclusion of posets
${\cal M}_n(k)\subset {\cal K}^{(n)}(k)$ which takes an object $A$ to the complete graph  on
$\{1,2,3,\dots,k\}$ with edges oriented and colored as follows:
$$a\stackrel{i}{\longrightarrow}b\qquad\mbox{if }a\Box_i b\mbox{ in }A$$
These inclusions define a map of operads ${\cal M}_n\to{\cal K}^{(n)}$.
\end{definition}

\begin{remark}
This definition departs slightly from that in \cite{Be2} in that Berger takes colorings with values in the
set $\{0,1,2,\dots\}$ rather than $\{1,2,3,\dots\}$.
\end{remark}

\begin{definition}
Define $\Gamma(k)$ to be the category whose objects are permutations in $\Sigma_k$ and with a
unique morphism between any two objects (which is hence an isomorphism).  The nerve of this category,
also abusively denoted $\Gamma(k)$, can be identified with the standard simplicial model $E\Sigma_k$ of
the total space of the universal principal $\Sigma_k$ bundle.  By rewriting the objects of $\Gamma(k)$
in the form
$$i_1\Box i_2\Box\dots\Box i_k,$$
instead of $[i_1,i_2,\dots,i_k]$, and appealing to MacLane's coherence theorem, we can identify $\Gamma(k)$
with a full subcategory of the free strict symmetric monoidal category on $k$ generators.  Thus the
spaces $\{\Gamma(k)\}_{k\ge0}$ can be naturally endowed with the structure of an operad which acts on the
nerves of strict symmetric monoidal categories.  The operad $\Gamma$ was extensively studied by Barratt
and Eccles \cite{BEc} and May \cite{May2} (who denotes the operad ${\cal D}$ instead).

Smith \cite{Sm} defined a filtration on $\Gamma$ as follows.  First of all he defined $\Gamma^{(n)}(2)$ to
be the $n-1$ skeleton of $\Gamma(2)$, which is easily identified as the standard $\Bbb Z/2$-equivariant
simplicial model of $S^{n-1}$. Then he defined a simplex in $\Gamma(k)$ to be in the $n$-th filtration
$\Gamma^{(n)}(k)$ if its images under all restriction maps
$$R_{a,b}:\Gamma(k)\longrightarrow\Gamma(2)$$
lies in $\Gamma^{(n)}(2)$ (cf. Remark \ref{octo:coherence}).  Equivalently an $r$-simplex
$$\sigma_0\longrightarrow\sigma_1\longrightarrow
\sigma_2\longrightarrow\dots\longrightarrow\sigma_r$$
lies in $\Gamma^{(n)}(k)$ if any pair of elements $a$, $b$ in $\{1,2,\dots, k\}$ change their
relative order in the given sequence of permutations at most $n-1$ times.  For example the 3-simplex
$$[1,2,3]\longrightarrow[2,1,3]\longrightarrow[2,3,1]\longrightarrow[2,1,3]$$
lies in $\Gamma^{(3)}(3)$ since the pair $(1,2)$ changes order once, the pair $(1,3)$ changes order
twice and the pair $(2,3)$ doesn't change order at all.  It is easy to see that $\Gamma^{(n)}(k)_{k\ge0}$
forms a suboperad of $\Gamma$.

The forgetful map $(\mu,\sigma)\mapsto\sigma$, which forgets the coloring, defines a functor and hence
a map of operads ${\cal K}\to\Gamma$.  It also preserves filtrations.  For given an $r$-simplex
$$(\mu_0,\sigma_0)\longrightarrow(\mu_1,\sigma_1)\longrightarrow(\mu_2,\sigma_2)\longrightarrow
\dots\longrightarrow(\mu_r,\sigma_r)$$
in ${\cal K}^{(n)}(k)$, any edge connecting two given vertices $a$ and $b$ can only change direction at
most $n-1$ times.  For every change in direction must correspond to an incrementation of the coloring of
that edge.

The composite
$${\cal M}_n\longrightarrow{\cal K}^{(n)}\longrightarrow\Gamma^{(n)}$$
can be identified with the map of operads arising from the fact that any symmetric monoidal category is
$n$-fold monoidal (cf. Remark \ref{symm:remark}).
\end{definition}

Smith \cite{Sm} conjectured that $\Gamma^{(n)}$ has the same homotopy type as the little $n$-cubes
operad ${\cal C}_n$, and thus could also be used to parametrize the structure of an $n$-fold loop space.
This conjecture was proved by Berger \cite{Be2}.  Our proof of Theorem \ref{graphs:theorem} below
gives an alternative proof of this conjecture.

\bigskip
\noindent
{\bf Proof Sketch of Theorem \ref{graphs:theorem}} 
 The diagram $F$ in the proof of Theorem \ref{main:theorem} can be expanded
in the evident way to a diagram on ${\cal K}^{(n)}(k)$ containing $F$ as a subdiagram of subspaces, and the
inclusion ${\cal M}_n\to{\cal K}^{(n)}$ can be shown to be an equivalence by an argument similar to the
proof we used above to prove that $\overline{\cal J}_n\subset{\cal M}_n$ is an equivalence. 
See \cite{Be2} for details.

To show that the map $p: {\cal K}^{(n)}\longrightarrow\Gamma^{(n)}$ is an equivalence we have to show that
for any simplex $\mathfrak S$ in $\Gamma^{(n)}$ the inverse image $p^{-1}({\mathfrak S})$ is contractible.
We prove this by induction on the dimension of $\mathfrak S$. If ${\mathfrak S}=\sigma$ is a vertex, then
$p^{-1}(\sigma)$ is a simplicial cone on the object $(\mu_0,\sigma)$, where $\mu_0$ is the coloring which
assigns to each edge the color 1.

Assume we have already shown the contractibility of inverse images for simplices of dimension less than that of
 $$\mathfrak S\quad=\quad\sigma_0\longrightarrow\sigma_1\longrightarrow\sigma_2\longrightarrow
 \dots\longrightarrow\sigma_r$$
We note that
$$p^{-1}({\mathfrak S})=T({\mathfrak S})\cup\bigcup_{i=0}^r p^{-1}({\mathfrak S}_i),$$
where $T({\mathfrak S})$ is the union of all simplices in ${\cal K}^{(n)}$ which map surjectively onto
$\mathfrak S$ and the ${\mathfrak S}_i$ are the codimension 1 faces of $\mathfrak S$.  To show that
this union is contractible, it suffices to show that all the intersections
\begin{eqnarray*}
\bigcap_{j\in J}p^{-1}({\mathfrak S}_j) &= &p^{-1}\left(\bigcap_{j\in J}{\mathfrak S}_j\right)\\
T({\mathfrak S})\cap\bigcap_{j\in J}p^{-1}({\mathfrak S}_j) &= 
&T({\mathfrak S})\cap p^{-1}\left(\bigcap_{j\in J}{\mathfrak S}_j\right)
\end{eqnarray*}
are contratible. Intersections of the first kind are contractible by induction hypothesis.  To see that intersections
of the second kind are contractible, we first consider the following distinguished simplex in $T({\mathfrak S})$:
$$(\mu_0,\sigma_0)\longrightarrow(\mu_1,\sigma_1)\longrightarrow(\mu_2,\sigma_2)\longrightarrow
\dots\longrightarrow(\mu_r,\sigma_r)$$
where the coloring $\mu_i$ assigns to the edges joining a pair of vertices $a$, $b$ the color which is 1 more
than the number of times that this pair of elements changes relative position in the subsimplex
$$\sigma_0\longrightarrow\sigma_1\longrightarrow\sigma_2\longrightarrow
 \dots\longrightarrow\sigma_i$$
of $\mathfrak S$.  Then it is easy to see that
$T({\mathfrak S})\cap p^{-1}\left(\bigcap_{j\in J}{\mathfrak S}_j\right)$
is a cone on the vertex $(\mu_m,\sigma_m)$ where $\sigma_m$ is the initial vertex of the face
$\bigcap_{j\in J}{\mathfrak S}_j$ of $\mathfrak S$, and is thus contractible.  This completes the
induction and proof.

\end{document}